\documentclass[preprint,sort&compress]{elsarticle}
\usepackage[utf8]{inputenc}

\usepackage{orcidlink}
\usepackage{amsmath,amssymb,amsthm}
\usepackage{bm}
\usepackage{multirow}
\usepackage[algoruled,boxed,lined]{algorithm2e}
\usepackage{algpseudocode}

\usepackage{hyperref}
\usepackage[misc]{ifsym} 

\usepackage{color}
\usepackage{graphicx}
\graphicspath{{./Figures/}}

\usepackage{amsfonts}
\usepackage{mathtools}

\usepackage[mathscr]{eucal}

\usepackage{caption}
\usepackage{subcaption}

\usepackage{booktabs}
\usepackage{array}

\newtheorem{theorem}{Theorem}
\newtheorem{corollary}{Corollary}
\newtheorem{remark}{Remark}
\newtheorem{assumption}{Assumption}

\providecommand{\rset}{\mathbb{R}} 

\providecommand{\E}[1]{{\ensuremath{\mathbb{E}}\mspace{-2mu}\left[#1\right]}}    
\providecommand{\var}[1]{{\ensuremath{\mathrm{Var}}\mspace{-2mu}\left[#1\right]}}

\providecommand{\nrs}{\ensuremath{M}}                         
\providecommand{\nelem}{\ensuremath{N}}          
\providecommand{\EstMLMC}{\mathcal{A}_{MLMC}}    

\providecommand{\QoI}{\ensuremath{Q}} 
\providecommand{\tol}{\ensuremath{\mathrm{TOL}}} 
\providecommand{\work}{\ensuremath{W}}           
\providecommand{\Ow}{{q_{w}}}                    
\providecommand{\Os}{{q_{s}}}                    
\providecommand{\splitting}{\theta}              
\providecommand{\confpar}{C_\xi}                 
\providecommand{\tolbsep}{{C}}                   

\newcolumntype{M}[1]{>{\centering\arraybackslash}m{#1}} 

\newcommand*\xoverline[2][0.75]{%
	\sbox{\myboxA}{$\m@th#2$}%
	\setbox\myboxB\null
	\ht\myboxB=\ht\myboxA%
	\dp\myboxB=\dp\myboxA%
	\wd\myboxB=#1\wd\myboxA
	\sbox\myboxB{$\m@th\overline{\copy\myboxB}$}
	\setlength\mylenA{\the\wd\myboxA}
	\addtolength\mylenA{-\the\wd\myboxB}%
	\ifdim\wd\myboxB<\wd\myboxA%
	\rlap{\hskip 0.5\mylenA\usebox\myboxB}{\usebox\myboxA}%
	\else
	\hskip -0.5\mylenA\rlap{\usebox\myboxA}{\hskip 0.5\mylenA\usebox\myboxB}%
	\fi}	

\begin{document}

\title{Goal-Oriented Adaptive Finite Element\\ Multilevel {M}onte {C}arlo with Convergence Rates}

\author[CEMSE]{Joakim~Beck\fnref{fn1}\,\orcidlink{0000-0001-8042-998X}}
\author[CEMSE]{Yang~Liu\corref{cor1}\fnref{fn1}\,\orcidlink{0000-0003-0778-3872}}
\ead{yang.liu.3@kaust.edu.sa}
\author[CEMSE]{Erik~von~Schwerin\fnref{fn1}\,\orcidlink{0000-0002-2964-7225}}
\author[CEMSE]{Ra\'{u}l~Tempone\fnref{fn1,fn2}\,\orcidlink{0000-0003-1967-4446}}
\address[CEMSE]{Computer, Electrical and Mathematical Sciences and Engineering,\\
  4700 King Abdullah University of Science and Technology (KAUST),\\
  Thuwal 23955-6900, Kingdom of Saudi Arabia}

\cortext[cor1]{Corresponding author}
\fntext[fn1]{KAUST SRI Center for Uncertainty Quantification in Computational Science and Engineering}
\fntext[fn2]{Alexander von Humboldt Professor in Mathematics for Uncertainty
  Quantification, RWTH Aachen University, 52062 Aachen, Germany.}

\markboth{Beck et al.}{Goal-Oriented Adaptive Multilevel {M}onte {C}arlo}
\begin{abstract}

{In this study, we} present an {adaptive multilevel} Monte Carlo (AMLMC) algorithm
for approximat{ing} deterministic, real-valued, bounded linear functionals {that depend} on the solution of a linear elliptic PDE with {a} lognormal diffusivity coefficient
 and geometric singularities in bounded domains of ${\mathbb{R}}^d$. 
 {Our AMLMC algorithm is built on the results of the} weak convergence rates 
 {in the work [Moon et al., BIT Numer. Math., 46 (2006), pp.~367--407]} for an adaptive algorithm {using} isoparametric d-linear quadrilateral finite element approximations and the dual
weighted residual error representation in {a} deterministic setting.
{Designed to suit} the geometric nature of the singularities in the solution, our AMLMC algorithm uses a sequence of deterministic, non-uniform auxiliary meshes as a building block.
The above-mentioned deterministic adaptive algorithm generates these meshes, corresponding to a geometrically decreasing sequence of tolerances. {In particular}, for a given realization of the diffusivity coefficient and accuracy level, AMLMC constructs its 
approximate sample using the first mesh in the hierarchy that satisfies {the} corresponding bias accuracy constraint.
This adaptive approach is particularly useful for the lognormal case treated here, which lacks uniform coercivity and thus produces functional outputs that vary over orders of magnitude when sampled.
Furthermore, we discuss iterative solvers and compare their efficiency with direct ones. To {reduce} computational work, we propose a stopping criterion for the iterative solver with respect to the quantity of interest, the realization of the diffusivity coefficient, and the desired level of AMLMC approximation.

From the numerical experiments, based on a Fourier expansion of the diffusivity coefficient field, we observe {improvements in efficiency} compared {with} both standard {Monte Carlo (MC)} and standard MLMC (SMLMC) {for a problem with} a singularity similar to that at the tip of a slit modeling a crack. 

\end{abstract}

\begin{keyword} 
  Multilevel Monte Carlo 
  \sep Goal-oriented adaptivity 
  \sep Computational complexity 
  \sep Finite elements
  \sep Partial differential equations with random data 
  \sep Lognormal diffusion
  \MSC[2020] 65C05 
  \sep 65N50 
  \sep 65N22 
  \sep 35R60 
\end{keyword} 
\maketitle

\section{Introduction}
\label{sec:intro}

Motivated by the importance of uncertainty quantification on {random partial differential equations} (RPDEs) \cite{Cetal2017,OBF2017,lord_powell_shardlow_2014}, we consider the adaptive computation and error control for {quantities of interest} (QoIs) of the form $\E{\QoI(u)},$ where $\QoI$ is a deterministic, real{-}valued, bounded linear functional of the stochastic solution $u$ to a class of linear elliptic PDE with random coefficients,
\begin{subequations}
	\label{eq:bvp_general}
	\begin{align}
		\label{eq:pde_general}
  	-\nabla \cdot \left( a(\mathbf x; \omega) \nabla u(\mathbf x; \omega) \right) 
  	&= f(\mathbf x; \omega) &&\text{for $\mathbf x  \in \mathcal D$,} 
    \\
  	u(\mathbf x; \omega) &=  0 &&\text{for $\mathbf x \in \partial\mathcal{D}_1$},
  	\\
    \partial_n u(\mathbf x; \omega) &=  0 &&\text{for $\mathbf x \in \partial\mathcal{D} -\partial\mathcal{D}_1$},
	\end{align}
\end{subequations}
where the differential operators, $\nabla\cdot$ and $\nabla$, are taken with respect to the spatial variable, $\mathbf x$, and  $\partial_n$ is the outward normal derivative operator. Here{,} $\omega$ corresponds to 
a complete probability space $(\Omega,\mathcal{F},\mathbb P)$. With respect to the spatial domain $\mathcal{D}$, {suppose it is} an open and bounded 
polygonal/polyhedral domain in $\rset^d$, $d\ge2$ with boundary $\partial\mathcal{D}$.
The boundary $\mathcal{D}$ splits into its Dirichlet and Neumann parts, and we {consider each part as} the union of a finite number of intervals or polygons.
{Naturally, the solution $u$ is stochastic {because of} the randomness induced by the stochastic diffusivity coefficient field $a(\mathbf x;\omega)$ and/or the stochastic forcing $f(\mathbf x; \omega)$.} {Using} properly-refined discretizations by isoparametric d-linear quadrilateral finite elements, {we propose the development and analysis of an adaptive multilevel} Monte Carlo (AMLMC) algorithm to approximate $\E{\QoI(u)}$ with error control.

An error expansion with a computable leading order term that is asymptotically accurate guides our AMLMC stochastic refinements. {First developed in~\cite{adFEM_our}, this error expansion uses} a dual weighted error representation and yields an error density.

{The error density characterizes} the optimality of refined meshes and the resulting complexity of AMLMC. {It also provides} insight {into} the comparison with MLMC {using} uniform discretizations, {thus} motivating the class of problems where AMLMC provides a computational advantage. A {specific} challenge to address in this work is {the lognormal distribution} in the diffusivity coefficient, $a$, {which implies that~\eqref{eq:bvp_general} lacks} uniform coercivity, {cf.~\cite{BNT_SIAM_review}.}
{Thus, it maps} to functional outputs that {extensively} vary, a particularly challenging feature that {requires addressing using} stochastic discretizations.

{ The analysis of the finite element method (FEM) with uniform meshes 
for numerically approximating PDEs is well established. See e.g. the 
books~\cite{OR2012,GOTZ1975,BWS2010,BCO1981,O1991}.}
FEM is {a major} numerical method for solving PDEs {because of} 
its accuracy, efficiency, stability, and versatility. However, isotropic or anisotropic, uniform meshes are not always optimal. For instance, {for} singularities induced by the problem geometry or when considering PDEs with highly spatially-varying coefficients, there are h-adaptive \cite{OP1999,OP2001,PO1999,OPD2005,KPCL2017,FLTP2006,adFEM_our,BURG2015125,PO2003,AO1997,Ainsworth_Oden_book_2003,oden_2018} and hp-adaptive \cite{DORH1989,ODRW1989,ROD1989,AO1992,O1994} strategies, based on local error indicators, {which} result in higher accuracies through {the} efficient allocation of the degrees of freedom (DoF) in the meshes. 
When approximating $\E{\QoI(u)}$, it is preferable to adapt the mesh with respect to the error of approximating 
this scalar QoI rather than the whole PDE solution. 
{This approach is} known as goal-oriented adaptivity~\cite{OP2002,OV2000,VO2001}. 
Our {aim} is to pursue goal-oriented adaptivity for multilevel Monte Carlo.
{For our approach, the} class of coefficients {requires sufficient} pathwise regularity. Problems whose coefficients are almost surely and everywhere rough {must} be treated with other techniques and are outside the scope of this {study}, see for {example} \cite{Hall_Hoel}, which studies approximation {using} uniform discretizations.


Heinrich \cite{HEINRICH1998,hs99} introduced multilevel Monte Carlo for applications {to} parametric integration. Motivated by applications in computational finance, Kebaier~\cite{kebaier05} introduced a two-level control variate technique in Monte Carlo (MC) sampling for the
weak approximation of stochastic differential equations (SDEs).
{Giles ~\cite{Giles_OpRes} extended t}his approach to the now-famous {m}ultilevel Monte Carlo (MLMC) using a full
hierarchy of discretizations with geometrically decreasing grid sizes.
By optimally {selecting} the number of samples on each level and sampling more {from} coarse, inexpensive levels,
and less {from} fine, expensive levels, the MLMC
method decreases the computational cost.
This cost reduction with respect to single{-}level MC usually goes beyond a constant factor{, unlike the cost reduction of} 
standard control variate techniques. MLMC can reduce the computational complexity to compute a
 solution with {an} error tolerance $\tol > 0$, cf. Theorem~\ref{thm:MLMC_compl}.
Central {l}imit results are useful {for estimating} and control{ling} the statistical error in the MLMC in terms of its variance.
These results, cf. \cite{Kebaier_CLT4MLMC,continuation_MLMC_our} and the generalization in \cite{Hoel_CLT}, are applications of the Lindeberg {c}entral {l}imit {t}heorem {because} the MLMC samples are not identically distributed across levels.

To achieve the optimal MLMC complexity rate, sufficiently accurate estimates of the variance on each level and the bias on the deepest level {are required}.
  
To preserve the {theoretical} complexity, any practical MLMC algorithm {must} produce variance and bias estimates without incurring {considerable} overhead {cost}. Giles proposed an algorithm \cite{Giles_OpRes} {that addresses this by} successively increas{ing} the number of levels and produc{ing} sample variance estimates across levels and {the} corresponding bias.
His approach uses an arbitrarily fixed accuracy splitting between the bias and the statistical errors. This choice appears in similar versions elsewhere in the literature, see~\cite{Teckentrup_2013,gs13c,gs13b,Cliffe_2011}.
Alternatively, the continuation MLMC (CMLMC) algorithm, proposed in~\cite{continuation_MLMC_our},
calibrates models for variance and weak convergence, and the average computational work per sample. 
 The CMLMC algorithm solves the given
problem for a sequence of decreasing tolerances, ensuring that the cost of the sequence of problems is dominated by the
{expense} of the problem we originally {intended} to solve.
Thus, tolerance {is the continuation parameter here}, and the algorithm stops when the
 critical error tolerance is achieved with a prescribed high probability level. 
Creating this auxiliary sequence of smaller problems {enhances the learning of}
 all necessary parameters to {efficiently} run MLMC and optimize the bias and variance contributions to minimize the computational work.
The area of MLMC methods is still very active. See~\cite{AcNum_MLMC} for a broader view.

%

MLMC has been successful in RPDEs. {Considerable} attention has been devoted to controlling the error in the challenging case of a lognormal distribution in the diffusivity coefficient, $a$, cf. \eqref{eq:pde_general}. Error estimates have been presented for {the} numerical approximation of solutions, or functionals thereof, to elliptic PDEs with random coefficients \cite{Charrier_2012,CD2013,cst13,Cliffe_2011,Teckentrup_2013,MLMC_FEM_Schwab,Hall_Hoel} with varying assumptions, in particular related to the random diffusivity coefficient, physical domain, and boundary conditions. Lognormal coefficient fields are considered in \cite{MLMC_FEM_Schwab} for the case of full spatial regularity and uniform coercivity of the random field, and in \cite{Charrier_2012,CD2013,cst13,Cliffe_2011} for limited spatial regularity and without assum{ing} uniform coercivity. Teckentrup et al. \cite{Teckentrup_2013} extended the error analysis provided by Charrier et al. \cite{cst13} to address challenges {because of} non-smooth physical domains. For MLMC, Charrier et al. \cite{cst13} derived uniform error bounds {for} approximating the expectation of a smooth functional of the random solution in $\mathcal{C}^2$ bounded Lipschitz domains. MLMC convergence results have been derived for real-valued, bounded linear functionals of the random solution of elliptic problems on $[0,1]^d$ domains {using} Dirichlet boundary conditions in \cite{Cliffe_2011} and on particular types of non-smooth domains in \cite{Teckentrup_2013}.
Teckentrup et al.~\cite[Theorem 2.3]{Teckentrup_2013} show that there is an optimal work convergence rate that 
is similar to the previously mentioned complexity rates, cf. Theorem \ref{thm:MLMC_compl}, but depends on the relation 
between the rate of variance convergence of the discretization method of the underlying equation and work complexity associated
with generating a single sample of the quantity of interest. In certain cases, 
the computational complexity {is of} the optimal rate, namely $\mathcal{O}\left(\tol^{-2}\right)$.

Multilevel Quasi-Monte Carlo (MLQMC) methods with uniform mesh refinements have been proposed in \cite{Graham_et_al_Num_Math_2015,HS2019,KSSSU2017} and are outside the scope of this {study}.

%
%

Within MLMC, the notion of adaptivity {takes} different forms. For instance, we may optimally distribute the MLMC hierarchy, an issue that {remains} relevant even for uniform discretization. {Thus,~\cite{optimal_hierarchies_our} contains} a general optimization of the MLMC parameters in the discretization hierarchy, {especially} the separation of levels and {optimized} the accuracy {of} budget allocation between the bias and the statistical
errors. In MLMC, we may consider {adaptively selecting} the time-stepping method to save computational work. This {is demonstrated in}
\cite{MLMC_Chernoff_2014, MLMC_reaction_split_2016,MLMC_split_2016} in the context of {s}tochastic {r}eaction {n}etworks. The notion of adaptivity exploited in this work is the mesh size control. Here, the step size or the mesh size {is selected} to minimize cost while achieving a prescribed accuracy. The work~\cite{adaptive_MLMC_1_our, adaptive_MLMC_2_our} first introduced MLMC stepsize goal-oriented adaptivity for  
It\^o SDEs, 
{ defining the MLMC levels in terms of accuracy requirements achieved by a suitable adaptive algorithm, rather than
in terms of uniform discretizations}. This idea is general and is also useful in MLMC for RPDEs.

%
%

Eigel et al.~\cite{MLMC_AdFEM_Eigel_2016}, the closest reference to this work, used goal-oriented adaptivity in the context of MLMC for~\eqref{eq:bvp_general}, assuming uniform coercivity, an assumption {for which in this work} we relax to treat the lognormal case. 
Their approach {differs} from ours {because} (i) it is based on batch adaptivity, (ii) the error indicators are different, do not yield an error density, and seem less sharp, cf.~\cite{MS2009}, (iii) they based their error control on the MSE error and did not separate the statistical error and bias errors as we {did}.
On a less related note, Scarabosio et al. \cite{SWOF2019} used goal-oriented model adaptivity based on a hierarchical control variate using two and three levels. {They} based their error estimates on 
{verifying} the current model output pathwise against a higher-fidelity model. For adaptive MLMC {using} norm error control we refer to \cite{MLMC_AdFEM_Korhhuber_2018} and the PhD thesis~\cite{MLMC_AdFEM_phdthesis_2018}. The goal in those studies was to perform pathwise adaptivity for functional approximation for the
expected value of the solution and not to compute a quantity of interest as we do here. 
The more recent work~\cite{LANG2020109692} is focused on adaptivity with multilevel stochastic collocation for solving elliptic PDEs with random data by combining adaptive (anisotropic) sparse Smolyak grid approximations in the stochastic space, using the Sparse grid Matlab kit (version 17-5) \cite{PT2022}, for each collocation point {by} applying the dual weighted residual method \cite{becker_rannacher_2001} for spatial goal-oriented adaptive mesh refinements.

 Our contributions in this work are as follows: we propose a novel AMLMC for 
\eqref{eq:bvp_general}
 with geometric singularities and a lognormal random diffusivity coefficient. We build our AMLMC on~\cite{adFEM_our}, which  
 developed weak convergence rates for an adaptive algorithm {using} isoparametric d-linear quadrilateral finite element approximations and the dual weighted residual error representation in {a} deterministic setting. This provides us with sharp error estimates and indicators {for creating} locally refined sequences of meshes tuned to the geometry-driven singularity at hand.
{Compared with} batch-adaptive MLMC approaches \cite{MLMC_AdFEM_Eigel_2016, oden_2018},
our AMLMC circumvents the {expense} of generating meshes on the fly. 
 Furthermore, {our application of} stochastic meshes {is} more efficient than {the} batch-adaptive algorithm {for cases without} uniform coercivity, like in the lognormal case.
{Moreover}, for cases where the solving cost is higher than the assembly cost, we discuss the use of iterative solvers and compare their efficiency with direct ones. To save computational work, we propose a goal-oriented stopping criterion for the iterative solver.

Theoretically, we characterize those problems where adaptivity provides a noticeable advantage, namely those where the error density blows up in $L^1_P(\mathcal{D}\times \Omega)$, as 
we refine the mesh around the singularities. We provide a slight generalization of the MLMC complexity theorem, cf. \eqref{eq:complexity_MLMC}, in Theorem~\ref{thm:MLMC_compl}, {thus} allowing the MLMC base level to converge to the exact solution, although at a relatively slow rate while essentially preserving the complexity in \eqref{eq:complexity_MLMC}. 
Corollary~\ref{cor:AMLMCComplexity} {further provides} an estimate for the computational work of AMLMC by identifying the multiplicative constants in terms of the quasi-norm of problem-dependent error density.
The pointwise convergence of the error density as $\tol\to 0$, based on proper local averages, is {fundamental} to {delivering} theoretical results on stopping, asymptotic accuracy, and efficiency, as proved in \cite{adFEM_our} and inherited naturally by our AMLMC {algorithm}. 
{Furthermore}, the efficient computation of these local 
averages in the error estimate is a contribution of this work, improving {the results in} \cite{adFEM_our}.

Our {2D} numerical experiments yield results consistent with the theoretical predictions, namely, (i) the a posteriori error estimates are sharp for bias and variance predictions {with} singularities and lognormal coefficients {present}, and (ii) AMLMC exhibits advantages compared to SMLMC, both on its complexity and ability to estimate the errors. More substantial computational gains are expected in higher dimensions.

The outline of this {study} is as follows: Section~\ref{sec:theory} {presents} the problem setting by first describing the boundary value problem under consideration, which is a particular class of elliptic PDEs with random coefficients. {Then,} the goal-oriented adaptive FEM strategy {is introduced and the} MLMC method with uniform refinements for estimating the expected values of linear functionals of the random PDE solution {is formulated}. {AMLMC, our goal-oriented adaptive finite element MLMC algorithm, is presented in Section~\ref{sec:adapt_MLMC}. Section~\ref{sec:numex} presents first the numerical implementation details of AMLMC, then} numerical results of {applying} the AMLMC algorithm to {a 2D boundary value problem of} an elliptic PDE with a stochastic coefficient field, not bounded away from zero nor bounded from above, and with {a} geometric singularity at the meeting point between Dirichlet and Neumann boundaries.  We observe that, for the same accuracy requirement, our AMLMC requires less work and is more reliable than MLMC on uniform meshes.
Finally, {Section~\ref{sec:concl} offers} the conclusions.


\section{Background material}
\label{sec:theory}

\subsection{Problem setting}\label{sec:problem_setting}

We consider the problem of approximating the expected value of a 
scalar quantity of interest (QoI), which is a function of the solution
to{~\eqref{eq:bvp_general},} a boundary value problem (BVP) for a linear second-order PDE with a stochastic coefficient field.
For a given $h$-adaptive FEM, we will use a 
sample-based method to approximate the expected value. 

\subsubsection{{BVP} for an elliptic PDE with smooth stochastic data}

Let $(\Omega,\mathcal{F},\mathbb P)$ be a complete 
probability space and let $\mathcal{D}$ be an open and bounded 
polygonal/polyhedral domain in $\rset^d$, $d\ge 2$, 
with boundary $\partial\mathcal{D}$.
With homogeneous Dirichlet and Neuman boundary conditions, 
the BVP with random data is to find 
$u: \mathcal{D} \times \Omega \to \rset$ that almost surely (a.s.)
for $\omega \in \Omega$ solves~\eqref{eq:bvp_general}. 
The model diffusivity coefficient considered here takes the form
\begin{align}
	a(x, \omega) & = \exp \left( \sum_{i\ge0} \xi_i(\omega) \sqrt{\lambda_i} \theta_i(x) \right), & x \in D,
	\label{eq:random_field_series0}
\end{align}
where $(\lambda_i)_{i\ge0}$ is a monotone decreasing sequence of non-negative numbers,  $(\theta_i)_{i\ge0}$ are trigonometric functions in $\rset^d$,
and $(\xi_i)_{i\ge0}$ are i.i.d. $\mathcal{N}(0, 1)$ random variables.
For the well-posedness and solution regularity for 
problem~\eqref{eq:bvp_general} with a lognormal 
	coefficient field \eqref{eq:random_field_series0}, we refer to~\cite{Charrier_2012,BNT_SIAM_review,Graham_et_al_Num_Math_2015,Bachmayr_et_al_M2AN_2017,Nguyen_Nuyens_M2AN_2021}.
	There, well posedness requires control on the pathwise coercivity cf. \cite{BNT_SIAM_review,Teckentrup_2013}, given by Assumption~\ref{assump1}.
\begin{assumption}[Well-posedness]\label{assump1}
To ensure the well posedness of the solution to~\eqref{eq:bvp_general} in 
$L^k_{\mathbb P}(\Omega;H^1(\mathcal{D}))$ for any $k\ge 1$, we assume that
  \begin{itemize}
  \item $a_{\min}(\omega) = \min_{\mathbf x \in \mathcal  D}a(\mathbf x; \omega) > 0$ a.s. and $1/a_{\min} \in  L^q_\mathbb{P}(\Omega)$, for all $q \in (0, \infty)$.
  \item $a \in L^q_{\mathbb P}(\Omega, C^1(\overline{\mathcal D}))$,
    for all $q \in (0, \infty)$.
  \item $f \in L^{q^\ast}_{\mathbb P}(\Omega, L^2(\mathcal D))$ for some
    $q^\ast > 2$.
  \end{itemize}
\end{assumption}
Here, $L^q_\mathbb{P}(\Omega)$ denotes the space of random variables with
  finite $q$-th moment in $(\Omega,\mathcal{F},\mathbb P)$. We will later also use 
  the notation $L^q_\mathbb{P}(\Omega\times\mathcal{D})$ when the integration is with 
  respect to the product measure of $\mathbb{P}$ and the Lebesgue measure.

Assumption~\ref{assump1} is satisfied by our numerical examples.
	However, pathwise regularity connects with the $p$-summability of sufficiently high order derivatives of the auxiliary sum 
	$\sum_{i\ge0}  \sqrt{\lambda_i} \theta_i(x)$.
Indeed, the fact that the coefficient function~\eqref{eq:random_field_series0} is not uniformly
bounded away from 0 makes the analysis more technically challenging. 
As will be elaborated on later,  this feature calls
for stochastic mesh selection because $\QoI(u)$ is unbounded.

The QoI is assumed to be given by a deterministic, real-valued, 
bounded linear functional $\QoI$ of the stochastic solution $u$, 
and the aim is to approximate its expected value, $\E{\QoI(u)}$, to a given 
accuracy $\tol$ and confidence level. 

\subsubsection{Goal-oriented adaptive FEM}

In many engineering applications, the primary goal of a computation
is to estimate expected values of functionals of the solution of the BVP. The standard deterministic technique for error estimation, 
applied pathwise in this work, is the use of a dual weighted error representation, cf. Section~2.1 in~\cite{adFEM_our}.
For more on the dual weighted residual (DWR) approach, see, 
e.g.,~\cite{becker_rannacher_2001,Ainsworth_Oden_book_2003,Bangerth_Rannacher_book_2003}.
The approximation error in a linear functional of the solution is then estimated by
$$|\QoI(u)-\QoI(u_\mathrm{FEM})| \leq \sum_{i}|r_i| |z_i|,$$
where $u_\mathrm{FEM}$ is the FEM approximation of $u,$ 
the sum is taken over all elements, and $r_i$ and $z_i$
are elementwise contributions from a residual and the corresponding dual solution, respectively.
Following the analysis in~\cite{adFEM_our}, this study uses 
an error expansion with computable leading order term, namely,
 $$\QoI(u)-\QoI(u_\mathrm{FEM}) \sim \sum_{i} \bar{\rho}_i h_i^{p+d},$$
cf.~\eqref{eq:rho_tilde}. Here, $d$ is the dimension of the domain, and $p$ is the approximation order of the FEM.
The primary advantage of this expansion 
is that the approximate error density $\bar{\rho}_i$, computed using both primal and dual solutions, is essentially independent of the mesh size and the theoretically optimal mesh size can be expressed in terms of 
the asymptotic error density, $\rho$, cf. Theorem 2.1 and Corollary 2.2 in \cite{adFEM_our}. 
Letting $h(x):\mathcal{D}\to\rset^+$ be the mesh size function, the error in the output functional can then be approximated as
\begin{align}
  \label{eq:leading_order_bias}
	\QoI(u)-\QoI(u_\mathrm{FEM}) & \simeq 
	\int_{\mathcal D} \rho(x)h(x)^{p}\,dx.
\end{align}
Similarly, the computational cost for the assembly and solution steps can be expressed using the mesh function.
Thus, a natural step is to formulate a variational problem to choose the optimal mesh function subject to an error budget constraint,
see Section~\ref{sec:adapt_MLMC}. 
This formulation provides a quantitative measure to understand the potential of adaptivity on a given 
problem and guides the construction of related adaptive algorithms. 
In this work, the deterministic adaptive algorithm developed in~\cite{adFEM_our} 
  is used pathwise. The pointwise convergence of the error density was proven in~\cite{adFEM_our} where,  
  using ideas from~\cite{Moon03,adOverview_our,Moon05_2}, it was exploited to show the stopping, accuracy, 
  and efficiency of the deterministic adaptive algorithm.
The error density-based goal-oriented adaptive algorithm we employ has a universal nature. Its success in optimal control problems, cf.~\cite{Karlsson2015}, is an example of this versatility. 

Note that this notion of error densities 
can be extended to nonlinear problems and nonlinear observables, 
as described in~\cite[Remark 2.4]{adFEM_our},
and higher-order approximations,~\cite[Section 2.4]{adFEM_our}.

In what follows, we will use these $h$-adaptive FEM algorithms to define the AMLMC levels, 
characterizing them by a sequence of error tolerances following  \cite{adaptive_MLMC_1_our, adaptive_MLMC_2_our}. These error tolerances
will control  both weak and
strong errors. Furthermore, we will focus on examples where 
locally refined meshes provide an advantage due to
geometry induced singularities, a commonplace in engineering
applications.

\subsubsection{Multilevel Monte Carlo}
\label{sec:std_MLMC}

MLMC uses hierarchical control variates to substantially reduce the computational cost of MC.
Thus, it uses a hierarchy of $L+1$ 
meshes defined by decreasing mesh sizes,
indexed by their levels $\ell=0,1,\dots,L$,
and the telescopic representation of the  expected value of the finest
approximation. Thus, letting $\QoI_{\ell}$ denote the 
approximation of $\QoI$ on level $\ell$,
\begin{align*}
  \E{\QoI_L} & =
  \E{\QoI_0} + \sum_{\ell=1}^L \E{\QoI_{\ell}-\QoI_{\ell-1}}.
\end{align*}
Then, the MLMC estimator approximates each 
expected value in the previous equation by a sample average, namely,
\begin{align}
  \label{eq:MLMC_estimator}
  \EstMLMC & = \frac{1}{\nrs_0}\sum_{n=1}^{\nrs_0}\QoI_0(\omega_{0,n}) +
  \sum_{\ell=1}^L\frac{1}{\nrs_\ell}\sum_{n=1}^{\nrs_\ell}
  \left(\QoI_{\ell}(\omega_{\ell,n})-\QoI_{\ell-1}(\omega_{\ell,n})\right),
\end{align}
where  $\{\omega_{\ell,n}\}_{n=1,\dots,\nrs_\ell}^{\ell=0,\dots,L}$
denote i.i.d. realizations of the mesh-independent random variables. 
In MLMC, it is fundamental to evaluate  both terms in the difference
\begin{align}
  \label{eq:Delta}
  \Delta\QoI_{\ell}(\omega_{\ell,n}) & = 
  \QoI_{\ell}(\omega_{\ell,n})-\QoI_{\ell-1}(\omega_{\ell,n})
\end{align}
with the same outcome of $\omega_{\ell,n}$. 
This implies that
$\var{\Delta\QoI_{\ell}}\to 0$, as $\ell\to\infty$, provided that the
numerical approximation $\QoI_{\ell}$ converges strongly. 

To simplify the presentation, introduce
the notation 
\begin{align}
  \label{eq:Var_ell}
  V_\ell & = 
  \begin{cases} 
    \var{\QoI_0}, & \ell=0,\\
    \var{\QoI_\ell - \QoI_{\ell-1}}, & \ell>0.
  \end{cases}
\end{align}

The average computational work for generating $\EstMLMC$ is
\begin{align}
  \label{eq:work_sum}
  \work_\mathrm{MLMC} & = \sum_{\ell=0}^{L} \nrs_\ell \work_\ell,
\end{align}
where $W_\ell$ is the average cost of generating $\Delta\QoI_\ell$.
We build the estimator~\eqref{eq:MLMC_estimator} to approximate $\E{\QoI}$ with accuracy $\tol>0$ with high probability.
Therefore, we introduce a splitting parameter $\theta\in(0,1)$ and enforce
    \begin{align}\label{eq:BIASC}
        \text{Bias Constraint: } |\mathrm{E}[\QoI] - \E{\QoI_L}|&\leq (1-\theta) \tol,\\
        \label{eq:STATC}
        \text{Statistical Constraint: } Var[\EstMLMC] =   \sum_{\ell=0}^L \frac{V_\ell}{\nrs_\ell} &\leq \left(\frac{\theta \tol}{\confpar}\right)^2,
    \end{align}
where the parameter $\confpar>0$ is chosen depending on the desired confidence level in the error control. Thus, the bias constraint
determines $L(\tol)$ whereas the statistical constraint, 
which is motivated by a central limit theorem approximation, cf.~\cite{Hoel_CLT}, 
determines the optimal number of samples across levels.
    
Given $L$ and $\splitting$, minimizing the work~\eqref{eq:work_sum}
subject to the constraint~\eqref{eq:STATC} and relaxing the integers into real values, leads to the
optimal number of samples per level in $\EstMLMC$,
\begin{align}
  \label{eq:opt_nrs}
  \nrs_\ell & =
  \left( \frac{\confpar}{\splitting\tol} \right)^2
  \sqrt{\frac{V_\ell}{\work_\ell}} \sum_{k=0}^L \sqrt{\work_k V_k}.
\end{align}

Substituting the ceiling of these optimal $\nrs_\ell$  in the total work \eqref{eq:work_sum} yields the upper bound:
\begin{align}
  \label{eq:total_work_MLMC}
  \work_\mathrm{MLMC} & =
  \left( \frac{\confpar}{\splitting\tol} \right)^2
  \left( \sum_{\ell=0}^L \sqrt{\work_\ell V_\ell} \right)^2 +  \sum_{\ell=0}^L {\work_\ell } .
\end{align}
Under the standard assumptions for MLMC, the first term in the last equation dominates the second one as $\tol \to 0$.
Thus, since $L(\tol)\to \infty$ as $\tol\to 0$, the different computational complexity regimes  of MLMC are determined by
the behavior of  $\sum_{\ell=0}^L \sqrt{\work_\ell V_\ell}.$ Indeed, if $\sum_{\ell=0}^\infty \sqrt{\work_\ell V_\ell}<\infty $ then the MLMC complexity
is $\tol^{-2}$, while if $\work_\ell V_\ell = \mathcal{O}(1)$, 
then  the MLMC complexity is the larger $\tol^{-2}L^2(\tol)$. 
Finally, if $\work_\ell V_\ell\to \infty$, the MLMC complexity deteriorates even further. 
For example, if the product $\work_\ell V_\ell$ grows geometrically, the MLMC complexity 
is dominated by the deepest level, namely, $\tol^{-2}\work_{L(\tol)} V_{L(\tol)}$, which still improves on the corresponding MC complexity.

Following Giles, in the MLMC standard analysis, 
one usually makes the even more precise assumptions:
\begin{assumption}\label{assu:MLMC}
Assume, for some $\kappa>1$, a variance convergence rate $V_\ell  \lesssim \kappa^{-\Os\ell}$, a bias convergence rate $|\mathrm{E}[\QoI] - \E{\QoI_\ell}| \lesssim \kappa^{-\Ow\ell}$ and a work per sample $W_\ell \lesssim \kappa^{\gamma\ell}$. Moreover, assume that
\begin{equation}\label{eq:Giles_constr}
\Ow\geq\min{\left(\Os,\gamma\right)}/2>0.
\end{equation}
\end{assumption}
When Assumption~\ref{assu:MLMC} is satisfied, the MLMC
computational work has, as $\tol\to 0$, the asymptotic behavior
\begin{align}
  \label{eq:complexity_MLMC}
  \work_\mathrm{MLMC} & \lesssim
    \begin{cases}
      \tol^{-2}, & \text{if $\Os>\gamma$},\\
      \tol^{-2}\left(\log{\tol^{-1}}\right)^2, & \text{if $\Os=\gamma$},\\
      \tol^{-2\left(1+\frac{\gamma-\Os}{2\Ow}\right)}, & \text{if $\Os<\gamma$},
    \end{cases}
\end{align}
see e.g. Theorem~3.1 in~\cite{Giles_OpRes}, or
Corollaries~2.1 and~2.2 in~\cite{optimal_hierarchies_our}.
Note that~\eqref{eq:Giles_constr} is used to ensure that the 
second term in~\eqref{eq:total_work_MLMC} is asymptotically negligible compared to the first term. 
	
We now state a slightly generalized version of the MLMC complexity theorem allowing the coarsest 
level to be refined as $\tol \to 0.$ This is particularly relevant when considering the complexity 
of adaptive MLMC algorithms because they usually rely on asymptotic properties, 
which become sharper as the largest element in the discretization is refined to ensure 
consistent approximation.

\begin{theorem}[MLMC Complexity with variable level $0$]
  \label{thm:MLMC_compl} 
  Take $0\le\epsilon<<1.$
  Let Assumption~\ref{assu:MLMC} hold and choose $L(\tol)$ to be the minimum integer that satisfies the bias constraint. If the MLMC base level $0$ is such that its average work per sample satisfies
  \begin{align}
    \label{eq:level0_MLMC2}
    W_0 & \lesssim
    \begin{cases}
      \left(\log\tol^{-1}\right)^{\epsilon}, & \text{if $\Os>\gamma$},\\
      \left(\log\tol^{-1}\right)^{2}, & \text{if $\Os=\gamma$},\\
      \tol^{-\frac{\gamma-\Os}{\Ow}}, & \text{if $\Os<\gamma$},
    \end{cases}
  \end{align}
the asymptotic MLMC computational work then satisfies
  \begin{align}
  \label{eq:complexity_MLMC2}
  \work_\mathrm{MLMC} & \lesssim
    \begin{cases}
      \tol^{-2} \left(\log\tol^{-1}\right)^{\epsilon}, & \text{if $\Os>\gamma$},\\
      \tol^{-2}\left(\log{\tol^{-1}}\right)^2, & \text{if $\Os=\gamma$},\\
      \tol^{-2\left(1+\frac{\gamma-\Os}{2\Ow}\right)}, & \text{if $\Os<\gamma$},
    \end{cases}
\end{align}
cf.~\eqref{eq:complexity_MLMC}.
\end{theorem}
{
\begin{proof} 
  The proof is a direct consequence of the standard MLMC complexity theorem.
  It amounts to bounding the relative increase of the MLMC work caused by
  dropping the first $k(\tol)$ levels, out of $L(\tol)$, given an arbitrary
  fixed $0$-level. Ignoring the second term of~\eqref{eq:total_work_MLMC}, 
  which is asymptotically negligible as $\tol\to 0$ by Assumption~\ref{assu:MLMC}, 
  this amounts to bounding the relative increase in 
  $S_L = \sum_{\ell=0}^L\sqrt{V_\ell W_\ell}$.
  To analyse this increase, choose $L = L(\tol)$ to satisfy the bias constraint,
  let $\widehat V_k = \var{\QoI_k} \approx V_0$ be the single
  level variance, for all $k<L$, and assume for simplicity of notation that the
  costs of sampling $\QoI_\ell$ and $\Delta\QoI_\ell$ are the same.

  The MLMC complexity remains unchanged provided that, for some $K>0$,
  \begin{align*}
    \sqrt{\widehat V_k W_k}+\sum_{\ell=k+1}^L\sqrt{V_\ell W_\ell}
    & \leq K
    \sum_{\ell=0}^L\sqrt{V_\ell W_\ell},
  \end{align*}
  or equivalently,
  \begin{align}
    \label{eq:demand}
    \sqrt{\widehat V_k W_k}
    & \leq (K-1)\sum_{\ell=0}^L\sqrt{V_\ell W_\ell} + \sum_{\ell=0}^k\sqrt{V_\ell W_\ell}.
  \end{align}
  For $\Os\leq\gamma$, when $S_L$ is unbounded as $\tol\to 0$, we can 
  satisfy~\eqref{eq:demand} while $k\to\infty$ as $\tol\to 0$ using the following
  specific choices.

  For $\Os<\gamma$, where the summands $\sqrt{V_\ell W_\ell}$ increase
  with $\ell$, choosing $k(\tol)$ such that
  \begin{align*}
    W_{k(\tol)} & \leq 
    {\widehat C}^2 V_{L(\tol)} W_{L(\tol)}
  \end{align*}
  satisfies~\eqref{eq:demand} with $K=1+\widehat C\sqrt{\widehat V_k}$. The fact that 
  $V_{L(\tol)} W_{L(\tol)}\to\infty$ as $\tol\to 0$, 
  confirms that it is a feasible to let $k(\tol)\to\infty$ as $\tol\to 0$.

  When $\Os=\gamma$ it holds that $\sum_{\ell=1}^L\sqrt{V_\ell W_\ell}=L\sqrt{V_1 W_1}$
  and consequently choosing $k(\tol)$ such that
  \begin{align*}
    W_{k(\tol)} & \leq \left(\widehat C L(\tol)\right)^2
  \end{align*}
  satisfies~\eqref{eq:demand} with $K=1+\widehat C\sqrt{\widehat V_k/(V_1 W_1)}$.
  Since $L(\tol)\simeq\log{\left(\tol^{-1}\right)}$ we can satisfy
  the asymptotic bound by letting 
  $W_{k(\tol)} \leq \widehat C( \log{\left(\tol^{-1}\right)})^2$
  for some other $\widehat C$.

  In the case $q_s > \gamma$, 
  where $S_\infty=\sum_{\ell=0}^\infty\sqrt{V_\ell W_\ell}$ is finite, 
  we need to accept a slightly increased complexity if we let $k\to\infty$ 
  as $\tol\to 0$. To quantify this increase, let $f:\rset^+\mapsto\rset^+$
  be a strictly increasing, unbounded, function and require that
  \begin{align*}
    W_{k(\tol)} & \leq f(\tol^{-1}).
  \end{align*}
  Since $\sum_{\ell=0}^{L(\tol)}\sqrt{V_\ell W_\ell}\to S_\infty$ as $\tol\to 0$,
  we can find a constant $K>1$, depending on both $\widehat V_k$ and $S_\infty$, 
  such that
  \begin{align*}
    \sqrt{\widehat V_k W_{k}} & \leq 
    \left(K-1\right)\sqrt{f(\tol^{-1})} \sum_{\ell=0}^{L(\tol)}\sqrt{V_\ell W_\ell},
  \end{align*}
  for all sufficiently small $\tol$. 
  Assuming further that $\tol<\left(f^{-1}(1)\right)^{-1}$, this implies
  \begin{align*}
    \sqrt{\widehat V_k W_k}+\sum_{\ell=k+1}^{L(\tol)}\sqrt{V_\ell W_\ell}
    & \leq K\sqrt{f(\tol^{-1})}
    \sum_{\ell=0}^{L(\tol)}\sqrt{V_\ell W_\ell},
  \end{align*}
  and consequently the MLMC complexity increases to $\tol^{-2}f(\tol^{-1})$, at worst.
  Choosing $f(x)\propto\left(\log{x}\right)^{\epsilon}$, $x\in[1,\infty)$, corresponds to the 
  statement of this theorem and completes the proof.
\end{proof}
}

\begin{remark}[MLMC with deterministic meshes]
When approximating the solution of an RPDE with finite elements on regular triangulations, the 
parameter $h>0$ refers to either the maximum element diameter or another characteristic length. 
Here, because we use quadrilateral isotropic elements, $h$ will be the side length.  
Denoting the corresponding
approximate solution by $u_h(\omega)$ and using
 piecewise linear or piecewise $d$-multilinear continuous finite element
approximations, and with the previous assumptions, 
it can be shown~\cite[Corollary 3.1]{Teckentrup_2013} that, in the absence of singularities, asymptotically as
$h \to 0$:
\begin{itemize}
\item $\left|\E{\QoI(u)-\QoI(u_h)} \right| \lesssim Q_W\, h^2$ for a constant $Q_W>0$.
\item $\var{\QoI(u)-\QoI(u_h)} \lesssim Q_S\, h^4$ for a constant $Q_S > 0$.
\end{itemize}
Furthermore, assuming an average cost per sample of the order of $h{-(d+\gamma)}$, 
with $1\le \gamma \le 3,$ a direct application of Theorem \ref{thm:MLMC_compl} yields the complexity 
$$
  \work_\mathrm{MLMC}  \propto \tol^{-\max(2,{d\gamma/2})}
(\log{\tol})^{2 \chi}
$$ 
with $\chi =1$ if $d\gamma-4=0$  and $\chi=0$ elsewhere. The corresponding MC method has the worse complexity
$$
  \work_\mathrm{MC}  \propto \tol^{-(2+{d\gamma/2})}.
$$ 
Observe that in the presence of singularities, both the weak and strong convergence rates decrease, 
and SMLMC yields a deteriorated complexity. This will be the case for the examples considered in 
Section~\ref{sec:example_problems}. There,
AMLMC will yield a complexity corresponding to smooth problems, provided the corresponding quasi-norm of the error density remains bounded. 
\end{remark}

\section{Adaptive Multilevel Monte Carlo}
\label{sec:adapt_MLMC}


For this analysis, we assume that the average work of generating one 
approximate sample of $\QoI$ is approximately proportional to the 
average number of elements $\nelem$, of the FEM discretization; either
because the linear solver is optimal or the assembly cost 
dominates the cost of the linear solver for the relevant range of
tolerances{, as in Example~2 in Section~\ref{sec:numex}.}

We write the approximation of the expected number of elements as
\begin{align}
	\label{eq:nelem_model}
	\E{\nelem} & = \int_\mathcal{D} \E{h^{-d}}
\end{align}
where $h:\mathcal \mathcal{D}\times\Omega\to\rset^+$ is a stochastic mesh size function used in the 
FEM
approximation; moreover, we assume that the bias is approximated
by the first term in the expansion, cf.~\eqref{eq:leading_order_bias},
\begin{align}
	\label{eq:bias_model}
	\int_\mathcal{D} \E{\rho h^p}
\end{align}
for a non-negative
error density $\rho\in L_\mathbb{P}^{\frac{d}{p+d}}(\mathcal{D}\times\Omega)$; 
cf. Theorem~2.1 in~\cite{adFEM_our}.
By assuming $\rho\geq 0$, we ignore any cancellation of error contributions at this stage;
both~\cite{adFEM_our} and the numerical experiments in Section~\ref{sec:numex}
involve $\tol$-dependent positive lower and upper bounds on $\rho$, 
which we omit here.
The bias model~\eqref{eq:bias_model} is an approximation both in the sense 
that higher-order terms in the error expansion are omitted and that 
unattainable continuous $h$ 
is used. However, it is a useful tool for examining 
what can be achieved using adaptive mesh refinement and sampling.

First, we will present 
the optimal, stochastic, continuously varying mesh 
function $h^\ast:\mathcal \mathcal{D}\times\Omega\to\rset^+$ as a function of the 
stochastic error density; i.e., $h^\ast$ minimizes the expected work required 
to obtain a given bias. 
Then we suggest an adaptive MLMC estimator as an approximation 
to the optimum, which is suitable when a singularity 
induced by the deterministic geometry is the primary reason to use
locally refined meshes.


\subsection{Optimal single level mesh distribution for a given bias tolerance}
\label{sec:opt_SMLMC}

In a single-level MC having a fixed split between bias and statistical errors, different methods differ only via their expected work per sample.
Thus, minimizing the expected work per sample, modeled as the expected number of 
elements~\eqref{eq:nelem_model}, subject to the constraint
that the bias model~\eqref{eq:bias_model} equals a prescribed tolerance 
$\tol_\mathrm{bias}$, via a standard Lagrange multiplier technique
(\ref{app:opt}) 
leads to the theoretically optimal stochastic mesh size
\begin{align}
	\label{eq:opt_h_stoch}
	h^\ast(x;\omega) & = 
	\frac{\tol_\mathrm{bias}^{1/p}}{\left(\int_\mathcal{D} \E{\rho^{\frac{d}{p+d}}}\right)^{1/p}}
	\rho(x;\omega)^{-\frac{1}{p+d}}
\end{align}
with an associated optimal mean work per sample proportional to
\begin{align}
	\label{eq:opt_work_h_stoch}
	\int_\mathcal{D} \E{(h^\ast)^{-d}} & = 
	\frac{\left(\int_\mathcal{D} \E{\rho^{\frac{d}{p+d}}}\right)^{\frac{p+d}{p}}}
	{\tol_\mathrm{bias}^{d/p}}
	=
	\left(\frac{\left\Vert\rho\right\Vert_{L_\mathbb{P}^{\frac{d}{p+d}}(\mathcal{D}\times\Omega)}}{\tol_\mathrm{bias}}\right)^{d/p}.
\end{align}
Thus, the optimal strategy, for single-level MC, does not evenly distribute 
the error contributions over the samples. Indeed, for this
strategy, the leading order error estimate satisfies
\begin{align}
	\label{eq:error_model_stop}
	\int_\mathcal{D} \rho(x;\omega) h^\ast(x;\omega)^p\,dx & =
	\frac{\tol_\mathrm{bias}}{\int_\mathcal{D} \E{\rho^{\frac{d}{p+d}}}}
	\int_\mathcal{D} \rho(x;\omega)^{\frac{d}{p+d}}\,dx
\end{align}
for each sample. Equation~\eqref{eq:error_model_stop} is the 
motivation for the AMLMC algorithm proposed in this section.

Similarly, minimizing the work model~\eqref{eq:nelem_model} under
the same bias constraint using only uniform meshes leads to
\begin{align}
	\label{eq:opt_h_uniform}
	h_\mathrm{uni}^\ast(\omega) & = 
	\tol_\mathrm{bias}^{1/p}\frac{\left(\int_\mathcal{D}\rho(x;\omega)dx\right)^{-\frac{1}{p+d}}}{\E{\left(\int_\mathcal{D}\rho\right)^{\frac{d}{p+d}}}^{1/p}}.
\end{align}
Thus, the work estimate corresponding to~\eqref{eq:opt_work_h_stoch} becomes
\begin{align}
	\nonumber
	\E{(h_\mathrm{uni}^\ast)^{-d}}\int_\mathcal{D}1\,dx  & = 
	\frac{\E{\left(\int_\mathcal{D} \rho\right)^{\frac{d}{p+d}}}^{\frac{p+d}{p}}}
	{\tol_\mathrm{bias}^{d/p}}
	\int_\mathcal{D}1\,dx \\
	\label{eq:opt_work_uni_h_stoch}
	& = 
	\left(\frac{\left\Vert\Vert\rho\Vert_{L^1(\mathcal{D})}\right\Vert_{L_\mathbb{P}^{\frac{d}{p+d}}(\Omega)}}
	{\tol_\mathrm{bias}}
	\right)^{d/p} 
	\int_\mathcal{D}1\,dx,
\end{align}
and the leading order error estimate per sample is given by
\begin{align}
	\label{eq:error_model_stop_uni}
	\left(h_\mathrm{uni}^\ast(\omega)\right)^p\int_\mathcal{D} \rho(x;\omega)\,dx & =
	\frac{\tol_\mathrm{bias}}{\E{\left(\int_\mathcal{D}\rho\right)^{\frac{d}{p+d}}}}
	\left(\int_\mathcal{D}\rho(x;\omega)\,dx\right)^{\frac{d}{p+d}}.
\end{align}

Suppose, in contrast, 
the samples are generated using a fixed, deterministic, 
mesh. The optimal mesh functions are then those obtained using 
the deterministic theory in~\cite{adFEM_our} with the expected 
value of the stochastic error density. Here, the 
optimized work model for deterministic adaptive meshes is given by
\begin{align}
	\label{eq:opt_work_h_det}
	\int_\mathcal{D}\left(h_\mathrm{det}^\ast(x)\right)^{-d}\, dx
	& =
	\left(\frac{\left\Vert\Vert\rho\Vert_{L_\mathbb{P}^1(\Omega)}\right\Vert_{L^{\frac{d}{p+d}}(\mathcal{D})}}
	{\tol_\mathrm{bias}}
	\right)^{d/p},
\end{align}
while that for deterministic uniform meshes is given by
\begin{align}
	\label{eq:opt_work_uni_h_det}
	(h_\mathrm{uni,det}^\ast)^{-d}\int_\mathcal{D}1\,dx & = 
	\left(\frac{\left\Vert\rho\right\Vert_{L_\mathbb{P}^{1}(\mathcal{D}\times\Omega)}}{\tol_\mathrm{bias}}\right)^{d/p}
	\int_\mathcal{D}1\,dx.
\end{align}

Repeated use of Jensen's inequality confirms that
\begin{align}
	\label{eq:norm_relations}
	\left\Vert\rho\right\Vert_{L_\mathbb{P}^{\frac{d}{p+d}}(\mathcal{D}\times\Omega)}^{d/p}
	& \leq 
	\begin{Bmatrix}
		\left\Vert\Vert\rho\Vert_{L^1(\mathcal{D})}\right\Vert_{L_\mathbb{P}^{\frac{d}{p+d}}(\Omega)}^{d/p} 
		\int_\mathcal{D}1\,dx \\~\\
		\left\Vert\Vert\rho\Vert_{L_\mathbb{P}^1(\Omega)}\right\Vert_{L^{\frac{d}{p+d}}(\mathcal{D})}^{d/p}
	\end{Bmatrix}
	\leq
	\left\Vert\rho\right\Vert_{L_\mathbb{P}^{1}(\mathcal{D}\times\Omega)}^{d/p}
	\int_\mathcal{D}1\,dx.
\end{align}
Consequently, these same relations hold between the work estimates 
in~\eqref{eq:opt_work_h_stoch}, \eqref{eq:opt_work_uni_h_stoch}, 
\eqref{eq:opt_work_h_det}, and \eqref{eq:opt_work_uni_h_det}, 
as is expected given the relative generality of the optimization
problems.
All four work estimates coincide in the extreme situation when
$\rho$ is constant in $\mathcal{D}\times\Omega$. 
Similarly, for any $\rho\in L_\mathbb{P}^{1}(\mathcal{D}\times\Omega)$, all
four complexity estimates will have the same asymptotic rate as 
$\tol_\mathrm{bias}\to 0$.  This shows that 
little is gained by going from uniform and deterministic meshes
to adaptively refined and/or sample dependent meshes for problems
where 
$\left\Vert\rho\right\Vert_{L_\mathbb{P}^{1}(\mathcal{D}\times\Omega)}$
is well behaved,  that is, it is uniformly bounded over refinement sequences.
However, the $L^{\frac{d}{p+d}}$ quasi-norms 
of error densities can remain bounded
although the corresponding $L^1$ norms are unbounded.
Canonical examples in the FEM setting for~\eqref{eq:bvp_general}
are: (i) deterministically induced singularities in $\rho$ by 
the geometry, where $h$-adaptive FEM uses highly local mesh refinements 
around singularities, or (ii) log-normal coefficient fields,
where the stochastic scaling factors multiplying $\tol_\mathrm{bias}$
in the right hand sides of~\eqref{eq:error_model_stop} 
and~\eqref{eq:error_model_stop_uni} vary hugely.
In such cases, this simplified complexity analysis does not directly apply, 
but the complexity can be estimated by introducing a parametric
regularization whose effect goes to zero as a suitable power of 
$\tol_\mathrm{bias}$, see~Example~4.1 in~\cite{adFEM_our}, which is
related to Example~0, 1, and~2 in this work.


\subsubsection{A fully stochastic adaptive algorithm} 

This approach is justified if new local mesh refinements are 
	required for each sample; for example: (i) if the error density has 
a singularity in a random location, or (ii) if the stochastic coefficient 
field has a short correlation length relative to the domain size. 
The algorithm can be developed by combining estimates of 
$\int_\mathcal{D} \E{\rho^{\frac{d}{p+d}}}$ with the deterministic
goal-oriented adaptive algorithm, where motivated 
by~\eqref{eq:error_model_stop} the stopping condition for the
adaptive refinements is set to
\begin{align}
	\label{eq:error_model_stop_XXX}
	\int_\mathcal{D} \rho(x;\omega) h(x;\omega)^p\,dx & \leq
	\tol_\mathrm{bias} 
	\frac{\int_\mathcal{D} \rho(x;\omega)^{\frac{d}{p+d}}\,dx}{\int_\mathcal{D} \E{\rho^{\frac{d}{p+d}}}}
\end{align}
for each sample, $\omega$, using an unsigned error density $\rho$. 
Here, the left-hand side is approximated 
using the error estimate of the adaptive FEM algorithm. The stochastic
scaling factor
\begin{align}
	\label{eq:scaling_factor}
	\frac{\int_\mathcal{D} \rho(x;\omega)^{\frac{d}{p+d}}\,dx}{\int_\mathcal{D} \E{\rho^{\frac{d}{p+d}}}}
\end{align}
corresponds to optimally distributed error contributions from different
samples. The numerator in this scaling factor is computed as a by-product 
of computing the error estimate for each sample. The denominator can be 
estimated using an (ML)MC estimator. With a view toward 
the AMLMC setting, we observe that this can be iteratively performed
in a continuation type algorithm, as described for other MLMC parameters 
in~\cite{continuation_MLMC_our}.

Note that the stochastic scaling factor~\eqref{eq:scaling_factor}
is unbounded for a log-normal coefficient field 
in~\eqref{eq:bvp_general}. 
A standard approach in MLMC is to simply distribute 
the error contributions evenly over the samples, i.e. to replace the scaling 
factor~\eqref{eq:scaling_factor} by the deterministic constant 1. That
is far from optimal here.

The problem setting targeted by this work is not 
of a kind where 
the fully stochastic adaptive algorithm is optimal; 
at least not when the computational 
overhead of generating new computational meshes for each sample is 
taken into account. 
Because the focus is on situations where (i) a 
geometry-induced singularity in a deterministic location justifies highly 
localized mesh refinements, and (ii) the stochastic coefficient
field is smooth with long correlation length but high variability,
we expect optimal meshes for all samples to be well
approximated by meshes selected from a deterministic sequence of 
adaptively refined meshes, as described in Section~\ref{sec:proposed_MLMC}.
Therefore, avoiding the computational overhead of mesh generation during 
the sampling phase leads to substantial computational gains.


\subsection{Multilevel hierarchies of optimal meshes}
\label{sec:optimal_MLMC_hierarchies}

We consider MLMC estimators obtained by introducing 
a decreasing sequence of bias tolerances, 
\begin{align}
	\label{eq:tol_sequence}
	\tol_{\ell} & 
	=\tolbsep^\ell\tol_{0}, 
	&& \ell=0,1,\dots,
\end{align}
for a given separation constant, $0<\tolbsep<1$, using
the optimal single level approximations above.
This choice fixes three parameters in Assumption~\ref{assu:MLMC}: 
$\kappa=\tolbsep^{-1}$, $\gamma=d/p$, and $q_w=1$. (Observe that defining 
the levels by decreasing tolerances makes $q_w$ independent of 
the FEM order, $p$.)
The MLMC complexity~\eqref{eq:complexity_MLMC} fundamentally 
depends on the fourth parameter, $q_s$, which we now turn to.

Here, we consider only smooth stochastic coefficient fields,  
whose pathwise regularity allows us to use the error
expansion from~\cite{adFEM_our} to obtain $q_s=2q_w=2$.
However, 
problems with rough coefficients require 
a different approach to the one considered here; see for 
instance~\cite{Hall_Hoel} for the case of uniform 
discretizations. 

Now we want to estimate the multiplicative constants in our models of 
the work and the strong convergence in terms of (quasi-) norms of the
error density.
Any family of mesh size functions parametrized by $\tol_\ell$ on the form
\begin{align}
	\label{eq:mesh_size_fam}
	h_\ell(x;\omega) & = \tol_\ell^{1/p}f(x;\omega),
\end{align}
such as the optimal meshes in~\eqref{eq:opt_h_stoch} 
and~\eqref{eq:opt_h_uniform} with $\tol_\mathrm{bias}=\tol_\ell$, 
satisfies
\begin{align*}
	\QoI_\ell(\omega) - \QoI_{\ell-1}(\omega) 
	& \simeq \int_\mathcal{D}\rho(x;\omega)
	\left(h_{\ell-1}(x;\omega)^p-h_\ell(x;\omega)^p\right)dx \\
	& = \tol_{\ell}\left(\tolbsep^{-1}-1\right)
	\underbrace{\int_\mathcal{D}\rho(x;\omega)f^p(x;\omega)dx}_{K_1(\omega)}
\end{align*}
and
\begin{align*}
	\int_\mathcal{D} h_\ell^{-d}(x;\omega)dx & \simeq 
	\tol_\ell^{-d/p}
	\underbrace{\int_\mathcal{D} f^{-d}(x;\omega)dx,}_{K_2(\omega)}
\end{align*}
such that 
\begin{align}
	\label{eq:optimal_variance_mesh_size_fam}
	\var{\QoI_\ell - \QoI_{\ell-1}} & \simeq\tol_{\ell}^2\left(\tolbsep^{-1}-1\right)^2\var{K_1}
\end{align}
and
\begin{align}
	\label{eq:optimal_work_mesh_size_fam}
	\int_\mathcal{D} \E{h_\ell^{-d}} + \int_\mathcal{D} \E{h_{\ell-1}^{-d}} & \simeq 
	\tol_\ell^{-d/p} \left(1+\tolbsep^{d/p}\right)\E{K_2}.
\end{align}
Consequently, as long as $\var{K_1}$ and $\E{K_2}$ are bounded, 
Assumption~\ref{assu:MLMC} is satisfied with $\kappa=\tolbsep^{-1}$, $\gamma=d/p$, 
$q_s=2$, and $q_w=1$.  

We can estimate the problem-dependent multiplicative
constants in the complexities, by letting $V_\ell$ and $W_\ell$ be given 
by~\eqref{eq:optimal_variance_mesh_size_fam} 
and~\eqref{eq:optimal_work_mesh_size_fam}, respectively, 
in the summands with $\ell\geq 1$ in the sum 
$\sum_{\ell=0}^{L}\sqrt{W_\ell V_\ell}$, 
which appears in the asymptotically dominant first term of the MLMC work 
model~\eqref{eq:total_work_MLMC}. Then
\begin{multline*}
	V_\ell W_\ell =
	\var{\QoI_\ell - \QoI_{\ell-1}}
	\left(\int_\mathcal{D} \E{h_\ell^{-d}} + \int_\mathcal{D} \E{h_{\ell-1}^{-d}}\right) 
	\\ \simeq
	\tol_0^{2-d/p}
	\left(\tolbsep^{-1}-1\right)^2 \left(1+\tolbsep^{d/p}\right)
	\var{K_1}\E{K_2} \tolbsep^{\left(2-d/p\right)\ell},\quad\ell\geq 1.
\end{multline*}
Summing $\sqrt{V_\ell W_\ell}$ over $\ell=1,\dots,L$, considering
the three cases $d<2p$, $d=2p$, and $d>2p$ for the geometric sum 
$\sum_{\ell=1}^L \tolbsep^{\left(1-\frac{d}{2p}\right)\ell}$, and
letting $\tol\to 0$ leads to the statement below on the AMLMC
complexity. The statement in the corollary uses that for the two
cases in~\eqref{eq:error_model_stop} and~\eqref{eq:error_model_stop_uni}
\begin{align*}
	\var{K_1} & = \frac{\var{Y}}{\E{Y}^2} = \left(c_\mathrm{V}(Y)\right)^2,
\end{align*}
where
$c_\mathrm{V}(Y)$ denotes the coefficient of variation of $Y$, 
and with 
$Y=\Vert\rho\Vert_{L^{\frac{d}{p+d}}(\mathcal{D})}^{\frac{d}{p+d}}$
and 
$Y=\Vert\rho\Vert_{L^1(\mathcal{D})}^{\frac{d}{p+d}}$, respectively, for the two cases.
%
%
\begin{corollary}[Adaptive MLMC Complexity]\label{cor:AMLMCComplexity}
	Assume the work and bias models~\eqref{eq:nelem_model} 
	and~\eqref{eq:bias_model} hold 
	and the sampling error control satisfies either
	\begin{description}
	\item[(i)] \eqref{eq:error_model_stop}, for 
		finite $\left\Vert\rho\right\Vert_{L_\mathbb{P}^{\frac{d}{p+d}}(\mathcal{D}\times\Omega)}$, 
		(fully adaptive case), or
	\item[(ii)] \eqref{eq:error_model_stop_uni}, for 
		finite $\left\Vert\Vert\rho\Vert_{L^1(\mathcal{D})}\right\Vert_{L_\mathbb{P}^{\frac{d}{p+d}}(\Omega)}$,
		(adaptive selection of uniform meshes).
	\end{description}
	Then, 
	the MLMC computational work model~\eqref{eq:work_sum} 
	with $M_\ell$ given by~\eqref{eq:opt_nrs} 
	and $L\in\rset^+$ implied by~\eqref{eq:BIASC} and~\eqref{eq:tol_sequence} 
	satisfies
	\begin{align}
  	\label{eq:complexity_MLMC_adaptive}
  	\work_\mathrm{MLMC} & \simeq
  	\begin{cases}
  		K\tol^{-2}, & \text{if $d<2p$},\\
    	K\tol^{-2}\left(\log{\tol^{-1}}\right)^2, & \text{if $d=2p$},\\
    	K\tol^{-d/p}, & \text{if $d>2p$},
  	\end{cases}
	\end{align}
	with $K=K\left(\rho;\tolbsep,\tol_0,\splitting,\confpar\right) 
	= K_3(\rho) K_4(\tolbsep,\tol_0,\splitting) K_5(\tolbsep,\splitting,\confpar)$
	for
	\begin{align*}
		K_3 & = \E{K_2}\var{K_1} \\ 
		& =
  	\begin{cases}
  		\left\Vert\rho\right\Vert_{L_\mathbb{P}^{\frac{d}{p+d}}(\mathcal{D}\times\Omega)}^{d/p}
  		\left(c_\mathrm{V}\left(\Vert\rho\Vert_{L^{\frac{d}{p+d}}(\mathcal{D})}^{\frac{d}{p+d}}\right)\right)^2
  		, & \text{case (i)},\\
  		\int_\mathcal{D}1\,dx 
  		\left\Vert\Vert\rho\Vert_{L^1(\mathcal{D})}\right\Vert_{L_\mathbb{P}^{\frac{d}{p+d}}(\Omega)}^{d/p}
  		\left(c_\mathrm{V}\left(\Vert\rho\Vert_{L^1(\mathcal{D})}^{\frac{d}{p+d}}\right)\right)^2
    	, & \text{case (ii)},
  	\end{cases}
	\end{align*}
	\begin{align*}
		K_4 & =
  	\begin{cases}
  		\tol_0^{-\frac{d}{p}}
  		\left(
  			\sqrt{\frac{V_0}{\var{K_1}}}
  			\frac{1}{\left(\tolbsep^{-1}-1\right) \sqrt{1+\tolbsep^{d/p}}}
  			+
  			\tol_0\frac{C^{1-\frac{d}{2p}}}{1-C^{1-\frac{d}{2p}}}
  		\right)^2
  		, & \text{if $d<2p$},\\
    	\left(\log{C}\right)^{-2}, & \text{if $d=2p$},\\
    	\left(  1-\splitting \right)^{\left(2-\frac{d}{p}\right)}
    	\left(1-C^{\frac{d}{2p}-1}\right)^{-2}, & \text{if $d>2p$},
  	\end{cases}
	\end{align*}
	and
	\begin{align*}
		K_5 & = 
		\left( \frac{\confpar}{\splitting} \right)^2
		\left(\tolbsep^{-1}-1\right)^2 \left(1+\tolbsep^{d/p}\right).
	\end{align*}
\end{corollary}
%
%

The proof is given in~\ref{app:proof_corollary}.

\begin{remark}[On integer constraints on $L$]
	\label{rem:integer_L}
	The simplification of allowing non-integer values of $L$
	in Corollary~\ref{cor:AMLMCComplexity} does not affect the result 
	for $d<2p$ and $d=2p$. When $d>2p$, imposing $L=\lceil L^\ast\rceil$, 
	where $L^\ast\in\rset^+$ is the value implied by~\eqref{eq:BIASC} 
	and~\eqref{eq:tol_sequence}, leads to
	\begin{align*}
		K \leq \liminf_{\tol\to 0^+} \frac{\work_\mathrm{MLMC}}{\tol^{-d/p}}
		& \leq \limsup_{\tol\to 0^+} \frac{\work_\mathrm{MLMC}}{\tol^{-d/p}}
		\leq K \tolbsep^{1-\frac{d}{2p}},
	\end{align*} 
	with $K$ defined in the corollary.
\end{remark}

\begin{remark}[On the use of the leading order approximation of bias]
	\label{rem:bias_model}
	Corollary~\ref{cor:AMLMCComplexity} is stated under the explicit
	assumption that the bias model~\eqref{eq:bias_model} holds, and it is used
	on all levels $\tol_\ell$, $\ell=0,1,\dots,L$. This model, however, 
	is only the leading order term of the bias as $\tol_\ell\to 0$.
	Theorem~\ref{thm:MLMC_compl} allows for letting $\tol_0\to 0$, slowly, 
	as $\tol\to 0$, which strengthens the connection to the
	asymptotic theory. Doing this comes at the cost of a logarithmic factor 
	in the complexity estimate when $d<2p$, but not when $d\geq 2p$. 

	A sharper analysis, relaxing~\eqref{eq:bias_model} and taking the effect 
	of the non-asymptotic levels in the complexity estimate into account,
	can be developed for this setting closely following the analysis 
	in~\cite{adaptive_MLMC_2_our}.
	We do not do it here for the sake of length of exposition.
\end{remark}

\begin{remark}[On the effect of $\tolbsep$]
	\label{rem:tolsep}
	In practical computations, it is unrealistic to assume 
	that $\var{\QoI_\ell - \QoI_{\ell-1}}\to 0$ as $C\to 1$. 
	The model~\eqref{eq:optimal_variance_mesh_size_fam}, with the
	factor $\left(\tolbsep^{-1}-1\right)$ is relevant when 
	$\tolbsep$ is far from 1; we use $\tolbsep=1/4$ in the numerical
	experiments in Section~\ref{sec:numex}. 
	For optimization of 
	$\tolbsep$, we can instead use the bound 
	$\var{\QoI_\ell - \QoI_{\ell-1}}\leq 
	2\var{\QoI - \QoI_\ell}+2\var{\QoI - \QoI_{\ell-1}}$, which leads
	to replacing $\left(\tolbsep^{-1}-1\right)$ by 
	$2\left(\tolbsep^2+1\right)$. This can be seen as a penalization
	on closely spaced tolerances.
\end{remark}

\begin{remark}[On the effect of $\tol_0$]
	When $d\geq 2p$, the choice of $\tol_0$ does not have any effect on the 
	asymptotic estimates in Corollary~\ref{cor:AMLMCComplexity}, 
	as long as the choice is fixed with respect to $\tol$.
	However, when $d<2p$ the choice of $\tol_0$ matters even
	asymptotically as $\tol\to 0$. The factor $K_4(C,\tol_0,\splitting)$
	allows for an explicit minimization with respect to $\tol_0$ to identify
	the unique minimizer in $(0,\infty)$.
\end{remark}

The constants in the corollary above were obtained under 
idealized assumptions, some of which were expanded
upon in Remarks~\ref{rem:integer_L},~\ref{rem:bias_model}, 
and~\ref{rem:tolsep}.
Considering these reservations, our particular numerical experiments in 
Section~\ref{sec:numex} show good agreement between the modeled 
and experimentally observed values of $V_\ell$ (Table~\ref{tab:ratio_variance_emperical}).

\begin{remark}[Dirichlet--Neumann boundary
	singularity (Needle problem) ]
	
	The gains of AMLMC with respect to SMLMC are dimensional dependent.
	To understand the potential advantages of AMLMC with respect to SMLMC, consider the problem described in Section~\ref{sec:example_problems}.
	This problem can be generalized to higher dimensions. For instance, to obtain the corresponding formulation in dimension $d=3,$ one uses a $3D$ cylindrical domain generated by rotating the  $2D$ domain depicted in Figure~\ref{fig:problem_settings}.a around its top boundary. In that case, the Neumann boundary condition is dropped, while the Dirichlet boundary condition corresponding to the axis of rotation, the line segment joining $(0,0)$ with $(1,0)$, is maintained, thus 
	yielding the penetrating needle shape. In this problem, all boundary conditions are zero Dirichlet. By its construction, this problem inherits the same singularity as the $2D$ problem and now the relevant distance to the singularity is the distance to the needle tip. 
	Here, one can argue that the error density quasi-norm is uniformly bounded and for AMLMC we have the full FEM order $p=2$. Thus $d<2p$ so that, assuming we use an optimal solver, the AMLMC yields $\tol^{-2}$ complexity while the SMLMC yields the worse $\tol^{-3}$ complexity in that case. The situation is even more favorable in $d=4$, where AMLMC and SMLMC yield $\tol^{-2}\log^2(\tol^{-1})$ and $\tol^{-4}$, respectively.
\end{remark}


\subsection{A Multilevel Monte Carlo estimator tailored to deterministic adaptivity}
\label{sec:proposed_MLMC}

In this section we demonstrated that AMLMC, guided by
an error density $\rho$, can provide substantial computational gain
compared with the standard MLMC on uniform discretizations, particularly
for problems where the norm 
$\left\Vert\rho\right\Vert_{L_\mathbb{P}^{1}(\mathcal{D}\times\Omega)}$
is not uniformly bounded over refinement sequences, while 
$\left\Vert\rho\right\Vert_{L_\mathbb{P}^{\frac{d}{p+d}}(\mathcal{D}\times\Omega)}$
or
$\left\Vert\Vert\rho\Vert_{L^1(\mathcal{D})}\right\Vert_{L_\mathbb{P}^{\frac{d}{p+d}}(\Omega)}$
do remain uniformly bounded; see case~(i) and~(ii) in Corollary~\ref{cor:AMLMCComplexity}.

For problems where highly localized adaptive mesh refinements are 
justified by deterministic effects, such as the needle problem with a lognormal coefficient field, 
we propose an approximation to the fully adaptive case~(i) using 
adaptive selection of graded meshes, generated by the $h$-adaptive 
FEM algorithm applied to a deterministic error density,  
for example the error density obtained when the random 
coefficient field is replaced by its expected value. Then we
use~\eqref{eq:error_model_stop} to select the appropriate deterministic
mesh for each outcome.

\subsubsection{Stochastic acceptance of deterministically generated meshes}

The idea in its single-level form is to divide the algorithm into
two parts:\\
(1) Using a decreasing sequence of tolerances, 
$\tol_{\mathrm{det},k}=C_1^k\tol_{\mathrm{det},0}$ for 
some 
$C_1\in(0,1)$, 
we generate FEM discretizations with mesh size functions 
$h_k:\mathcal{D}\to\rset$, $k=0,1,\dots$, using the goal-oriented 
FEM algorithm applied to the deterministic problem with 
\emph{unsigned} approximations to the error density.\\
(2) Given a bias tolerance, $\tol_\mathrm{bias}$, 
and an estimate $R\approx\int_\mathcal{D} \E{\rho^{\frac{d}{p+d}}}$,
we generate approximate samples of $\QoI$ by starting on the coarsest
level, $k=0$, computing the primal and dual solutions using mesh $h_k$
and the corresponding approximate \emph{signed} error density 
$\tilde\rho_k(x;\omega)$. Repeat for increasing $k$ until 
\begin{align}
	\label{eq:model_stopping}
	\left|\int_\mathcal{D} \tilde\rho_k(x;\omega) h_k(x)^p\,dx\right| & \leq
	\frac{\tol_\mathrm{bias}}{R}
	\left|\int_\mathcal{D} \tilde\rho_k(x;\omega)^{\frac{d}{p+d}}\,dx\right|.
\end{align}

The MLMC version, defined as before by a sequence 
(possibly distinct from $\tol_{\mathrm{det},k}$) of 
tolerances $\tol_\ell=\tolbsep^\ell \tol_0$, generates sample pairs
$\left(\QoI_\ell,\QoI_{\ell-1}\right)$ by applying the procedure in 
part~(2) to the bias tolerances $\left(\tol_\ell,\tol_{\ell-1}\right)$.

Note that the goal-oriented FEM algorithm only has to be called
once for each mesh $h_k$, which is then saved for all future samples. 
This sequence of meshes can be accessed on demand by generating 
additional refined meshes only when~\eqref{eq:model_stopping} fails to hold
for all previously generated meshes. Therefore, in practice 
parts (1) and (2) interlace.

This approach is particularly useful in the lognormal case, 
cf. Section~\ref{sec:numerical_results}, 
where we must use a stochastic mesh to efficiently address the lack of 
uniform coercivity.

\section{Implementation and numerical examples}
\label{sec:numex}

This section presents the algorithms and their implementation for SMLMC and AMLMC{; moreover, it} shows a numerical comparison between them for elliptic problems with random coefficients. {Here, we} consider both {lognormal constant and a spatially-correlated lognormal fields (Section~\ref{sec:example_1}} and Section~\ref{sec:example_2}), in a {2D} physical domain with a geometrically driven singularity. 

\subsection{Implementation}
\label{sec:implement}

{First, we present} the implementation of the numerical approximation of the elliptic problem for a fixed realization of the random coefficient field, followed by SMLMC and AMLMC {algorithms} and their implementation.

The numerical experiments are written in \texttt{C++} and run on a workstation with {a} 2.10GHz Intel Xeon Gold 6230R CPU and 252 GByte RAM. The computational work is measured in FLOPS using {the} performance application programming interface (PAPI), Version 6.0.0~\cite{PAPI_paper}. The FEM implementation uses the \texttt{deal.II} Library, Version 9.2~\cite{dealII92}.

\begin{description}

	\item[Compute the Primal and Dual Solution] The primal and dual solutions {$u$ and $\varphi$, respectively,} to problem~\eqref{eq:bvp_general} are computed on the same quadrilateral mesh (possibly with hanging nodes, handled by {the} \texttt{deal.II} {library}~\cite{dealII92}), using the same first-order bilinear quadrilateral finite element class $Q_1$, see~\cite{hughes2012finite}, for instance. 
	
	\item[The Iterative Solver Control] {The primal and dual solutions obtained by an iterative solver are denoted by $u_{\textrm{iter}}$ and $\varphi_{\textrm{iter}}$, respectively}. The algorithm iteratively proceeds, {thus simultaneously} advancing the pair of primal solvers and dual solvers in each iteration step, {checking} against the solver stopping tolerance, see Algorithm~\ref{alg:primal_dual_solver}. For each sample,  we {select} the iterative solver tolerance $\textrm{TOL}_{\textrm{iter}}$ to be $\frac{1}{10}$ of the stopping tolerance {such} that the iterative solver error does not dominate the error. More details can be seen in Algorithm~\ref{alg:sample_DQ}, where we use iterative solvers to compute samples in AMLMC. {Here,} the preconditioned conjugate gradient method {is chosen} as our iterative solver. 
	
	\item[Compute the Error Density] The error density $\rho$ is computed using the finite element solutions $u_h$ and $\varphi_h$. An approximation of the error density, on a quadrilateral cell $K$, $\tilde{\rho}_K$, is derived in \cite{adFEM_our} and written as
	\begin{align}
	\tilde{\rho}_K := \frac{1}{48} \sum_{j=1}^4 (a_{11}^\ast \overline{D_1^2 u_h} \, \overline{D_1^2 \varphi_h} + a_{22}^\ast \overline{D_2^2 u_h} \, \overline{D_2^2 \varphi_h}) (x_j^K),
	\label{eq:rho_tilde}
	\end{align}
	where $j = 1,2,3,4$ are the four vertices of the cell $K$, $D_1^2 u_h$ is the second-order difference quotient over the reference direction $x_1$ and the notation $\overline{D_1^2 w}$ means the averaged difference quotient $D_1^2 w$. This averaging is {required} to make the difference quotient uniformly converge to the corresponding derivatives. For discussions on its motivation and fast computation, we refer to Remark 2.1 and Section 2.3 of \cite{adFEM_our}. {From} our implemented data structure (see Algorithm~\ref{alg:tree_lines}), the complexities of computing difference quotients and error densities are {provided} in Table~\ref{tab:numerical_complexities_each_part}, which improved {on} the complexities in~\cite{adFEM_our}. 

	We adopt{ed} the strategy proposed in~\cite{adFEM_our}, {which} uses a lower bound for the error density, i.e.
	\begin{align}
		\bar{\rho} = \textrm{sgn}\left( \tilde{\rho}  \right) \max\{\lvert \tilde{\rho} \rvert, \delta\}
			\label{eq:approximate_error_density}
	\end{align}
	where $\textrm{sgn}(\cdot)$ is the sign function, $\lvert \tilde{\rho} \rvert$ is the unsigned computed error density from~\eqref{eq:approximate_error_density}, and $\delta = \frac{\lVert \tilde\rho \rVert_{L^{\frac{1}{2}}(\mathcal{D})}}{(\int_{\mathcal D} dx)^2} \sqrt{\textrm{TOL}}$, with $\textrm{TOL}$ being the tolerance used when generating the meshes (refer to Algorithm~\ref{alg:mesh_generation}). The lower bound {ensures} that the maximum element size {asymptotically} converges to 0, and the ratio of the error density at two consecutive refinement levels is close to 1~\cite{adFEM_our}. {Because} the adaptive algorithm does not include a coarsening step, we impose a $\tol$-dependent upper bound on the approximate error density to prevent over-refinement {because of} poor estimates on coarse meshes, cf. \cite{adFEM_our}.

	The \texttt{deal.II} library uses a tree structure to represent finite element meshes. To efficiently compute the difference quotients and subsequent error densities and error estimates, we implemented Algorithm~\ref{alg:tree_lines}, which extracts horizontal and vertical lines structures from a tree structure. 

	\item[SMLMC and AMLMC] SMLMC uses sequentially refined uniform meshes on the levels and uses a direct solver UMFPACK~\cite{umfpack_davis2004algorithm} for linear systems. AMLMC {comprises} two phases. In phase I, an adaptive mesh hierarchy is generated~(Algorithm~\ref{alg:mesh_generation}). In phase II, the MLMC sampling is performed on such meshes, to compute the MLMC estimate~\eqref{eq:MLMC_estimator}. The same parameters $\theta$ and $C_\xi$ {are used} in all AMLMC and SMLMC computations, with $\theta = 0.5$ and $C_\xi = 1.96$. 

	{Next, a} sample $(Q_\ell - Q_{\ell - 1})(\omega)$ is generated {using} Algorithm~\ref{alg:sample_DQ}, and the direct solver UMFPACK is applied in the coarsest mesh in AMLMC. The quantity $\E{\rho^{\frac{d}{p+d}}}$ in~\eqref{eq:model_stopping} is roughly approximated using the coarsest mesh in the hierarchy (Figure~\ref{fig:Uniform_Adaptive_Meshes}). {Note that} this estimate can be improved using an MLMC estimator. 

\end{description}

The complexities {of} the aforementioned implementation methods are summarized in Table~\ref{tab:numerical_complexities_each_part}. 

\begin{table}[ht]
	\centering
	\begin{tabular}{l|l}
		\hline
		& Complexity\\
		\hline
		Assembly & $\mathcal{O}(N^2)$ \\
		Direct Solver & $\mathcal{O}(N^3)$\\
		Iterative Solver & $\mathcal{O}(I_{\mathrm{iter}}N)$\\
		Assemble Line(Uniform) & $\mathcal{O}(N\log N)$\\
		
		Difference Quotient(Uniform) & $\mathcal{O} (N)$\\
		Error Estimation (Uniform) & $\mathcal{O}(N \log^2 \sqrt{N})$\\
		Assemble Line(Adaptive) & $\mathcal{O}(N^2)$\\
		Difference Quotient(Adaptive) & $\mathcal{O} (N)$\\
		Error Estimation (Adaptive) & $\mathcal{O} (N \log^2 N)$\\
		\hline
	\end{tabular}
	\caption{Numerical complexities of each simulation part. {Note that} the complexity for the iterative solver depends on the number of iterations, $I_{\mathrm{iter}}$. The algorithm for {the} Assemble Line is Algorithm~\ref{alg:tree_lines}. The provided complexity for {the} adaptive mesh is considered for the {worst}-case scenario. The difference quotient and error estimation follows a natural implementation {in} the data structure created by Algorithm~\ref{alg:tree_lines}. } 
	\label{tab:numerical_complexities_each_part}
\end{table}

\subsection{Example problems}\label{sec:example_problems}

\begin{figure}[htbp]
	\begin{subfigure}[b]{0.45\textwidth}
		\centering
		\includegraphics[width=\textwidth]{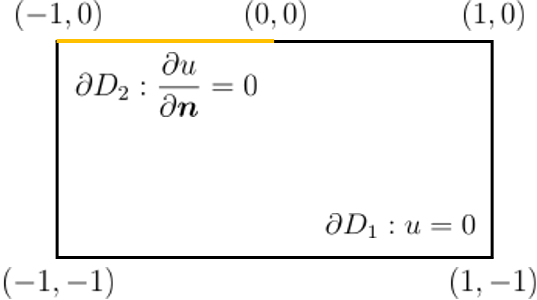}
		\caption{Our problem domain $[-1, 1] \times [-1, 0]$, with Dirichlet boundary condition on $\partial D_1$ and Neumann boundary condition on $\partial D_2$.}
		\label{fig:problem_domain}
	\end{subfigure}
	\hfill
	\begin{subfigure}[b]{0.4\textwidth}
		\centering
		\includegraphics[width=\textwidth]{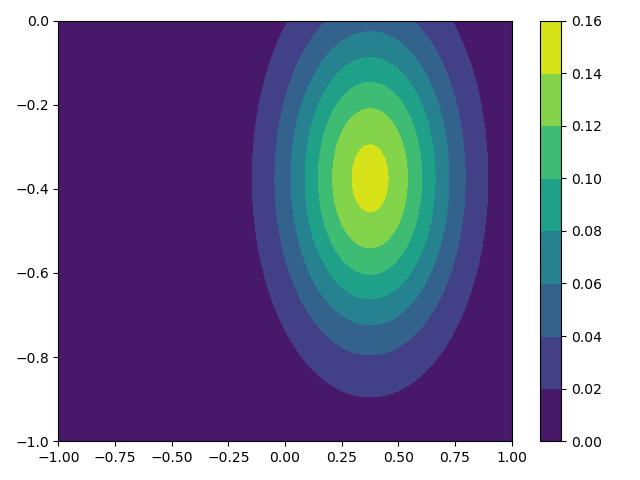}
		\caption{{Right-hand side of the dual problem}: indicator function $\mathbf{1}_{[0.25, 0.5]\times [-0.5, -0.25]}$ {convolved with the 2D} Gaussian kernel \eqref{eq:gaussian_kernel}.}
		\label{fig:dual_rhs}
	\end{subfigure}
	
	\caption{{The figure presents: (a) The problem domain and boundary conditions, and (b), a plot of the dual problem right-hand side function.}}
	\label{fig:problem_settings}
\end{figure}

We solve BVP~\eqref{eq:bvp_general} on the domain $\mathcal D = [-1, 1] \times [-1, 0] \subset \mathbb{R}^2$ (see Figure~\ref{fig:problem_domain}), with $\partial \mathcal D_1 = \partial \mathcal D - \partial \mathcal{D}_2$ and $\partial \mathcal{D}_2 = [-1, 0] \times (0)$, {for} $f \equiv 1000$.  Studies on the effect of the singularity {because of} the geometry on the Laplace solution regularity \cite{Stephan_Whiteman_1988,Li_Lu_2000,grisvard1992singularities, Stephan_Whiteman_1988, wahlbin1984sharpness, Johnson_book_1987,Babuska_book_2010} {show} that the regularity of the solution $u$ {reduces} to $H^{\frac{3}{2} - \epsilon}(\mathcal D)$ for $\epsilon > 0$, and the pointwise finite element approximation error for uniform meshes cannot be better than $\mathcal{O}(C(\omega) hr^{\frac{1}{2}})$, where $r$ is the distance to the origin, where the singularity occurs.  

In this example, the QoI is the weighted average
$Q(u) = (g * \mathbf{1}_{D_0}, u)$, 
where $*$ is the convolution operator, $\mathcal{D}_0 = \left[0.25, 0.5\right] \times \left[-0.5, -0.25\right]\subset \mathcal{D}$, and $g(\cdot)$ is the Gaussian density in $\mathbb{R}^2$ {given by}
\begin{equation}
	g(x) = \frac{1}{\sqrt{(2\pi)^2 |\Sigma|}} \exp(-\frac{1}{2}(x - \mu)^T \Sigma^{-1} (x - \mu) )
	\label{eq:gaussian_kernel}
\end{equation}
with mean $\mu = 0$ and covariance $\Sigma = \frac{1}{16} \bm{I}$.
The corresponding dual problem is~\eqref{eq:bvp_general}, with $f = \mathbf{1}_{D_0} * g$.

{Figure~\ref{fig:dual_rhs} gives} a visualization of {the right-hand side of the} dual problem. The different numerical examples below differ {in} the choice of diffusivity coefficient, $a$.

\subsection{Numerical results}\label{sec:numerical_results}

\subsubsection{Example 0}
\label{sec:example_0}
\begin{figure*}
	\centering
	{
	\renewcommand{\arraystretch}{1.2}
	\begin{tabular}{cM{18mm}M{18mm}M{18mm}M{18mm}}
		\toprule
		$k$ & 0 & 1 & 2 & 3\\
		\midrule
		UNIFORM &
		\includegraphics[width=0.16\textwidth]{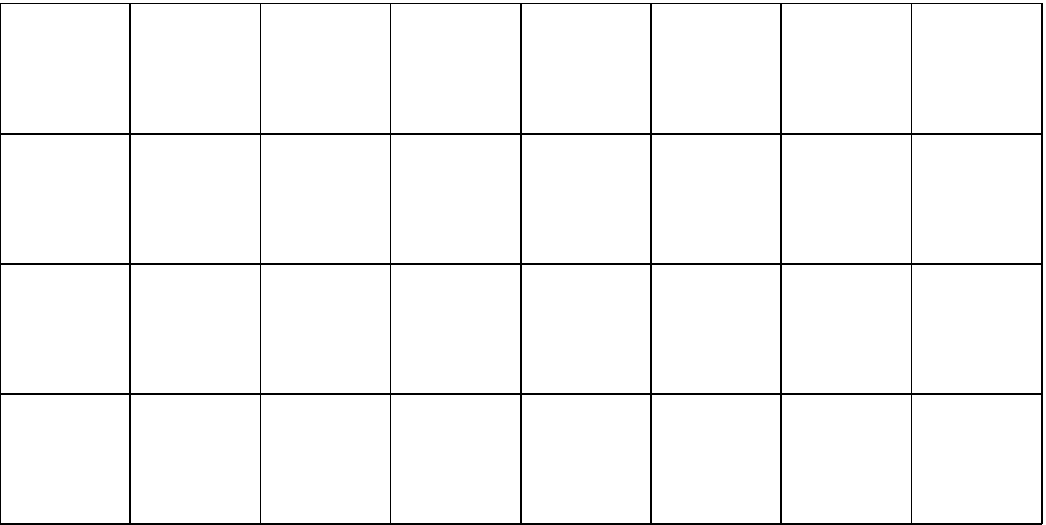}&
		\includegraphics[width=0.16\textwidth]{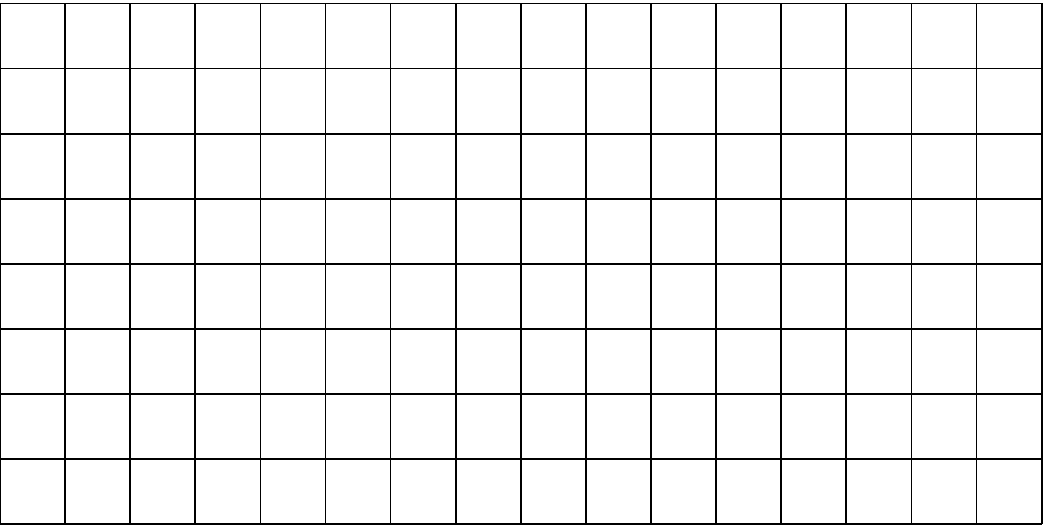}&	
		\includegraphics[width=0.16\textwidth]{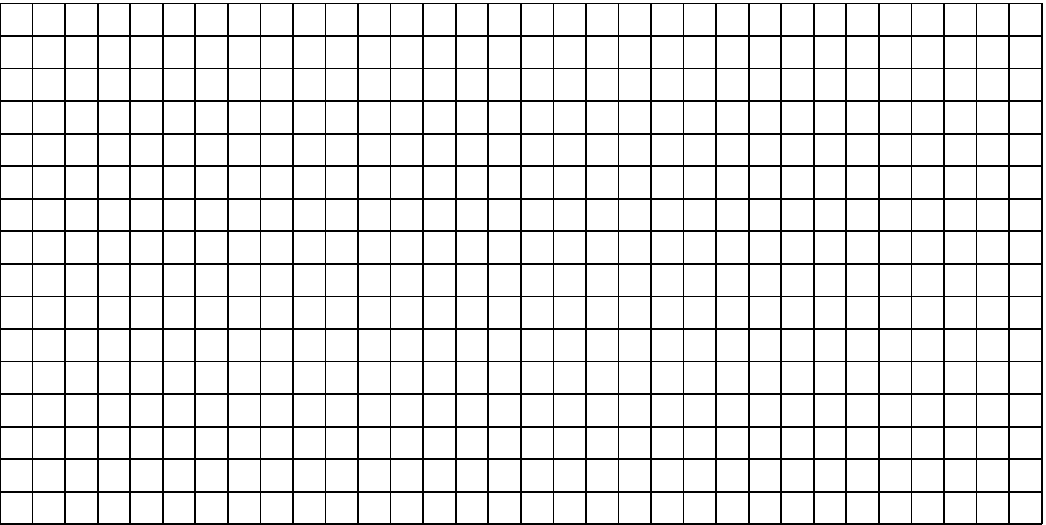}&
		\includegraphics[width=0.16\textwidth]{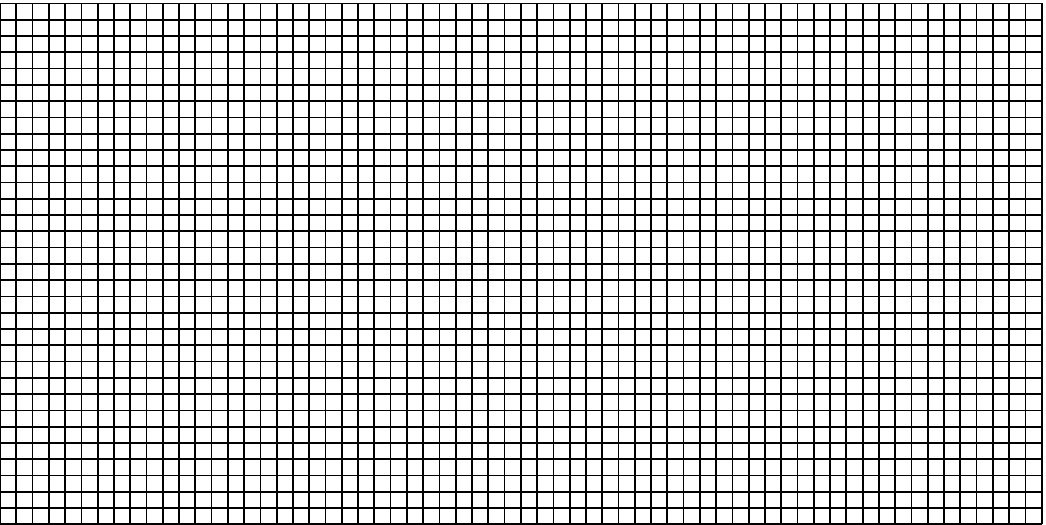}\\
		ADAPTIVE &
		\includegraphics[width=0.16\textwidth]{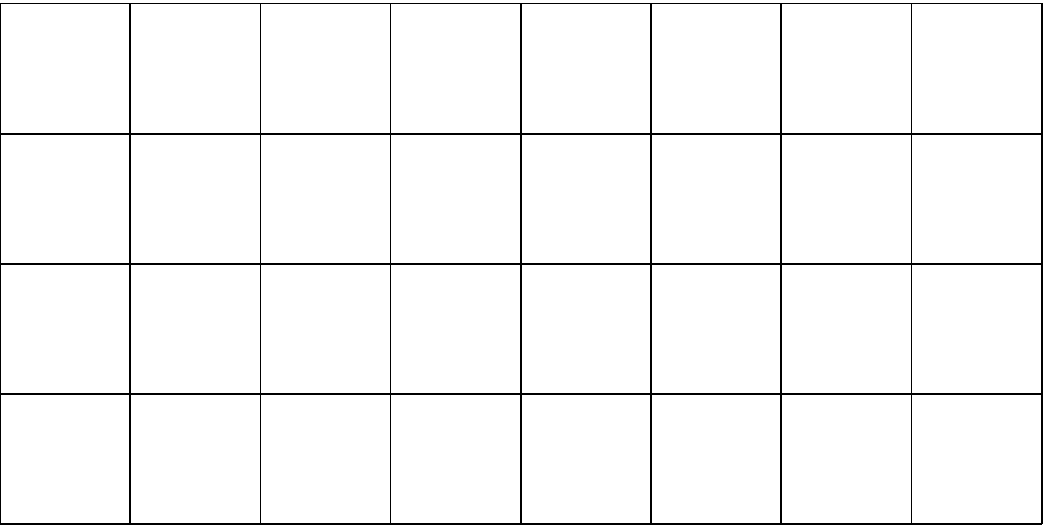}&
		\includegraphics[width=0.16\textwidth]{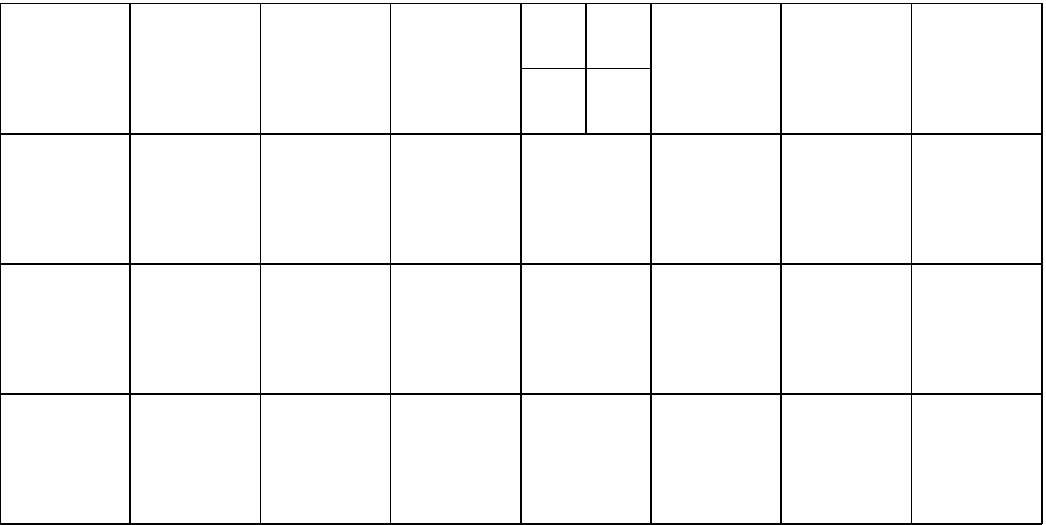}&
		\includegraphics[width=0.16\textwidth]{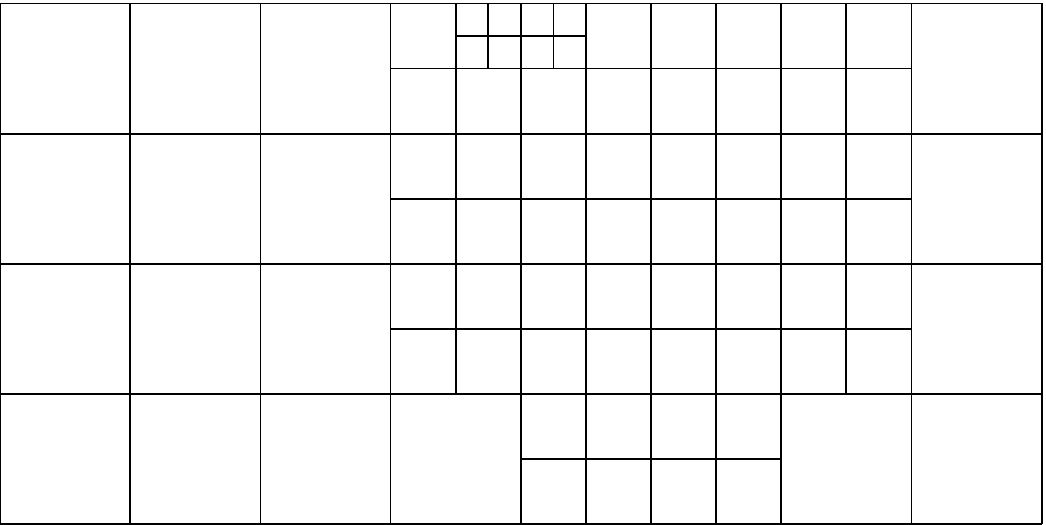}&
		\includegraphics[width=0.16\textwidth]{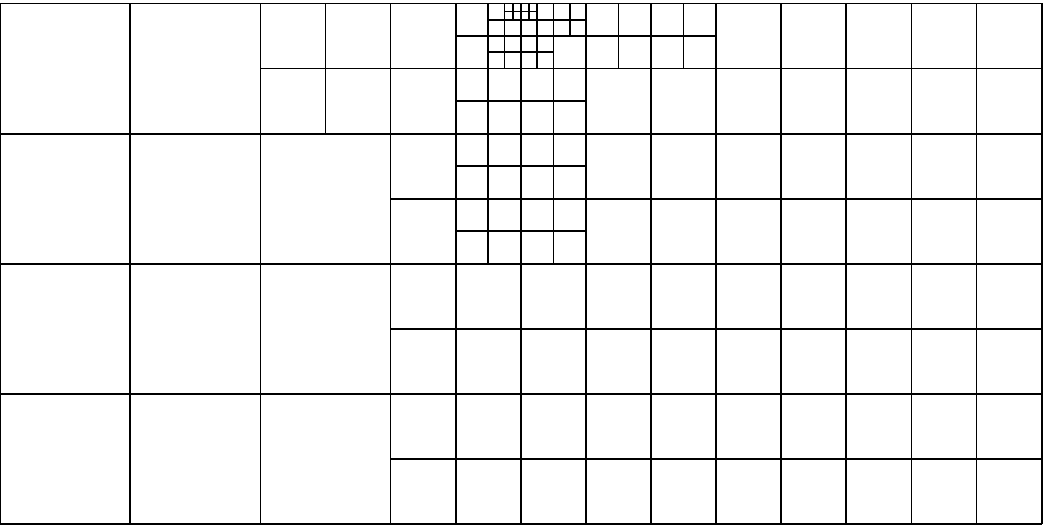}\\
		\toprule
		$k$ & 4 & 5 & 6 & 7\\
		\midrule
		ADAPTIVE &
		\includegraphics[width=0.16\textwidth]{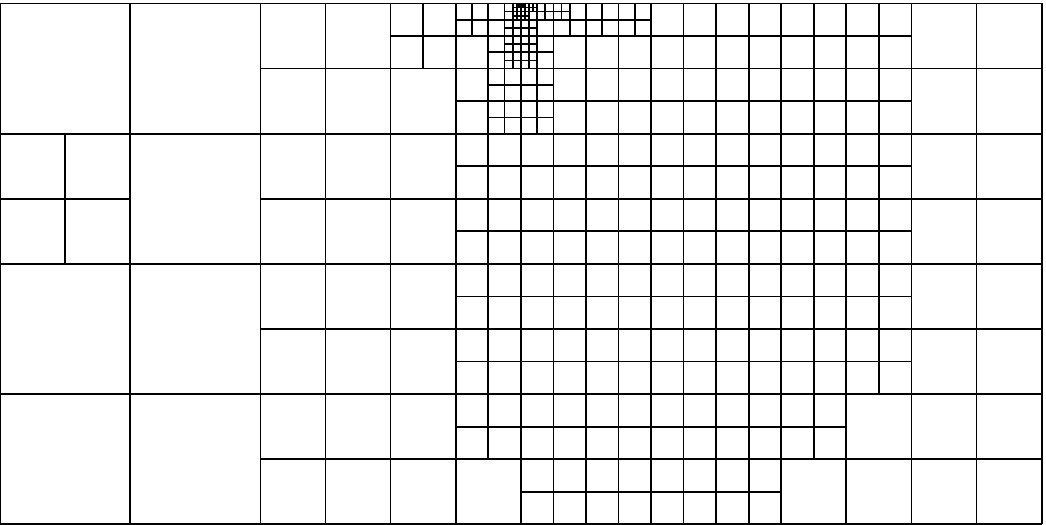}&
		\includegraphics[width=0.16\textwidth]{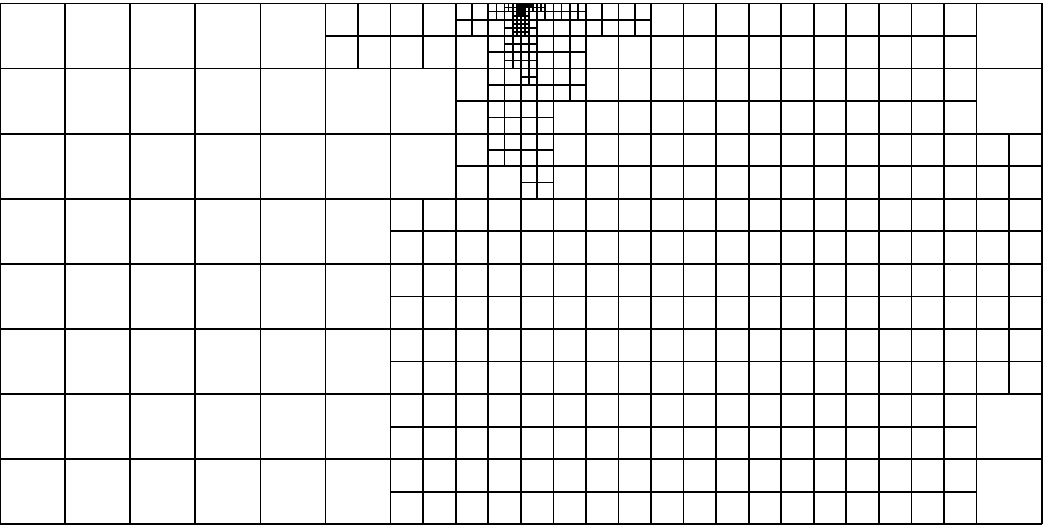}&
		\includegraphics[width=0.16\textwidth]{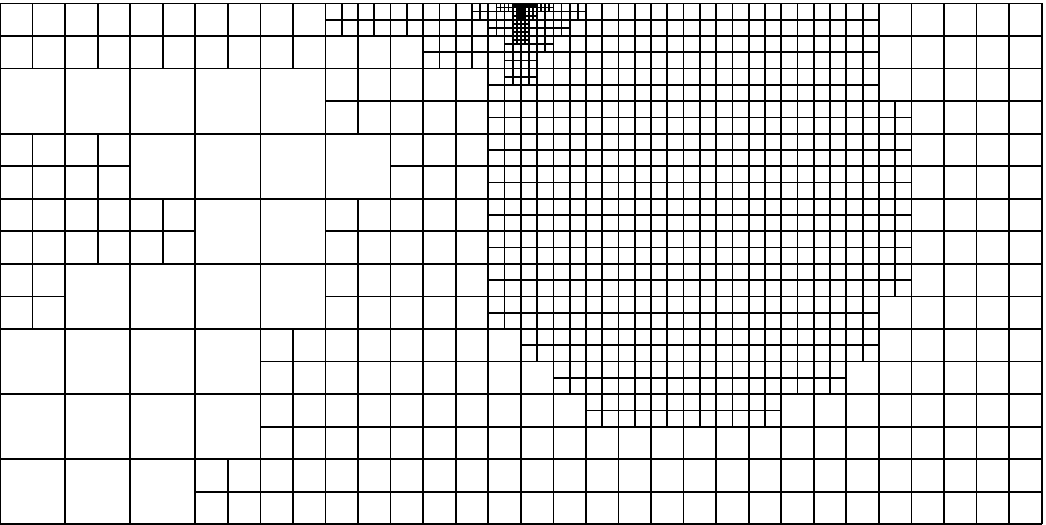}&
		\includegraphics[width=0.16\textwidth]{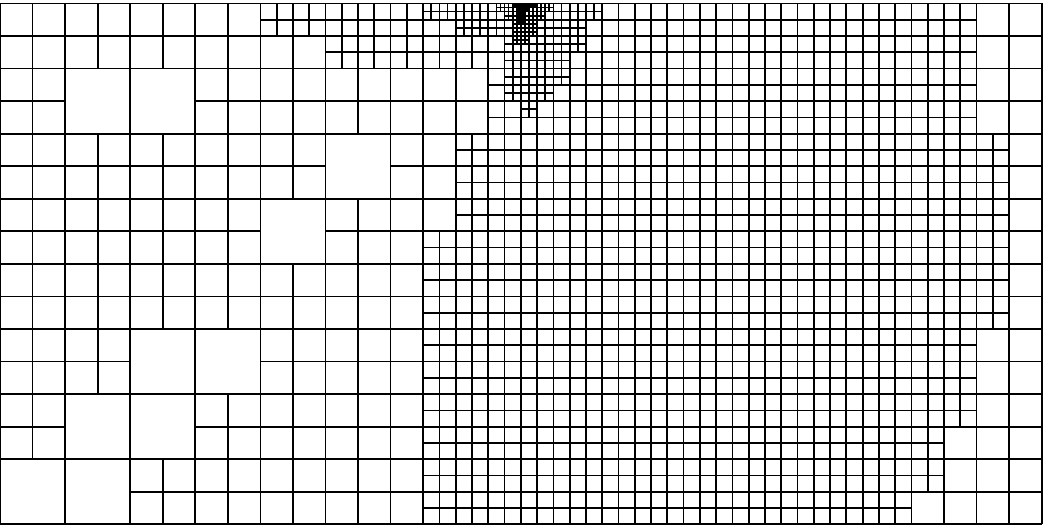}\\
		\bottomrule
	\end{tabular}
	}
	\caption{Example 0: The uniform and adaptive mesh hierarchies. Both mesh hierarchies start from the same mesh level $0$. As expected, the adaptive meshes are refined towards the singularity from the boundary conditions at $(0, 0)$ and {the region most affecting} the QoI. The adaptive meshes are generated {using} Algorithm~\ref{alg:mesh_generation} with $\tol_k = 2^{-(k+5)}$ and the {remaining} the parameters mentioned in Section~\ref{sec:example_0}. }
	\label{fig:Uniform_Adaptive_Meshes}
\end{figure*}

{Here, we treat} problem~\eqref{eq:bvp_general} with a deterministic, constant coefficient $a\equiv\exp{(2)}$. 
We consider the problem {for} uniform and adaptive meshes. The sequence of uniform meshes are globally refined once in each direction {per} time. The hierarchy of adaptive meshes is generated {using} Algorithm~\ref{alg:mesh_generation} with the sequence of $\tol$s being $\tol_k = 2^{-(k+5)}$, {$k=0,1,\ldots,$} {at} $C_R = 2.5$, $C_S = 3$, and $c = 2$. Figure~\ref{fig:Uniform_Adaptive_Meshes} shows the generated uniform and adaptive meshes. 

\begin{figure}[htbp]
	\centering
	\includegraphics[width=0.90\textwidth]{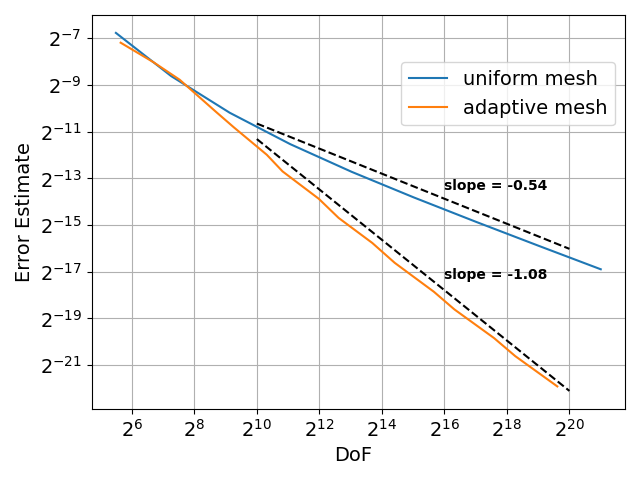}
	\\
	\includegraphics[width=0.48\textwidth]{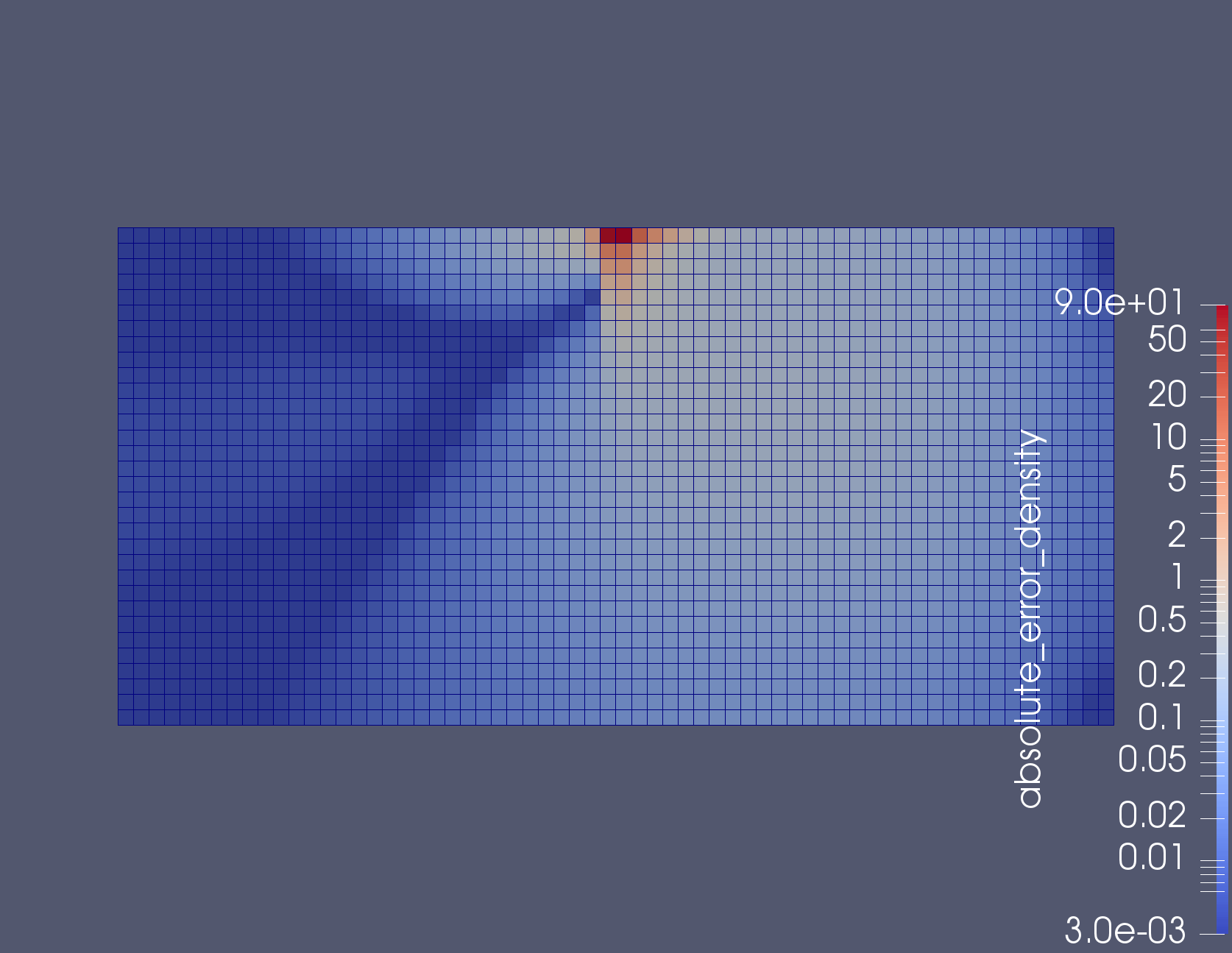}
	\hfill
	\includegraphics[width=0.48\textwidth]{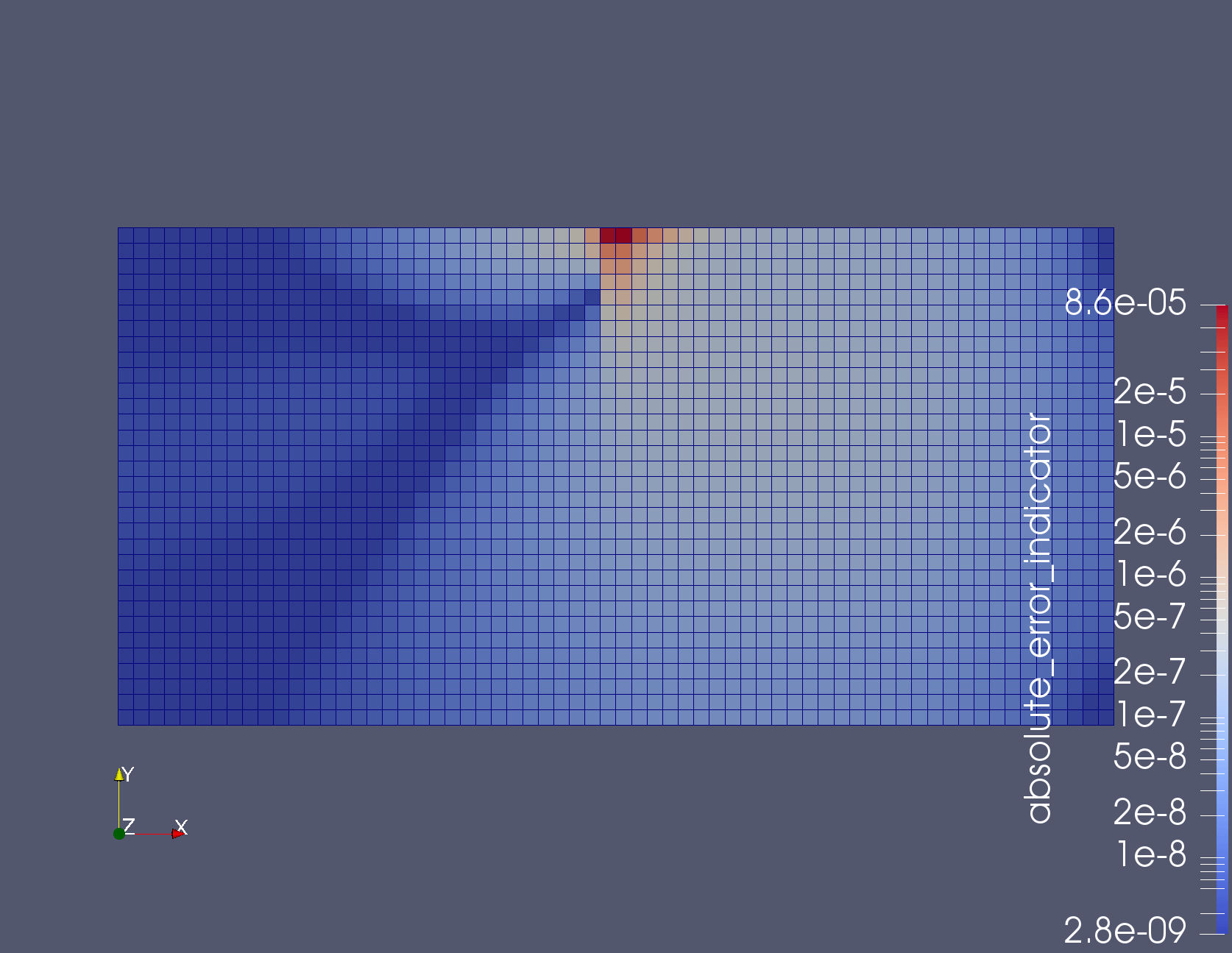}
	\\
	\includegraphics[width=0.48\textwidth]{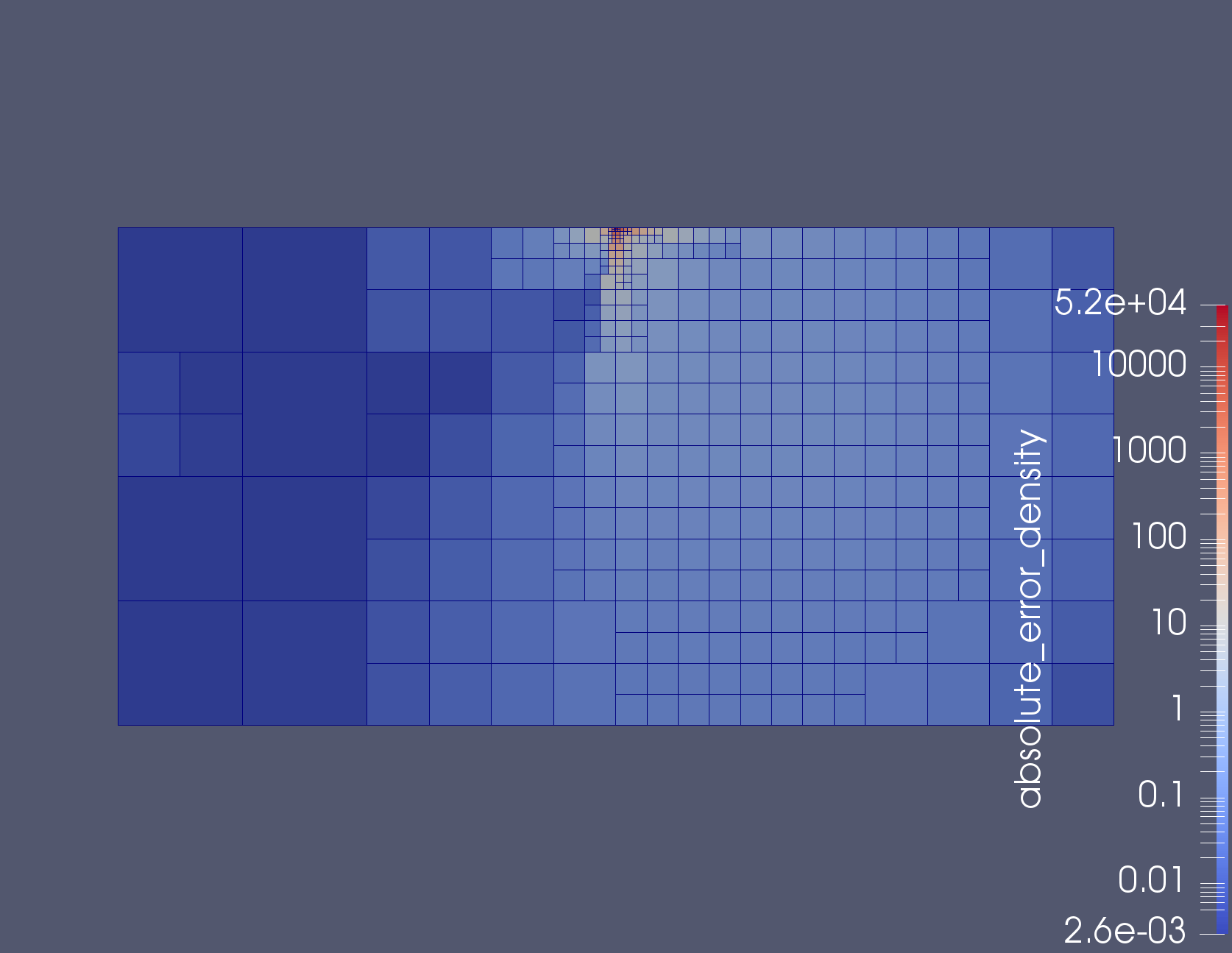}
	\hfill
	\includegraphics[width=0.48\textwidth]{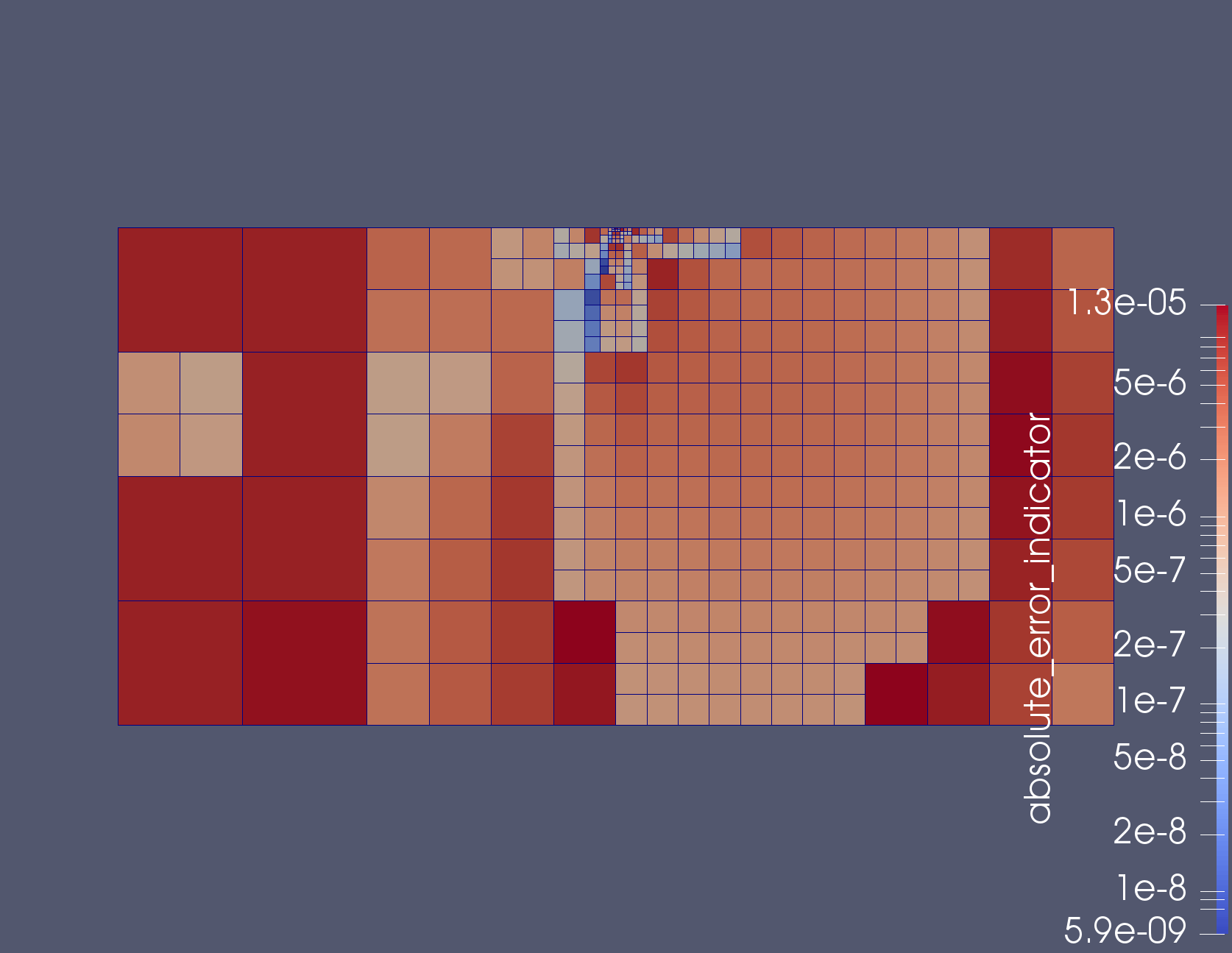}
	\caption{Example 0: The absolute error estimates $e_{\textnormal{est, err}}$~\eqref{eq:abs_error_estimate}, absolute error densities $\lvert \bar{\rho} \rvert$ from~\eqref{eq:approximate_error_density}, and absolute error indicators $\lvert r \rvert$ from~\eqref{eq:error_indicator} against \#DoF for the uniform and adaptive mesh hierarchies (Figure~\ref{fig:Uniform_Adaptive_Meshes}). Top: The convergence of the error estimate against \#DoF for uniform and adaptive meshes. Middle: The error density (left) and the absolute error indicator (right) for one uniform mesh. Bottom: The absolute error density (left) and absolute error indicator (right) for one adaptive mesh. }
	\label{fig:uniform_adaptive_error_estimate_DoFs}
\end{figure}

Figure~\ref{fig:uniform_adaptive_error_estimate_DoFs} shows the decay of the error estimate against {the number of DoF (\#DoF)} for uniform and adaptive mesh hierarchies. 
{ The observed convergence rate of the error estimate using adaptive meshes is twice the rate obtained with uniform meshes.}
{Adaptively refined meshes are more resolved around the singularity and the region most affecting the QoI, 
thus distributing the error contributions more equally over the cells.}

\begin{figure}[htbp]
	\centering
	\includegraphics[width=0.48\textwidth]{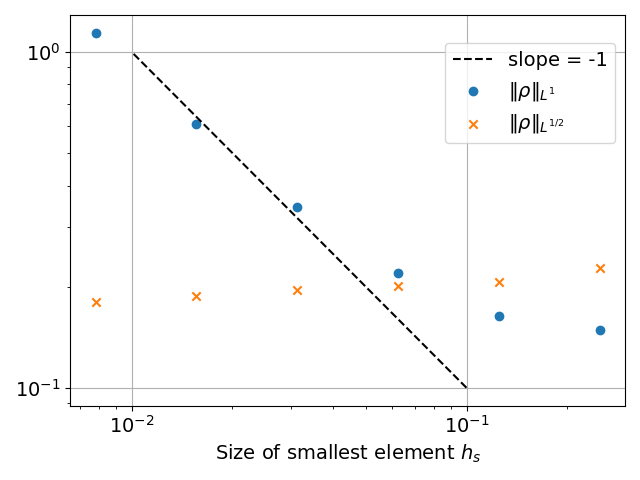}
	\hfill
	\includegraphics[width=0.48\textwidth]{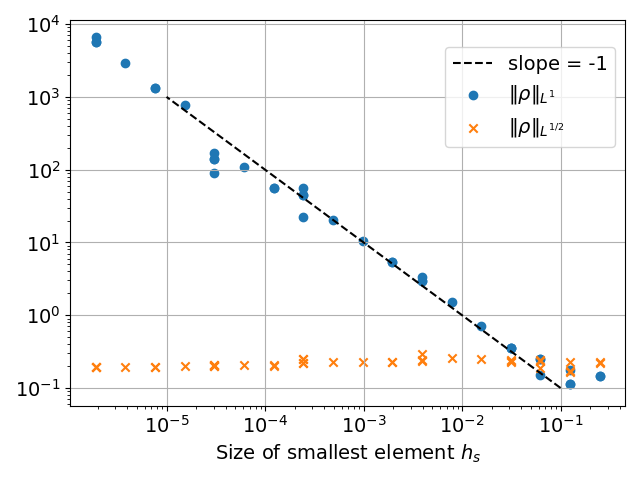}
	\\
	\caption{Example 0: {Graphs of} $L^1(\mathcal{D})$ norm and $L^{\frac{1}{2}}(\mathcal{D})$ quasi-norms of the error density $\rho$ {against} the size of smallest element $h_s$ on the uniform (left) and adaptive meshes (right)~\ref{fig:Uniform_Adaptive_Meshes}. $\left\Vert {\rho}\right\Vert_{L^1(\mathcal{D})}$ grows proportional to $h_s^{-1}$, while $\left\Vert {\rho}\right\Vert_{L^{\frac{1}{2}}(\mathcal{D})}$ remains unchanging throughout {the} meshes. The finest uniform and adaptive meshes in this plot have 33,153 and 13,061 DoF, respectively.}
	\label{fig:ex0_norm_quasi_norm}
\end{figure}

Figure~\ref{fig:ex0_norm_quasi_norm} plots the $L^1(\mathcal{D})$ norm and $L^{\frac{1}{2}}(\mathcal{D})$ quasi-norm of the error density $\rho$ against the size of the smallest element $h_s$ for uniform and adaptive meshes. The $L^{\frac{1}{2}}(\mathcal{D})$ quasi-norm is defined by
\begin{align}
	\lVert \rho \rVert_{L^\frac{1}{2}(\mathcal{D})} = \left(\int_{\mathcal D} \rho^{\frac{1}{2}} \right)^2.
\end{align}
{Because of} the slit singularity, the smallest element is the one at the origin $(0, 0)$. The growth of $\lVert \rho \rVert_{L^1(\mathcal{D})}$ with slope 1 and the constant $\lVert \rho \rVert_{L^\frac{1}{2}(\mathcal{D})}$ matches the Sobolev regularity $u, \varphi \in H^{\frac{3}{2} - \epsilon}(\mathcal{D})$, for $\epsilon > 0${, see}~\cite{Stephan_Whiteman_1988}. 

The asymptotic growth behavior of $\lVert \rho \rVert_{L^1(\mathcal{D})}$ for adaptive meshes starts occurrs at about the same value of the smallest mesh size, $h_s$, value as uniform ones (Figure \ref{fig:ex0_norm_quasi_norm}). With comparable \#DoF, adaptive meshes have {considerably} smaller $h_s$ than uniform meshes. {They more efficiently capture} the singularity than the uniform ones, {thereby providing} more reliable error estimates.  Observe that the numerical results for Example 1, which considers a lognormal constant field $a$, represent this example as well. When the coefficient $a$ is constant over the domain $\mathcal{D}$,  the QoI (i.e., a linear functional of u) scales linearly with $1/a$, as seen {in}
\begin{equation}
  	-\nabla \cdot \left( \nabla u(\mathbf x; \omega) \right) = \dfrac{f(\mathbf x)}{a(\omega)} \text{ for $\mathbf x  \in \mathcal D$.} \nonumber
\end{equation}
Therefore, we postpone the comparison between error and tolerances to Section~\ref{sec:example_1}. 

\subsubsection{Example 1}
\label{sec:example_1}
In this section, we consider the lognormal constant coefficient field
\begin{align}
	a(x) & \equiv \exp{\left(Y\right)}, & x\in D,
\end{align}
{where} $Y$ {is} a Gaussian random variable with mean 0 and variance $\sigma^2$. 

\begin{figure}[htbp]
	\centering
	\includegraphics[width=0.49\textwidth]{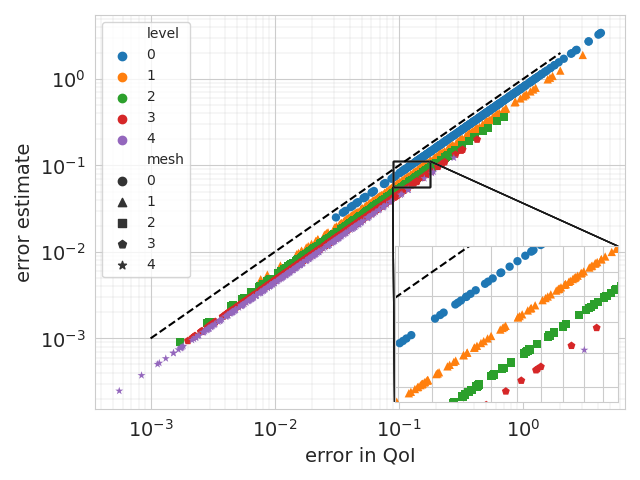}
	\hfill
	\includegraphics[width=0.49\textwidth]{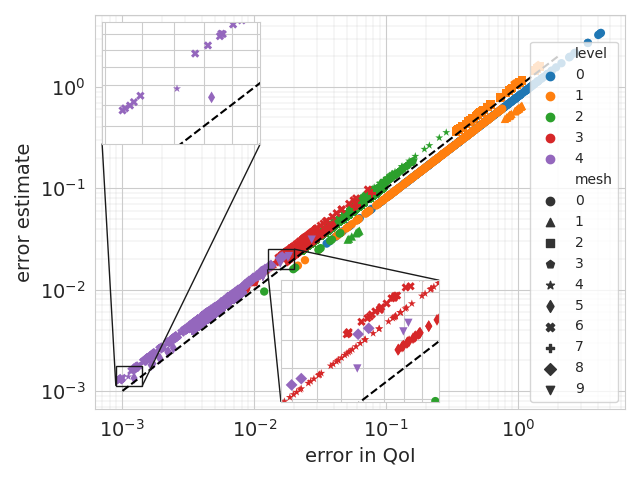}
	\caption{Example 1 with $\sigma^2 = 1$: Errors in QoI (compared to a reference solution) versus error estimates of $Q_\ell$ for SMLMC (left) and AMLMC (right) for their first 5 levels respectively. The dashed line {denotes} where the error estimate of $Q_\ell$ equals the error in $Q_\ell$. The dots are samples with colors representing at which MLMC level they belong {to}, whereas symbols {were} used to distinguish between auxiliary mesh levels. {Observe that} \#DoF {grew} with the index of the symbols.}
	\label{fig:logrv_err_est}
\end{figure}

Figure~\ref{fig:logrv_err_est} shows the error estimates against the computational errors in $Q_\ell$ for SMLMC and AMLMC. The errors are computed by comparing with a reference solution obtained on a refined adaptive mesh with 801,430 DoF. 

In SMLMC (left plot of Figure~\ref{fig:logrv_err_est}), all the samples in one level are obtained from the same mesh. Observe that each level forms a line pattern, and these line patterns share the same slope but {are at} an offset from each other. It {follows} that on the same mesh, the error estimate is different from the error (in $Q_\ell$) by a constant factor, which {is independent of} the random sample. As we increase the levels, the error estimates stagnate at an offset below the dashed line, {exposing} that the uniform meshes underestimate the errors.

In AMLMC case (right plot of Figure~\ref{fig:logrv_err_est}), one level contains the samples {that} satisfy a given criterion (refer to Algorithm~\ref{alg:sample_DQ}) and are distributed into a range of (adaptive) auxiliary meshes. Several line patterns {having} the same slope but different offsets can also be observed, with each line consisting of samples of the same mesh but belonging to different levels. {Because of} the variation of the stochastic coefficient, computing to a given accuracy incurs a higher cost for {certain} samples than {for} others. The {magnified} views {demonstrate} that the {proposed} stochastic mesh selection {makes} the error contribution smaller where it is {inexpensive} to do so and larger where it is more expensive to reduce the error. The contrast between cheap and expensive samples {follows} from using our sample-dependent stopping scheme in AMLMC. The error estimates {were} slightly overestimated in the last level of the shown adaptive meshes.

In {Examples} 1 and 2, we {fixed} a geometric sequence  $(\tol_\ell)_{\ell\ge0}$, with $\tol_0$ of AMLMC selected {from} (2.5) in~\cite{AcNum_MLMC}. In Example 1, the $\tol_\ell$ used for $\sigma^2 = 1$ and $\sigma^2 = 4$ are $\tol_\ell = 2 \cdot 4^{-\ell}$ and $\tol_\ell = 4^{1 - \ell}$, respectively. 

\begin{figure}[htbp]
	\centering
	\includegraphics[width=0.49\textwidth]{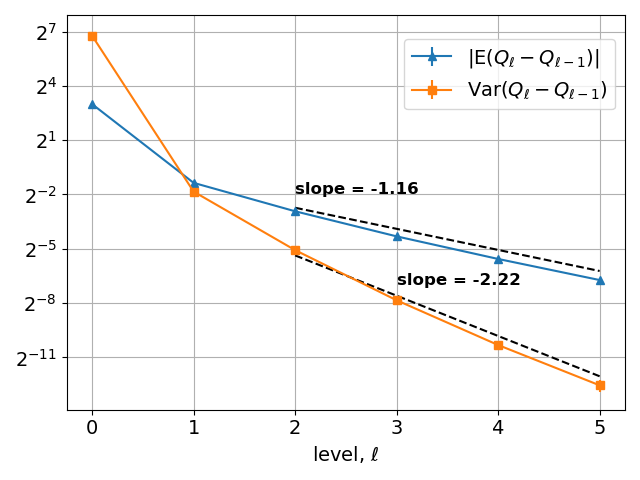}
	\hfill
	\includegraphics[width=0.49\textwidth]{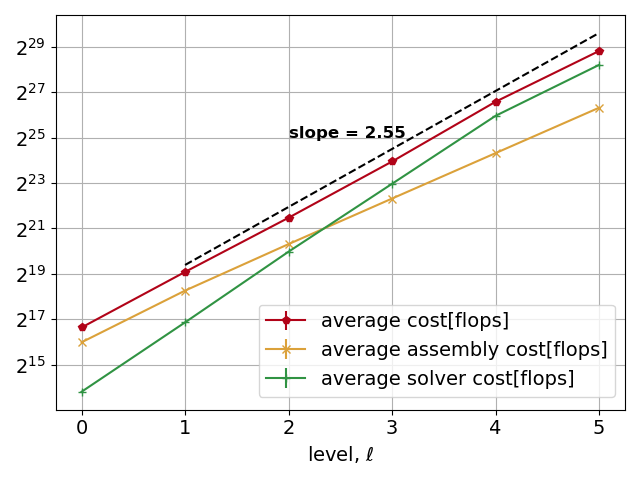}
	\\
	\caption{Example 1 with $\sigma^2=1$ using SMLMC. Left: Mean and variance of $Q_\ell - Q_{\ell-1}$ for level $\ell$, with a 95\% confidence interval. Right: The average work per sample (measured in flops) for level $\ell$.}
	\label{fig:ex1_sig_1_standard_MLMC_mean_var_cost}
\end{figure}

\begin{figure}[htbp]
	\centering
	\includegraphics[width=0.49\textwidth]{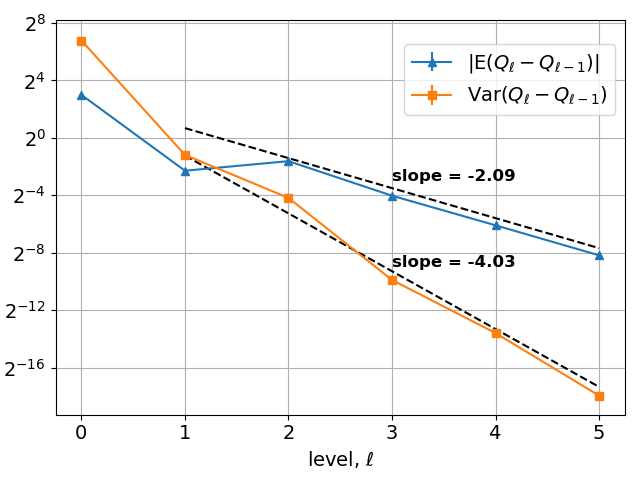}
	\hfill
	\includegraphics[width=0.49\textwidth]{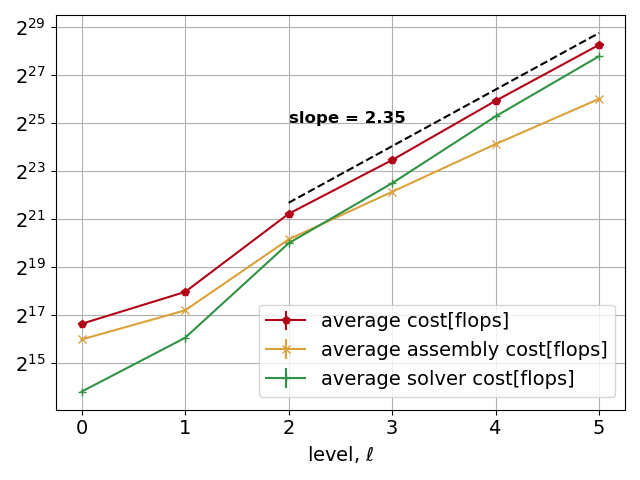}
	\\
	\caption{Example 1 with $\sigma^2=1$ using AMLMC. Left: Mean and variance of $Q_\ell - Q_{\ell-1}$ for level $\ell$, with a 95\% confidence interval. Right: The average work per sample (measured in flops) for level $\ell$.}
	\label{fig:ex1_sig_1_adaptive_MLMC_mean_var_cost}
\end{figure}

Figures~\ref{fig:ex1_sig_1_standard_MLMC_mean_var_cost} and \ref{fig:ex1_sig_1_adaptive_MLMC_mean_var_cost} show the convergence of mean $E_\ell = \lvert \mathbb{E} (Q_\ell - Q_{\ell - 1}) \rvert$ and variance $V_\ell$ and the growth of $W_\ell$ for SMLMC and AMLMC in Example 1 with $\sigma^2 = 1$. The mean $E_\ell$ and variance $V_\ell$ {of telescopic differences} of AMLMC converge twice {faster than the ones of} SMLMC with respect to $\ell$; {however, the work} $W_\ell$ {remains} essentially the same. {Note} that in our setting, {the sequence of $\tol$s} in AMLMC decays with a factor of 1/4, {controlling} the decay of the sharp error estimate of each sample. {From} Figure~\ref{fig:uniform_adaptive_error_estimate_DoFs}, the convergence rate of the error estimate for a fixed realization of the coefficient $a$ on adaptive meshes is proportional to the growth {in} \#DoF, {which is} twice the rate of the uniform meshes. Thus, with the level increments, the estimated weak error {in AMLMC} should decay by a factor of 4; {however,} the cost {increases} by a factor of 4, which is the same cost {required} to globally refine a uniform mesh. In both schemes, the assembly cost is the dominating cost in the coarser levels, and is surpassed by the solver cost as we go deeper on the level $\ell$.

\begin{figure}[htbp]
	\centering
	\includegraphics[width=0.49\textwidth]{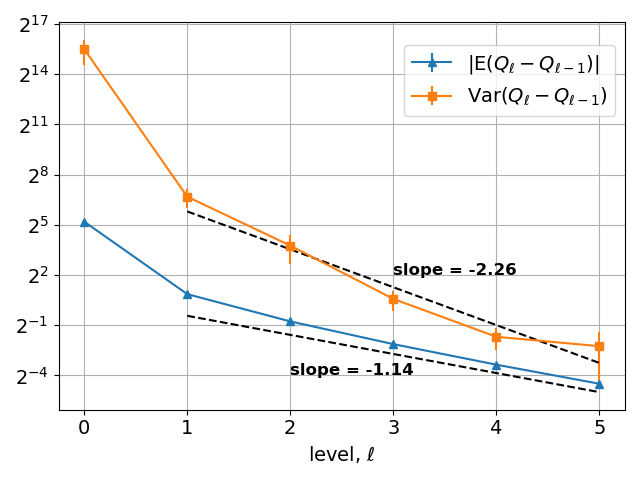}
	\hfill
	\includegraphics[width=0.49\textwidth]{Constant_field_sig_10e-1_lam_10e-1_UNIFORM_standard_MLMC_Gauss_25e-2_cost.png}
	\\
	\caption{Example 1 with $\sigma^2=4$ using SMLMC. Left: Mean and variance of $Q_\ell - Q_{\ell-1}$ for level $\ell$, with a 95\% confidence interval. Right: The average work per sample (measured in flops) for level $\ell$. 
	The cost of SMLMC on this example does not depend on $\sigma$, 
	but for easy reference we duplicate the graphs from Figure~\ref{fig:ex1_sig_1_standard_MLMC_mean_var_cost} here.} 
	\label{fig:ex1_sig_2_standard_MLMC_mean_var_cost}
\end{figure}

\begin{figure}[htbp]
	\centering
	\includegraphics[width=0.49\textwidth]{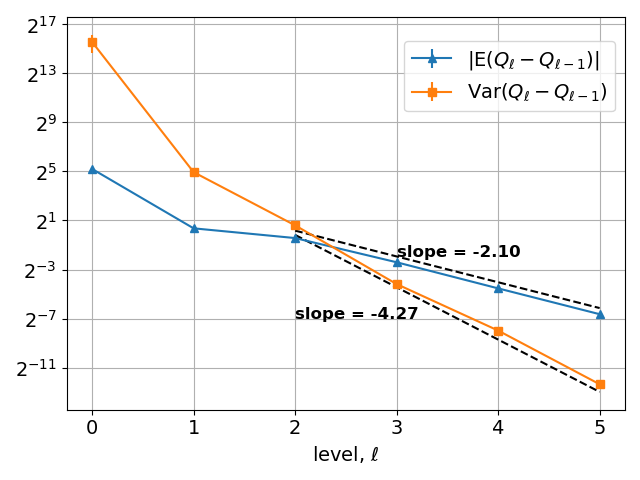}
	\hfill
	\includegraphics[width=0.49\textwidth]{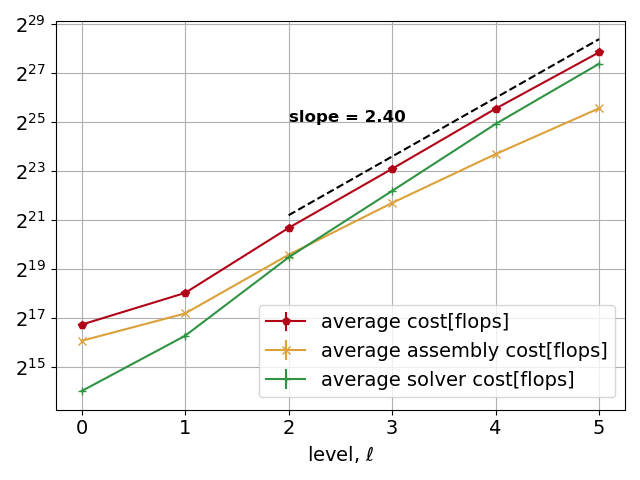}
	\\
	\caption{Example 1 with $\sigma^2=4$ using AMLMC. Left: Mean and variance of $Q_\ell - Q_{\ell-1}$ for level $\ell$, with a 95\% confidence interval. Right: The average work per sample (measured in flops) for level $\ell$.}
	\label{fig:ex1_sig_2_adaptive_MLMC_mean_var_cost}
\end{figure}

Figures~\ref{fig:ex1_sig_2_standard_MLMC_mean_var_cost} and \ref{fig:ex1_sig_2_adaptive_MLMC_mean_var_cost} show the convergence of mean $E_\ell$ and variance $V_\ell$ and the growth of $W_\ell$ for SMLMC and AMLMC in Example 1 with $\sigma^2 = 4$. We observed the mean and variance {of AMLMC} converges twice {faster than} SMLMC again. {Moreover,} with the same number of samples, AMLMC {displays} a more stable variance convergence than SMLMC.

\begin{figure}[htbp]
    \begin{subfigure}[b]{0.49\textwidth}
		\centering
		\includegraphics[width=\textwidth]{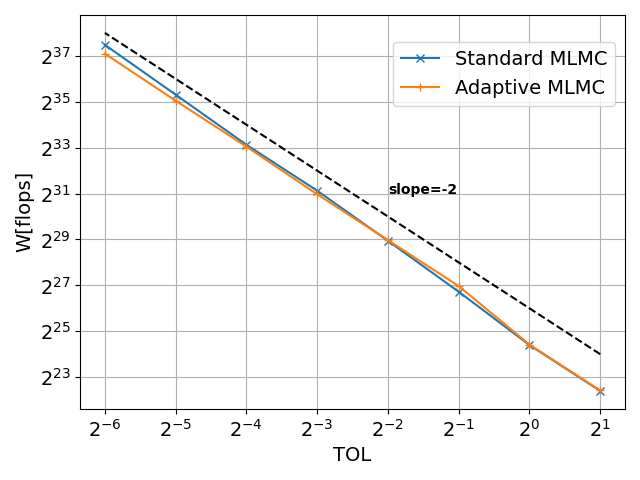}
		\caption{$\sigma^2 = 1$}
		\label{fig:constant_field_sig_10e-1_cost}
	\end{subfigure}
	\hfill
	\begin{subfigure}[b]{0.49\textwidth}
		\centering
		\includegraphics[width=\textwidth]{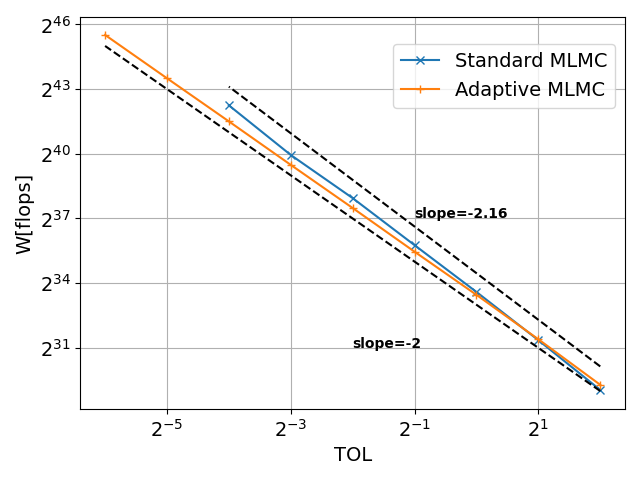}
		\caption{$\sigma^2 = 4$}
		\label{fig:constant_field_sig_20e-1_cost}
	\end{subfigure}
	\\
    \begin{subfigure}[b]{0.49\textwidth}
		\centering
		\includegraphics[width=\textwidth]{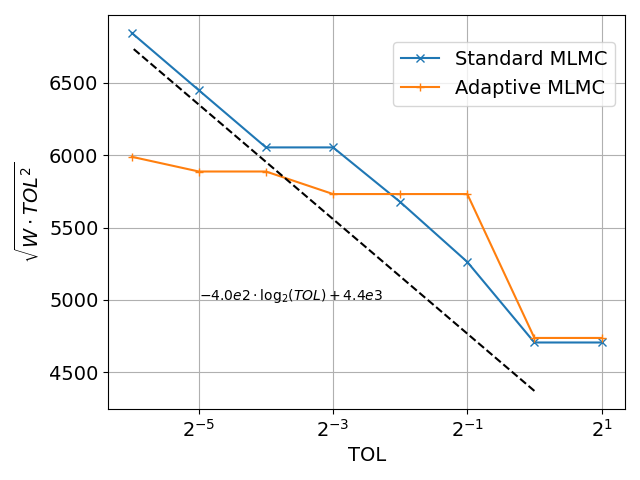}
		\caption{$\sigma^2 = 1$}
		\label{fig:constant_field_sig_10e-1_cost_log}
	\end{subfigure}
	\hfill
	\begin{subfigure}[b]{0.49\textwidth}
		\centering
		\includegraphics[width=\textwidth]{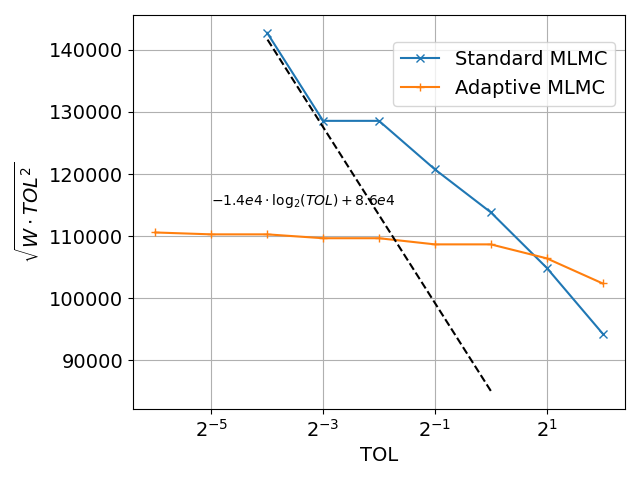}
		\caption{$\sigma^2 = 4$}
		\label{fig:constant_field_sig_20e-1_cost_log}
	\end{subfigure}	
	\caption{Example 1: Estimated MLMC work $W$ (in flops) (\ref{fig:constant_field_sig_10e-1_cost}, \ref{fig:constant_field_sig_20e-1_cost}) and $\sqrt{W\cdot \tol^2}$ (\ref{fig:constant_field_sig_10e-1_cost_log}, \ref{fig:constant_field_sig_20e-1_cost_log}) versus error tolerance $\tol$ for SMLMC and AMLMC for two variances in Example 1.}
	\label{fig:Ex1_MLMC_WORK_TOL}
\end{figure}

Figure~\ref{fig:Ex1_MLMC_WORK_TOL} shows the estimated MLMC work complexity against the given tolerance $\tol$. Figure~\ref{fig:constant_field_sig_10e-1_cost} and \ref{fig:constant_field_sig_20e-1_cost} show the average MLMC work $W$ (measured in flops) against $\tol$ for $\sigma^2 = 1$ and $\sigma^2 = 4$. With {increase in} variance, {there is a fine} distinction between the work complexity of the {both} schemes. To reveal or {exclude} the logarithmic term in the MLMC complexity, the quantity $\sqrt{W \cdot \tol^2}$ against $\tol$ is shown in Figure~\ref{fig:constant_field_sig_10e-1_cost_log} and \ref{fig:constant_field_sig_20e-1_cost_log}. When $W \propto \tol^{-2} \log^2(\tol)$, $\sqrt{W \cdot \tol^2} \propto \log(\tol)$. In both Figure~\ref{fig:constant_field_sig_10e-1_cost_log} and \ref{fig:constant_field_sig_20e-1_cost_log}, the line plot of SMLMC {grows} linearly. {The results show that} the slope for the case $\sigma^2 = 4$ is larger than for $\sigma^2 = 1$, {whereas} the line plot of AMLMC does not demonstrate a dependence on $\tol$ as $\tol$ decreases.  

\begin{figure}[htbp]
	\centering
	\includegraphics[width=0.48\textwidth]{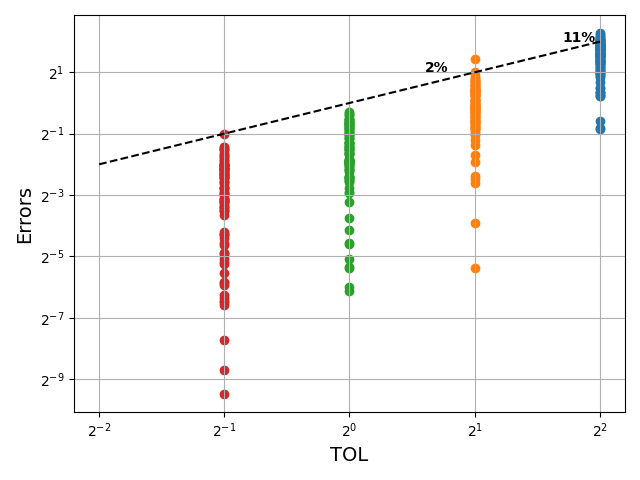}
	\includegraphics[width=0.48\textwidth]{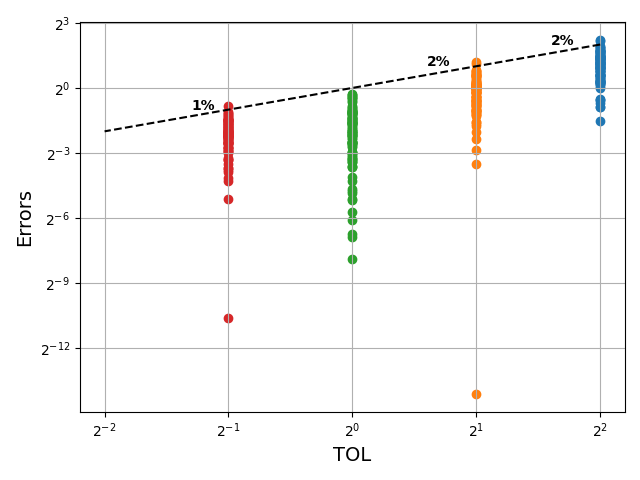}	
	\caption{Example 1 with $\sigma^2=4$: Estimated error for four $\tol$ values with SMLMC (Left) and AMLMC (Right). The reference line is for {when} $\tol$ equals to the estimated error, and the percentages are calculated from the number of estimator realizations for which the estimated error is above $\tol$.}
	\label{fig:Ex_1_MLMC_error_TOL}
\end{figure}

In Figure~\ref{fig:Ex_1_MLMC_error_TOL}, we {verified} the accuracy of our MLMC algorithm by comparing the errors of the MLMC estimator and the given $\tol$. The errors {were} computed against a reference value, which is approximately 38.7, given by an AMLMC estimator with $\tol = 2^{-6}$. {For all} cases except the largest {considered} $\tol$ for SMLMC, the percentage of samples having larger errors than {their prescribed} $\tol$ {was $<$} $5\%$. The result is consistent with our choice of 95\% success probability. 

\subsubsection{Example 2}
\label{sec:example_2}

\begin{figure*}[]
	\centering
	\begin{tabular}{M{55mm}|M{55mm}}
		\toprule
		$\sigma^2 = 1$ & $\sigma^2 = 4$ \\
		\midrule
		\includegraphics[width=0.4\textwidth]{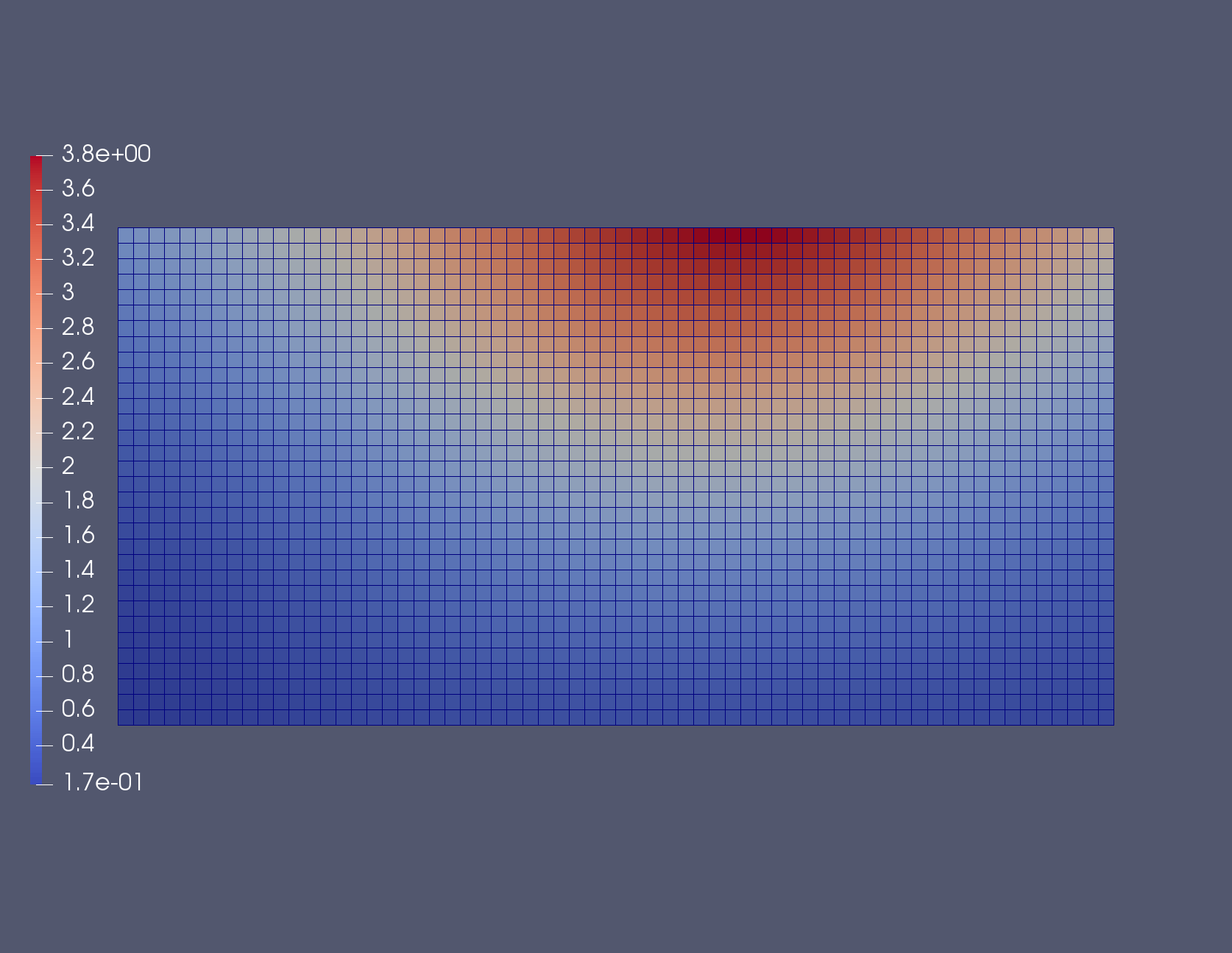} &
		\includegraphics[width=0.4\textwidth]{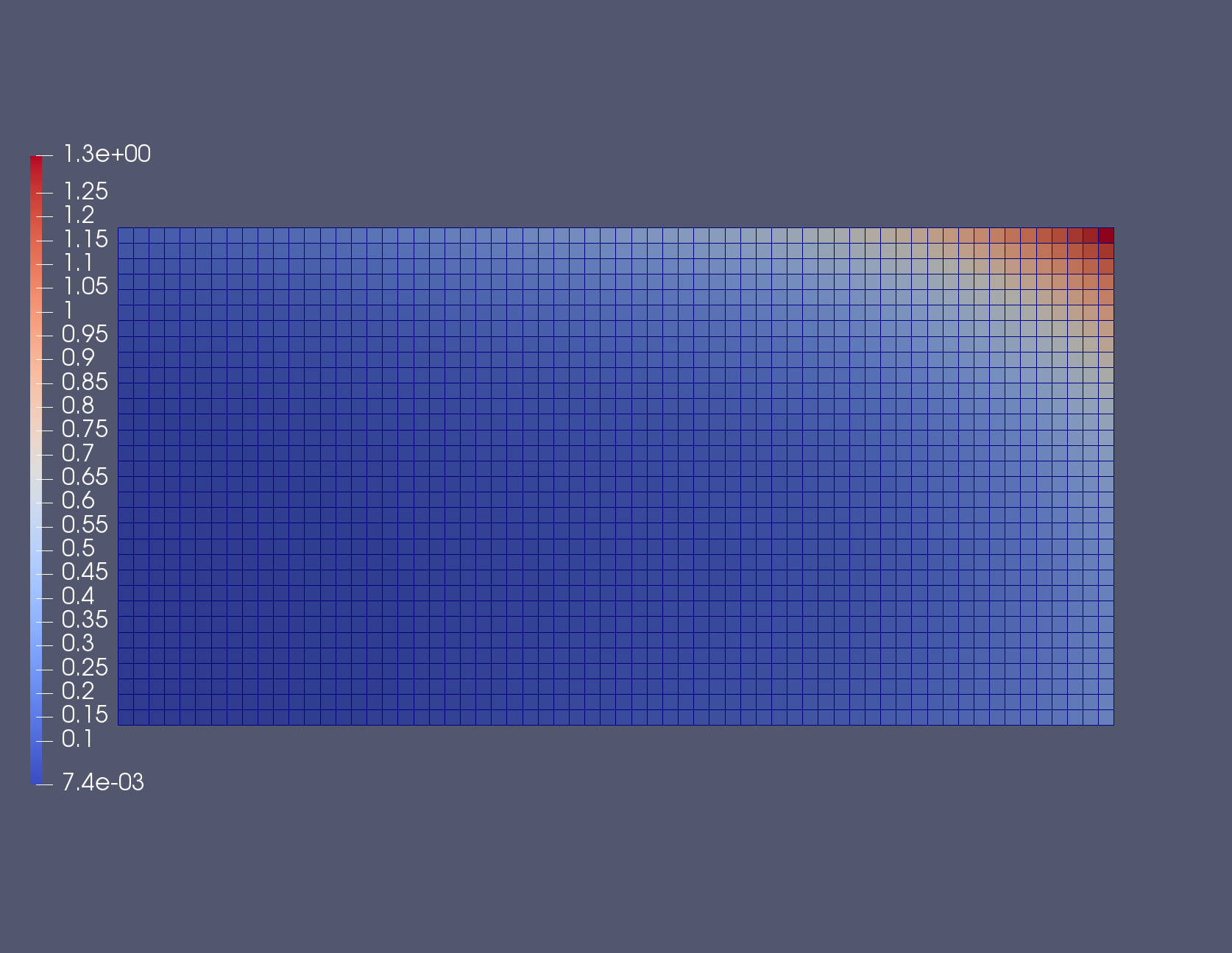}
		\\
		\midrule
		\includegraphics[width=0.4\textwidth]{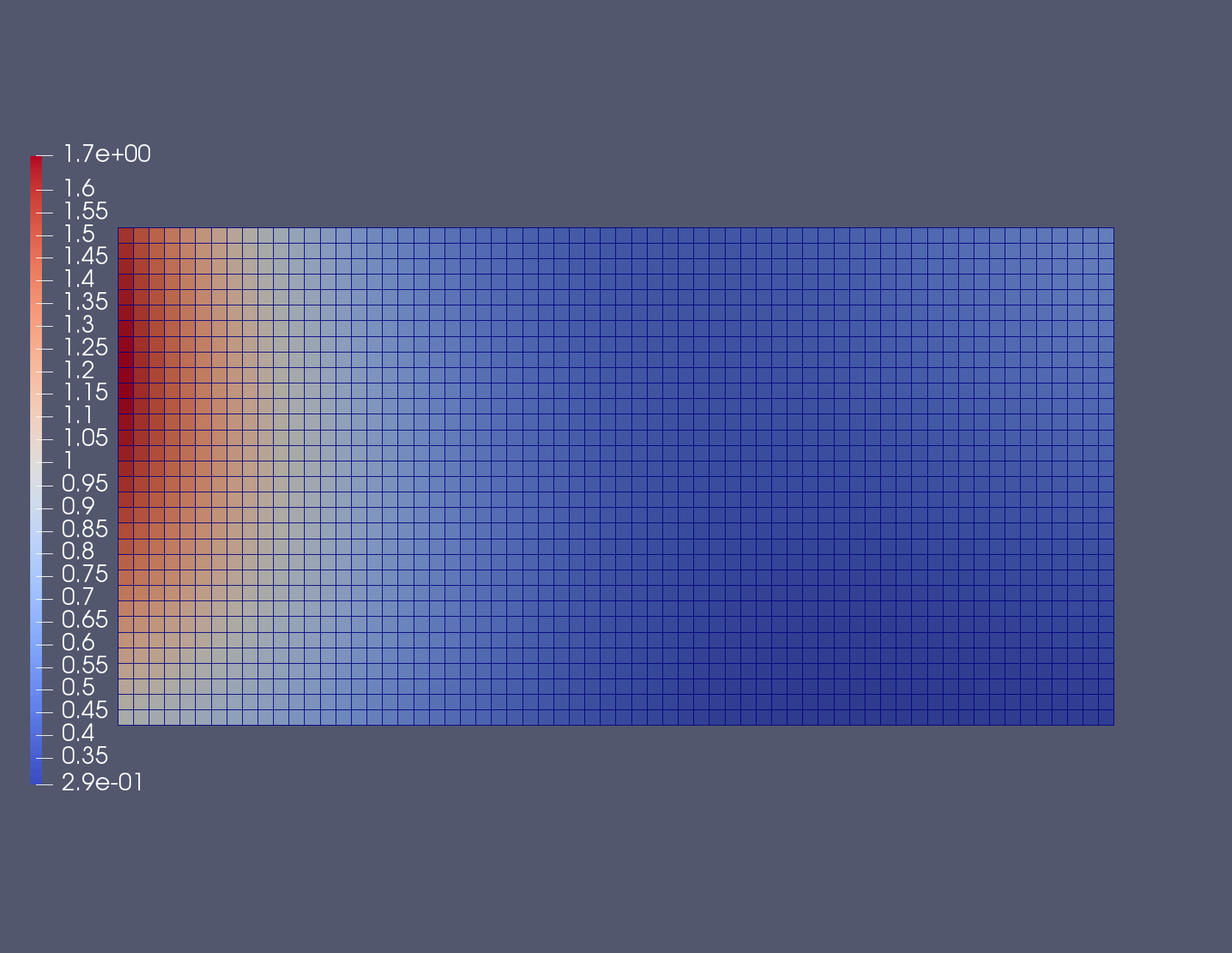} &
		\includegraphics[width=0.4\textwidth]{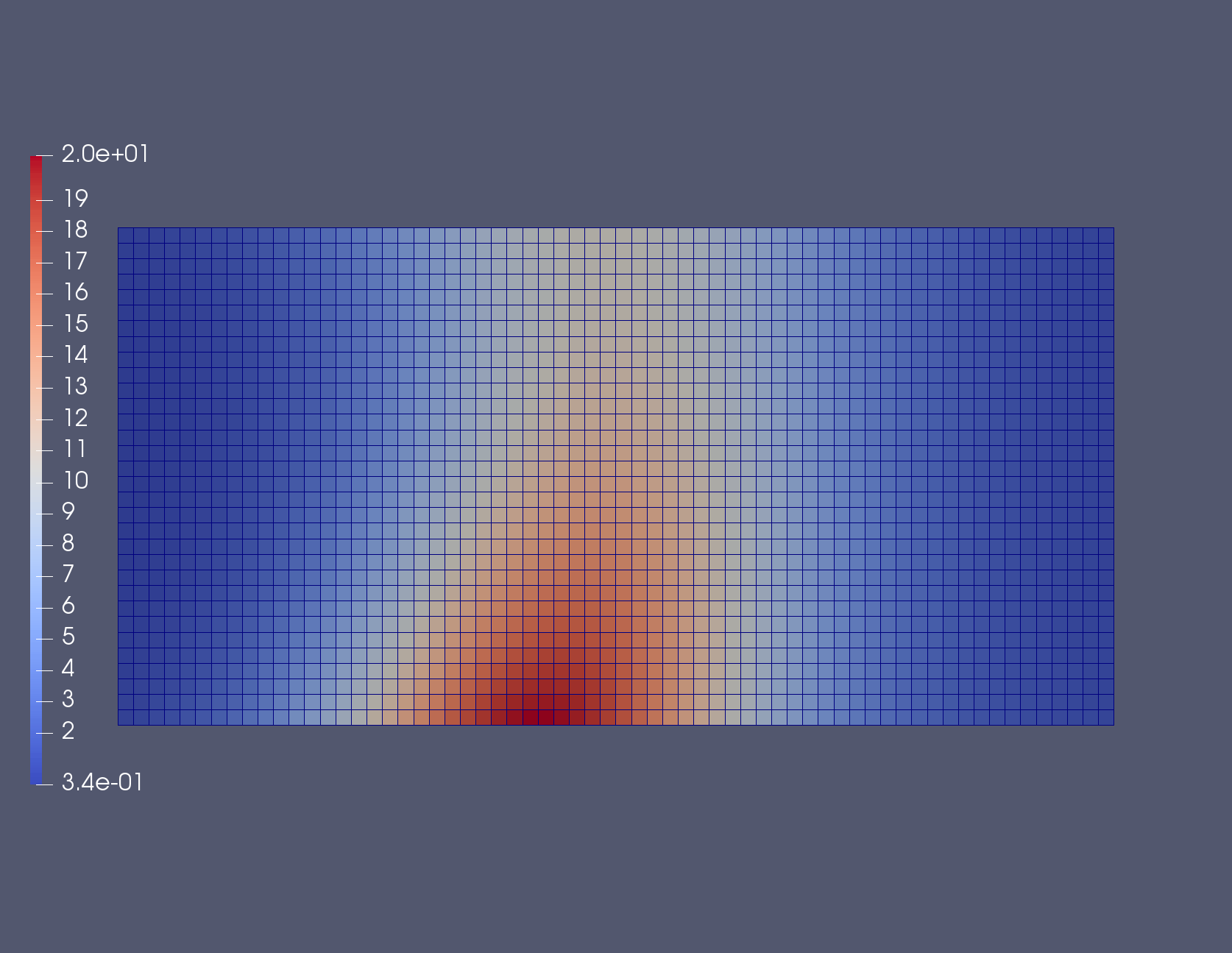}
		\\
		\bottomrule
	\end{tabular}
	
	\caption{Example 2 with $\sigma^2 = 1$ (Left) and $\sigma^2 = 4$ (Right): A few realizations of the random field $a$~\eqref{eq:random_field_series}}
	\label{fig:matern_random_field_realizations}
\end{figure*}

Now let $a$ be a lognormal random field
\begin{align}
	a(x, \omega) & = \exp \left( \sum_{i=0}^{I - 1} \xi_i(\omega) \sqrt{\lambda_i} \theta_i(x) \right), & x \in D,
	\label{eq:random_field_series}
\end{align}
with $I=256$ terms, where $\xi_i$ are i.i.d. $\mathcal{N}(0,1)$ and
$\lambda_i$ are Fourier coefficients corresponding to trigonometric 
basis functions $\theta_i$ in the Fourier expansion of the Mat\'ern covariance
\begin{align}
	C(h) & = \frac{\sigma^2}{2^{\nu - 1} \Gamma(\nu)} 
	\left( \sqrt{2\nu} \frac{h}{r}\right)^{\nu} 
	K_{\nu} \left( \sqrt{2\nu} \frac{h}{r}\right), 
	& h\geq 0,
	\label{eq:cov_function}
\end{align}
with correlation length $r=1$ and smoothness parameter $\nu = 6.5$,
where $\Gamma$ is the gamma function and $K_{\nu}$ is the modified 
Bessel function of the second kind. 
The Fourier coefficients, arranged in descending order, have the decay 
$\lambda_i \leq C i^{1 + \frac{2\nu}{d}}$; see~\cite{bachmayr2018representations}.
The connection to the Mat\'ern covariance with these choices of $\nu$ and $r$ is 
{holds because} the latter {yields} the random field pathwise $C^3(\mathcal{D})$ regularity; 
see~\cite{scheuerer2010regularity}. This in turn ensures that 
$u(\cdot, \omega)$ and $\varphi(\cdot, \omega)$ $\in C^3(\mathcal{D})$, as required 
{for computing} the error density~\cite{adFEM_our, adOverview_our}. 
Figure~\ref{fig:matern_random_field_realizations} shows a few random field realizations 
for $\sigma^2 = 1$ and $\sigma^2 = 4$. 

Here, the pathwise regularity of the coefficient $a(\cdot, \omega)$ is {sufficient} to make the numerical quadrature error negligible in our computations, even on the coarsest mesh. 

{Many studies} on sampling Gaussian random fields {have been reported, see}~\cite{HARBRECHT2012, lord_powell_shardlow_2014, Feischl_H_matrices, LITVINENKO2019, SCHWAB2006100, Ghanem_Spanos, Graham_et_al_SINUM_2018, LLSY2019}. 
The $\theta_i$, $\lambda_1$, and $\xi_i$ used are {obtained} from the Fourier series expansion method as proposed in~\cite{bachmayr2018representations} {and} considered as a continuous version of the circulant embedding approach~\cite{Graham_et_al_SINUM_2018, dietrich1997fast}. The method first periodically extends the covariance operator $C$ and then the Fourier coefficients $\{\lambda_i\}$ of the covariance operator can be efficiently computed. 

In Example 2, we chose $\tol_\ell = 2 \cdot 4^{-\ell}$ and $\tol_\ell = 4^{1 - \ell}$ for the AMLMC cases $\sigma^2 = 1$ and $\sigma^2 = 4$, respectively. 

\begin{figure}[htbp]
	\centering
	\includegraphics[width=0.49\textwidth]{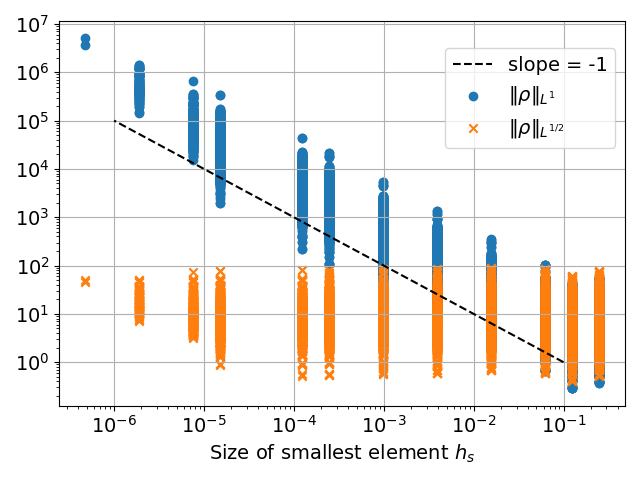}
	\hfill
	\includegraphics[width=0.49\textwidth]{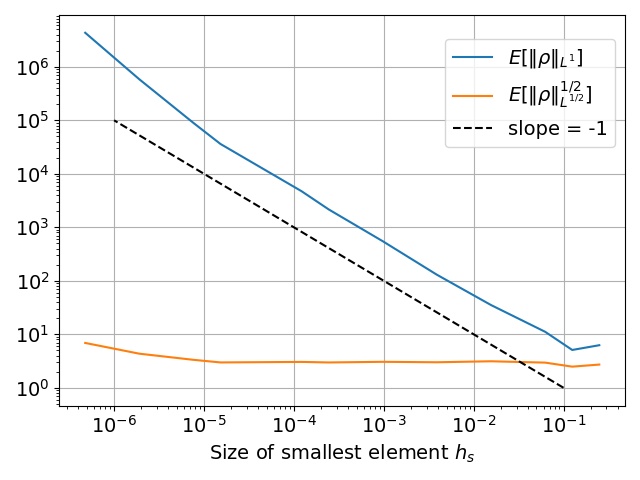}
	\\
	\caption{Example 2 with $\sigma^2 = 1$: $L^1(\mathcal{D})$ norm and $L^{\frac{1}{2}}(\mathcal{D})$ quasi-norms of the error density $\rho$ (left) and the expected value of the norms and quasi-norms on each mesh (right) {against} the size of {the} smallest element $h_s$ on a hierarchy of adaptive meshes (Figure~\ref{fig:Uniform_Adaptive_Meshes}).}
	\label{fig:ex2_sig1_adaptive_norms_quasi_norms}
\end{figure}

\begin{figure}[htbp]
	\centering
	\includegraphics[width=0.49\textwidth]{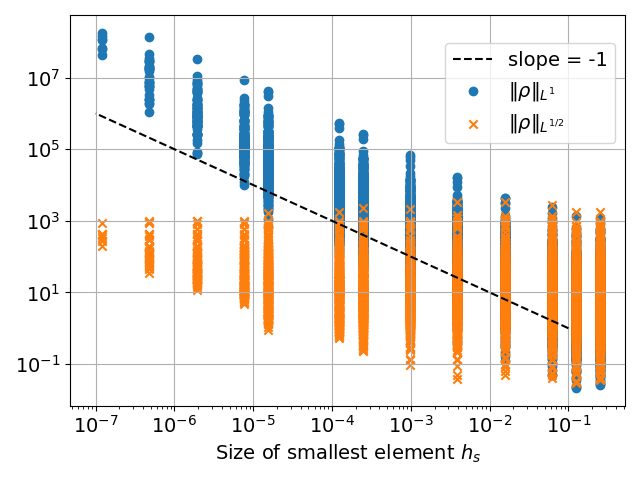}
	\hfill
	\includegraphics[width=0.49\textwidth]{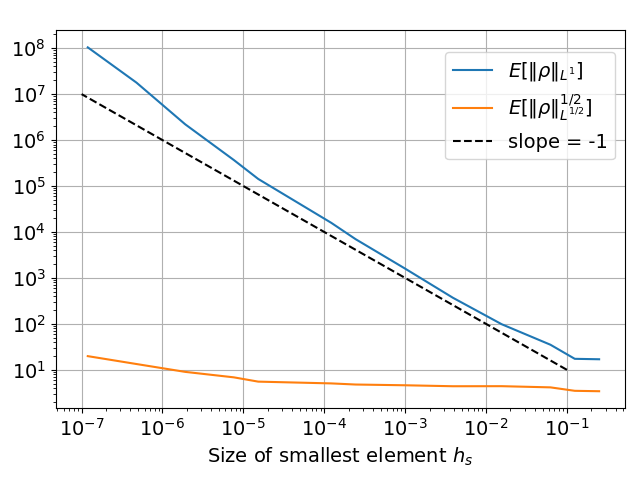}
	\\
	\caption{Example 2 with $\sigma^2 = 4$: $L^1(\mathcal{D})$ norm and $L^{\frac{1}{2}}(\mathcal{D})$ quasi-norms of the error density $\rho$ (left) and the expected value of the norms and quasi-norms on each mesh (right) {against}  the size of {the} smallest element $h_s$ on a hierarchy of adaptive meshes (Figure~\ref{fig:Uniform_Adaptive_Meshes}).}
	\label{fig:ex2_sig2_adaptive_norms_quasi_norms}
\end{figure}

Figures~\ref{fig:ex2_sig1_adaptive_norms_quasi_norms} and \ref{fig:ex2_sig2_adaptive_norms_quasi_norms} show the $L^1(\mathcal{D})$ norms and $L^{\frac{1}{2}}(\mathcal{D})$ quasi-norms of the samples from Example 2 with $\sigma^2 = 1$ and $\sigma^2 = 4$, respectively. {As shown in} Figure~\ref{fig:ex0_norm_quasi_norm}, the $L^\frac{1}{2}(\mathcal{D})$ quasi-norms remain stable across meshes, while the $L^{1}(\mathcal{D})$ norms grow with the rate $h_s^{-1}$.

\begin{figure}[htbp]
	\centering
	\includegraphics[width=0.49\textwidth]{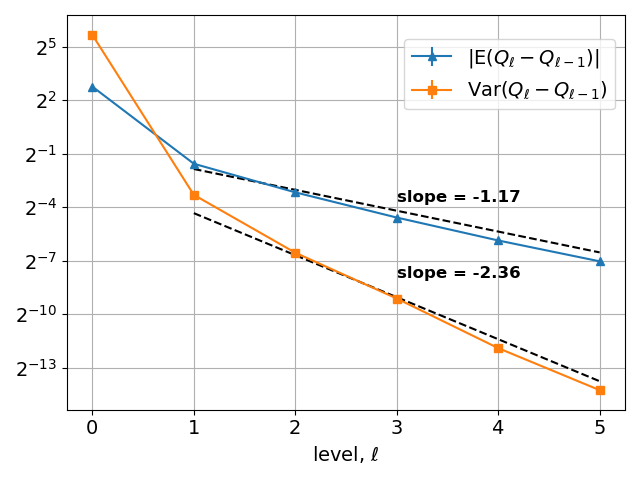}
	\hfill
	\includegraphics[width=0.49\textwidth]{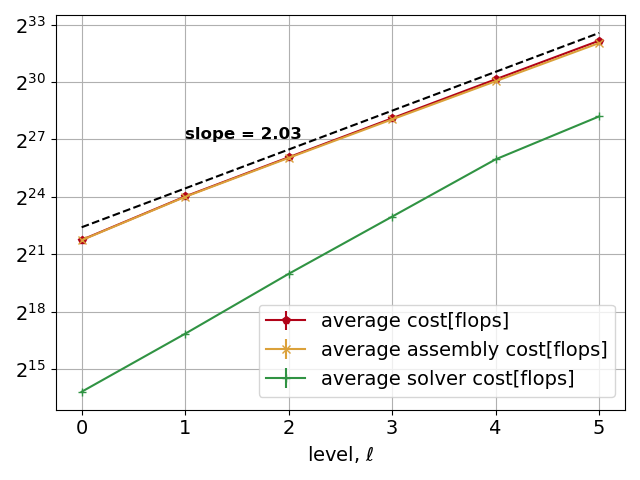}
	\\
	\caption{Example 2 with $\sigma^2=1$ using SMLMC. Left: Mean and variance of $Q_\ell - Q_{\ell-1}$ for level $\ell$. Right: The average work per sample (measured in flops) for level $\ell$.}
	\label{fig:ex2_sig_1_uniform_MLMC_mean_var_cost}
\end{figure}

\begin{figure}[htbp]
	\centering
	\includegraphics[width=0.49\textwidth]{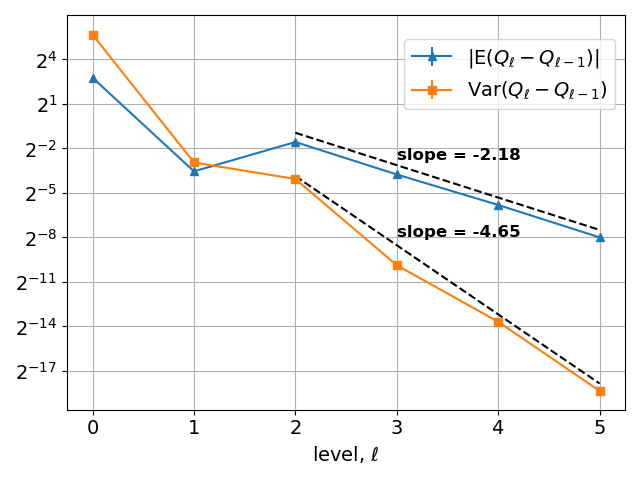}
	\hfill
	\includegraphics[width=0.49\textwidth]{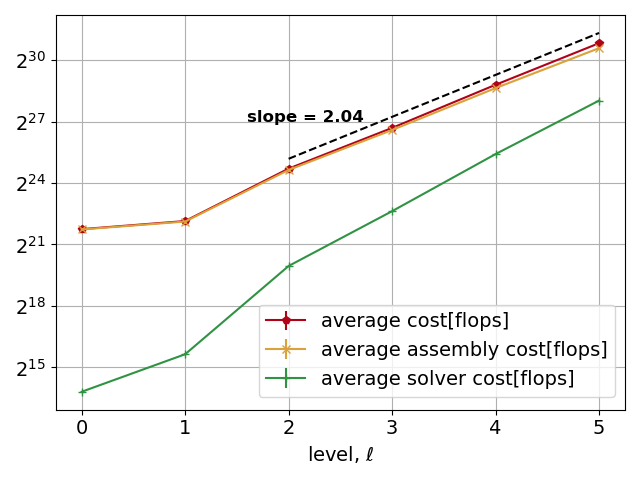}
	\\
	\caption{Example 2 with $\sigma^2=1$ AMLMC. Left: Mean and variance of $Q_\ell - Q_{\ell-1}$ for level $\ell$. Right: The average work per sample (measured in flops) for level $\ell$.}
	\label{fig:ex2_sig_1_adaptive_MLMC_mean_var_cost}
\end{figure}

\begin{figure}[htbp]
	\centering
	\includegraphics[width=0.49\textwidth]{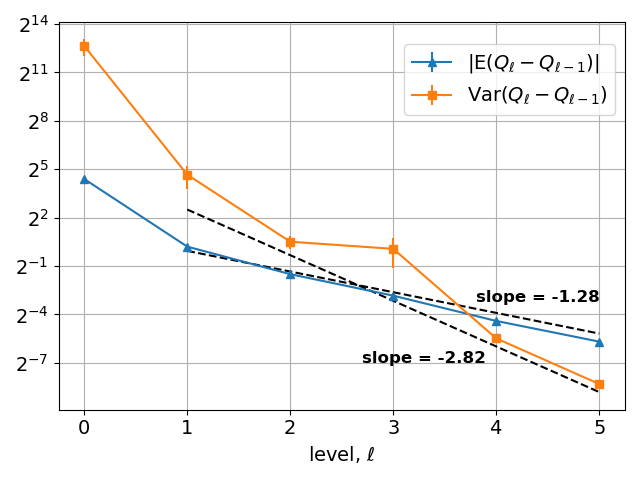}
	\hfill
	\includegraphics[width=0.49\textwidth]{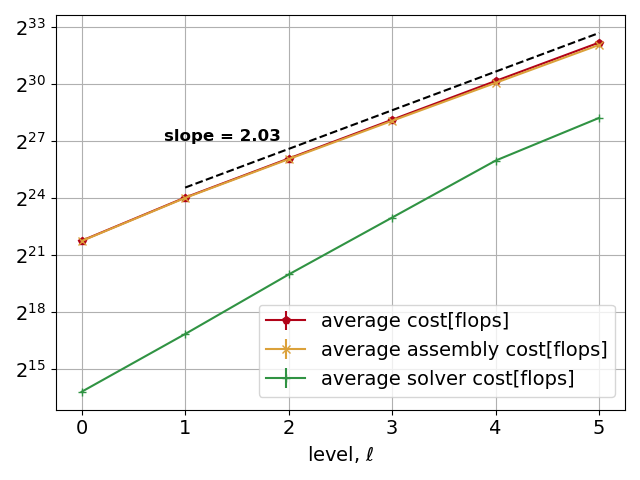}
	\\
	\caption{Example 2 with $\sigma^2=4$ using SMLMC. Left: Mean and variance of $Q_\ell - Q_{\ell-1}$ for level $\ell$. Right: The average work per sample (measured in flops) for level $\ell$.}
	\label{fig:ex2_sig_2_uniform_MLMC_mean_var_cost}
\end{figure}

\begin{figure}[htbp]
	\centering
	\includegraphics[width=0.49\textwidth]{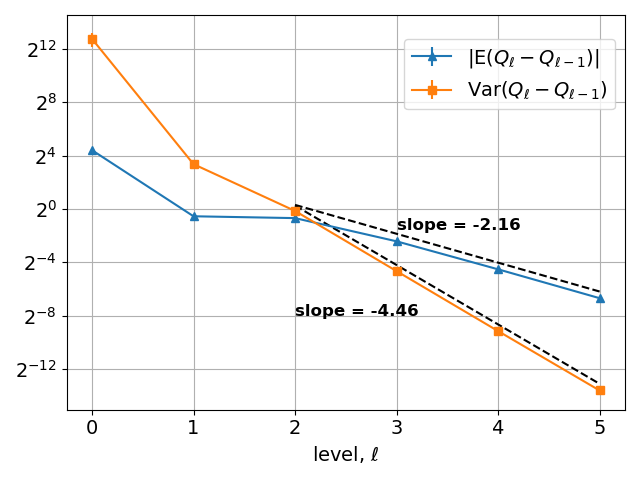}
	\hfill
	\includegraphics[width=0.49\textwidth]{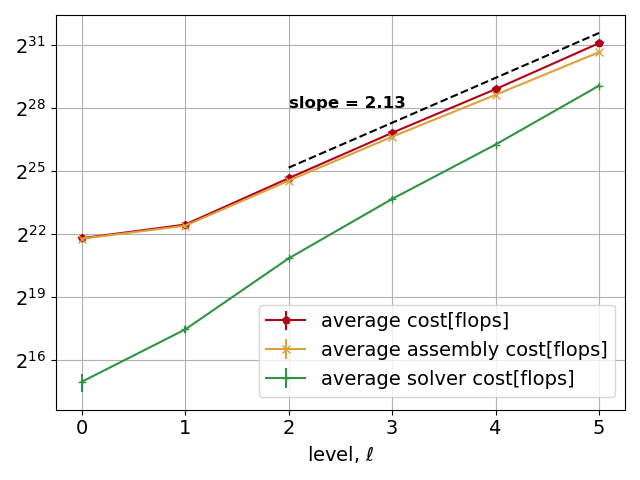}
	\\
	\caption{Example 2 with $\sigma^2=4$ using AMLMC. Left: Mean and variance of $Q_\ell - Q_{\ell-1}$ for level $\ell$. Right: The average work per sample (measured in flops) for level $\ell$.}
	\label{fig:ex2_sig_2_adaptive_MLMC_mean_var_cost}
\end{figure}

Figures~\ref{fig:ex2_sig_1_uniform_MLMC_mean_var_cost} and \ref{fig:ex2_sig_1_adaptive_MLMC_mean_var_cost} {show} the convergence of mean $E_\ell$ and variance $V_\ell$, and the cost growth $W_\ell$ in Example 2 with $\sigma^2 = 1$ for SMLMC and AMLMC. Figure~\ref{fig:ex2_sig_2_uniform_MLMC_mean_var_cost}, and \ref{fig:ex2_sig_2_adaptive_MLMC_mean_var_cost} give the result with $\sigma^2 = 4$. In Example 2, the cost is dominated by the assembly {caused by the} random fields {that were generated}. 

Figures~\ref{fig:ex2_sig1_adaptive_norms_quasi_norms}-\ref{fig:ex2_sig_2_adaptive_MLMC_mean_var_cost} are comparable to the corresponding figures in Example~1, {demonstrating that} SMLMC and AMLMC have similar behaviors as in Example~1.

\begin{figure}[htbp]
	\centering
	\includegraphics[width=0.48\textwidth]{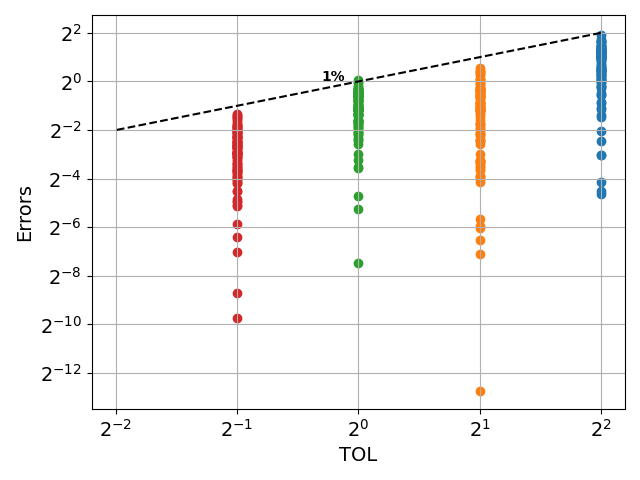}
	\includegraphics[width=0.48\textwidth]{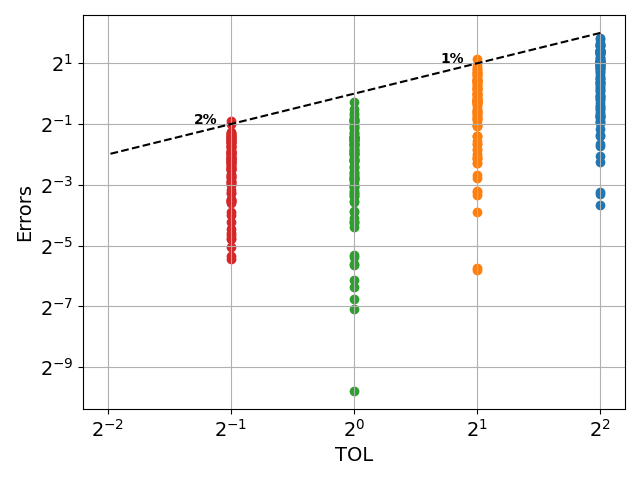}
	\caption{Example 2 with $\sigma^2=4$. Estimated error for four $\tol$ values with SMLMC (left) and AMLMC (right). The reference line {shows when} $\tol$ equals the estimated error, and the percentages are calculated {using} the number of estimator realizations for which the estimated error is above $\tol$.}
	\label{fig:Ex_2_MLMC_error_TOL}
\end{figure}

{Similar to Figure~\ref{fig:Ex_1_MLMC_error_TOL}, Figure~\ref{fig:Ex_2_MLMC_error_TOL} shows a graph of} the errors computed against a reference value ($\sim$ 22.6), given by an AMLMC estimator with $\tol = 2^{-7}$. {For} all the cases except the largest {considered} $\tol$ for AMLMC, the percentage of samples having larger errors than {their prescribed} $\tol$ is $< 5\%$. This is consistent with our choice of 95\% success probability.

\begin{figure}[htbp]
	\begin{subfigure}[b]{0.49\textwidth}
		\centering
		\includegraphics[width=\textwidth]{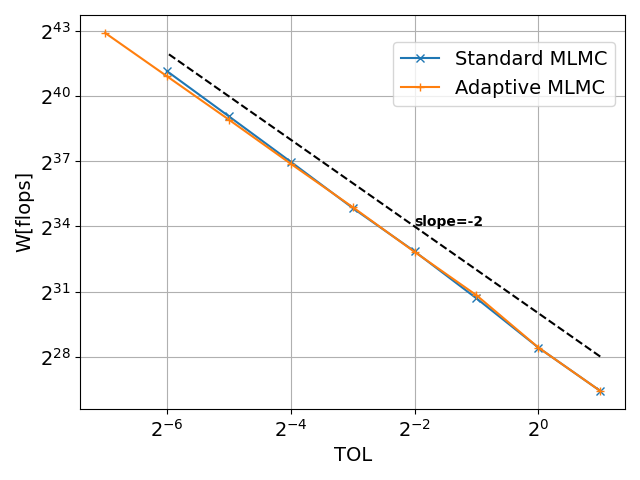}
		\caption{$\sigma^2 = 1$}
		\label{fig:matern_field_sig_10e-1_cost}
	\end{subfigure}
	\hfill
	\begin{subfigure}[b]{0.49\textwidth}
		\centering
		\includegraphics[width=\textwidth]{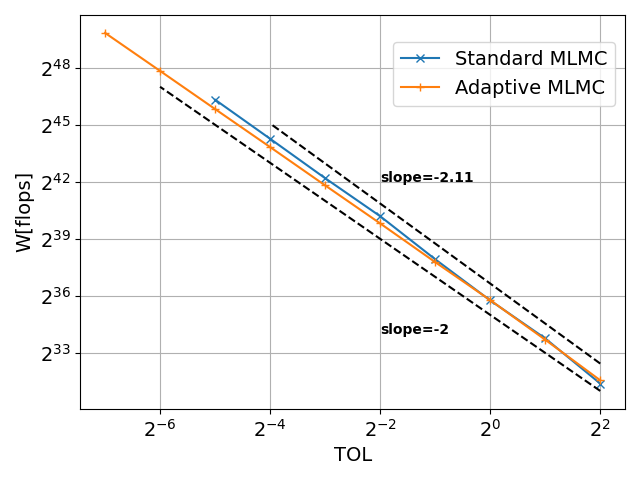}
		\caption{$\sigma^2 = 4$}
		\label{fig:matern_field_sig_20e-1_cost}
	\end{subfigure}
	\\
\begin{subfigure}[b]{0.49\textwidth}
	\centering
	\includegraphics[width=\textwidth]{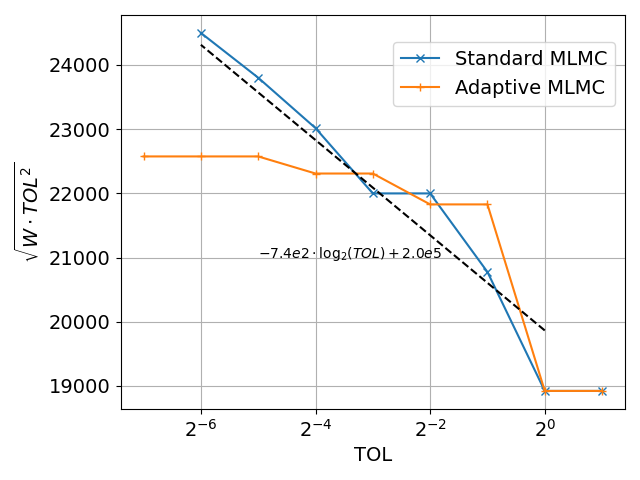}
	\caption{$\sigma^2 = 1$}
	\label{fig:matern_field_sig_10e-1_cost_log}
\end{subfigure}
	\hfill
	\begin{subfigure}[b]{0.49\textwidth}
		\centering
		\includegraphics[width=\textwidth]{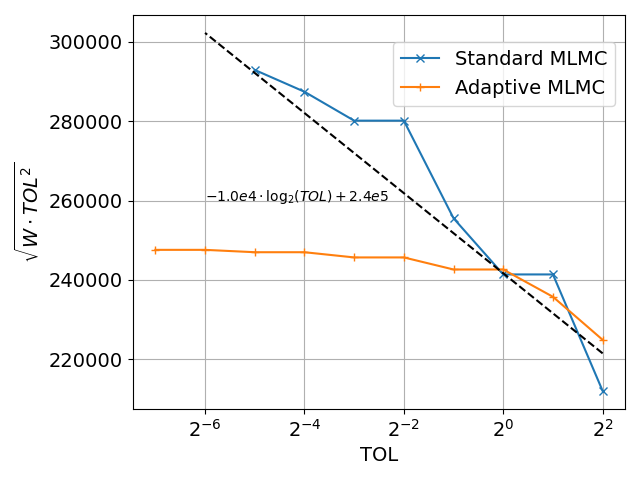}
		\caption{$\sigma^2 = 4$}
		\label{fig:matern_field_sig_20e-1_cost_log}
	\end{subfigure}	
	\caption{Example 2. Estimated MLMC work $W$ (in flops) (\ref{fig:matern_field_sig_10e-1_cost}, \ref{fig:matern_field_sig_20e-1_cost}) and $\sqrt{W\cdot \tol^2}$ (\ref{fig:matern_field_sig_10e-1_cost_log}, \ref{fig:matern_field_sig_20e-1_cost_log}) {against} error tolerance $\tol$ for SMLMC, and AMLMC for two variances $\sigma^2 = 1$ and $\sigma^2 = 4$.}
	\label{fig:Ex2_MLMC_WORK_TOL}
\end{figure}

Figure~\ref{fig:Ex2_MLMC_WORK_TOL} shows the estimated MLMC work complexity against the given tolerance $\tol$. Figures~\ref{fig:matern_field_sig_10e-1_cost} and \ref{fig:matern_field_sig_20e-1_cost} show the work $W$ (measured in flops) against $\tol$ for $\sigma^2 = 1$ and $\sigma^2 = 4$. {The increase in variance shows a fine} distinction between the work complexity of both schemes. {To reveal} or rule out the logarithmic term in the MLMC complexity, the quantity $\sqrt{W \cdot \tol^2}$ against $\tol$ is shown in Figures~\ref{fig:matern_field_sig_10e-1_cost_log} and \ref{fig:matern_field_sig_20e-1_cost_log}. When $W \propto \tol^{-2} \log^2(\tol)$ {and} $\sqrt{W \cdot \tol^2} \propto \log(\tol)$. In both figures, the line plot of SMLMC grows linearly. {The results show that the slope for the case} $\sigma^2 = 4$ is larger than $\sigma^2 = 1$, {whereas} the line plot of AMLMC does not demonstrate a dependence on $\tol$.

\begin{figure}[htbp]
	\centering
	\includegraphics[width=0.48\textwidth]{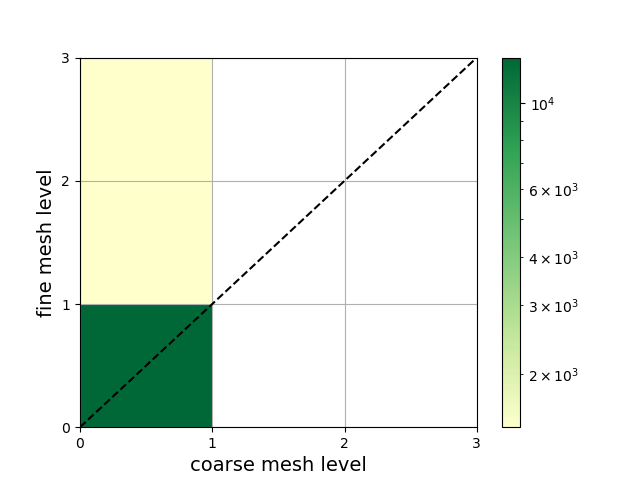} 
		\includegraphics[width=0.48\textwidth]{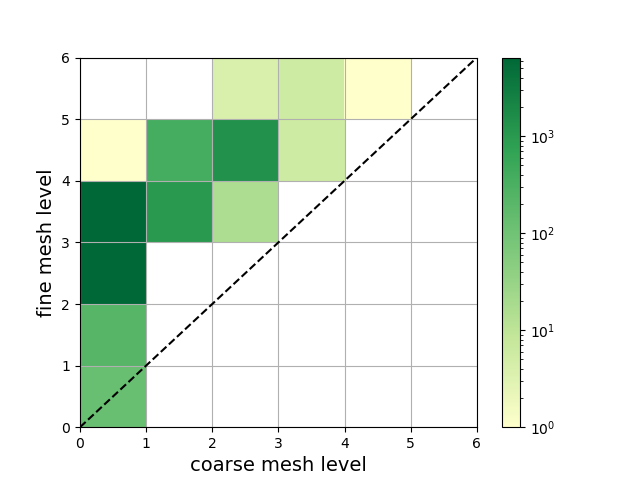} \\
			\includegraphics[width=0.48\textwidth]{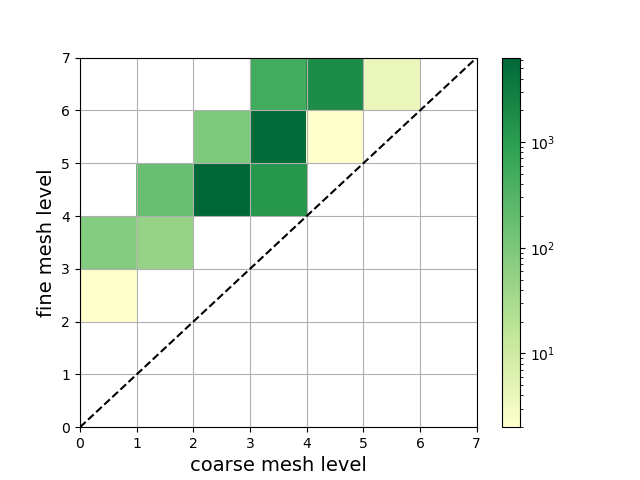} 
				\includegraphics[width=0.48\textwidth]{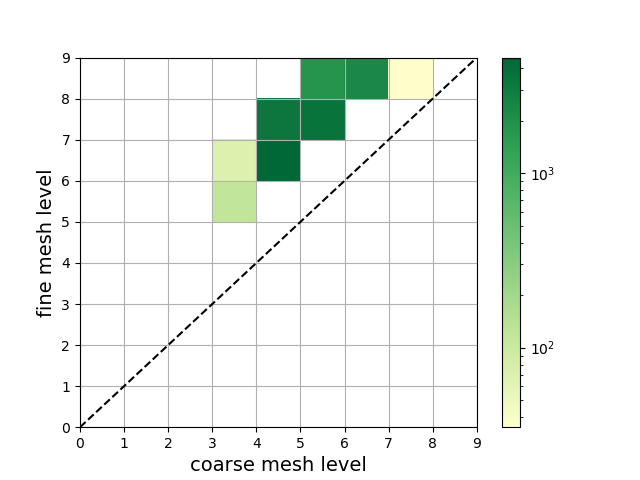} \\
					\includegraphics[width=0.48\textwidth]{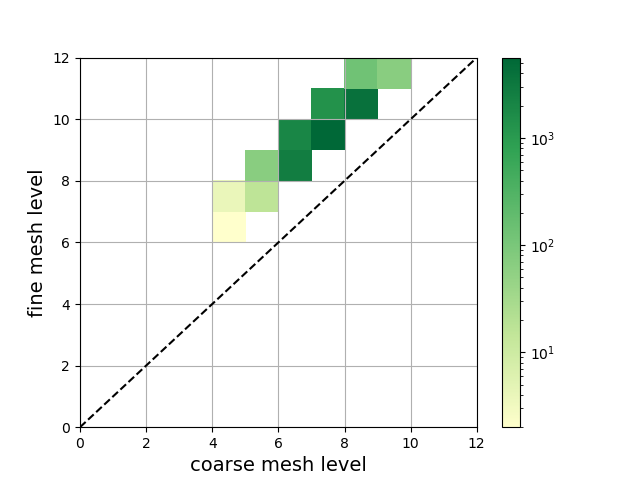}
	\caption{Example 2 with $\sigma^2=1$ using AMLMC: Stochastic stopping of coarse-fine pairs for {the} mesh level $\ell$. Top: $\ell=1,2$ from left to right. Middle: $\ell=3,4$ from left to right. Bottom: $\ell=5$. In each subplot, the dashed diagonal line {denotes} where the coarse mesh level equals the fine mesh level.}
	\label{fig:ex2_sig1_coarse_fine_pair}
\end{figure}

Figure~\ref{fig:ex2_sig1_coarse_fine_pair} plots the histograms of coarse-fine mesh level of AMLMC in Example 2 with $\sigma^2 = 1$ for level $\ell = 1, 2, 3, 4, 5$. {Because} the coarse-fine pair in our MLMC scheme uses the same random field, $\Delta Q_\ell (\omega)= (Q_\ell - Q_{\ell - 1}) (\omega)$ if the coarse mesh level equals the fine mesh level, thus reducing the efficiency of MLMC. With increasing $\ell$, the fraction of zero samples ($\Delta Q_\ell (\omega) = 0$) reduces and disappears starting from $\ell = 2$. 

\begin{figure}[htbp]
	\centering
	\includegraphics[width=0.48\textwidth]{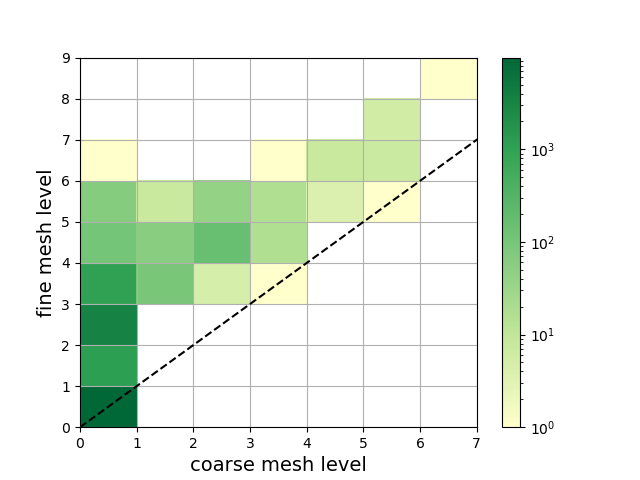} 
		\includegraphics[width=0.48\textwidth]{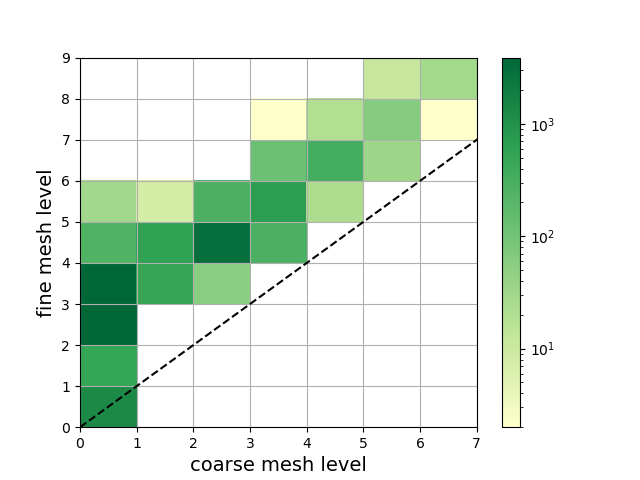} \\
			\includegraphics[width=0.48\textwidth]{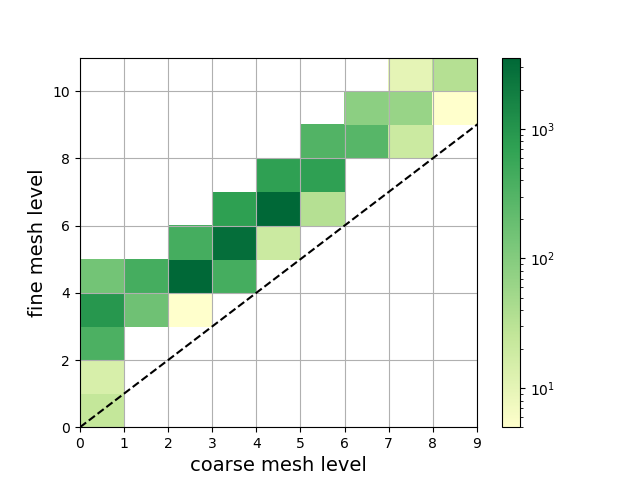}
				\includegraphics[width=0.48\textwidth]{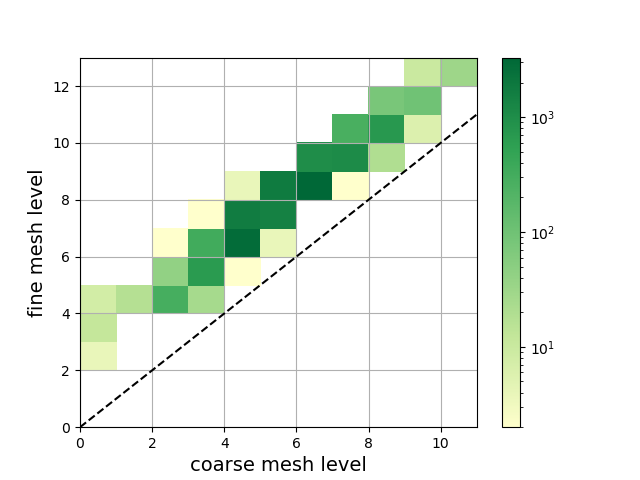} \\
					\includegraphics[width=0.48\textwidth]{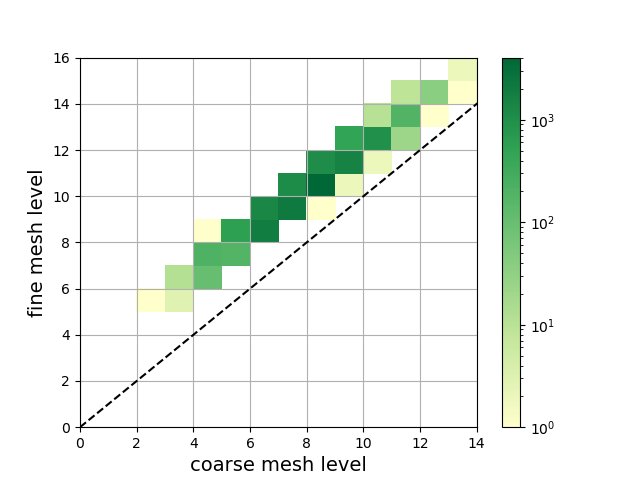}
	\caption{Example 2 with $\sigma^2=4$ using AMLMC: Stochastic stopping of coarse-fine pairs for {the} mesh level $\ell$. Top: $\ell=1,2$ from left to right. Middle: $\ell=3,4$ from left to right. Bottom: $\ell=5$.}
	\label{fig:ex2_sig2_coarse_fine_pair}
\end{figure}

Figure~\ref{fig:ex2_sig2_coarse_fine_pair} plots the histograms of coarse-fine mesh level of AMLMC in Example 2 with $\sigma^2 = 4$ for level $\ell = 1, 2, 3, 4, 5$. Similar to the previous example, the fraction of zero samples ($\Delta Q_\ell (\omega) = 0$) reduces with an increasing $\ell$, and disappears starting from $\ell = 3$. The samples satisfying a particular error estimate {are extensively} distributed into many auxiliary meshes, {indicating} that computing samples on a fixed mesh, as in SMLMC or {using} a batch adaptive approach, is less efficient.

\begin{table}[ht]
	\centering
	\begin{tabular}{lllllll}
		\hline
		Level, $\ell$ & 0 & 1 & 2 & 3 & 4 & 5 \\
		\hline
		SMLMC  & 2.48e+4 & 2.50e+3 & 1.66e+3  & 1.07e+3  & 9.04e+2  & 8.40e+2  \\
		AMLMC & 2.51e+4 & 1.67e+3 & 2.19e+3 & 5.67e+2 & 2.95e+2 &  1.19e+2 \\
		\hline
	\end{tabular}
	\caption{Example 2 with $\sigma^2=1$: The values $\sqrt{V_\ell W_\ell}$ for the considered MLMC estimators.}
  \label{tab:ex2_sig1_standard_adaptive_MLMC_work}
\end{table}

\begin{table}[ht]
	\centering
	\begin{tabular}{lllllll}
		\hline
		Level, $\ell$ & 0 & 1 & 2 & 3 & 4 & 5 \\
		\hline
		SMLMC & 3.15e+5 & 3.65e+4 & 1.54e+4  & 9.84e+3  & 8.08e+3  & 1.06e+4  \\
		AMLMC  & 3.38e+5 & 1.41e+4 & 6.09e+3 & 2.75e+3 & 1.13e+3 &  5.17e+2 \\
		\hline
	\end{tabular}
	\caption{Example 2 with $\sigma^2=4$: The values $\sqrt{V_\ell W_\ell}$ for the considered MLMC estimators.}
  \label{tab:ex2_sig2_standard_adaptive_MLMC_work}
\end{table}

\begin{table}[ht]
	\centering
	\begin{tabular}{llllll}
		\hline
		$\ell$ & 1 & 2 & 3 & 4 & 5\\
		\hline
		$\sigma^2 = 1$  & $0.31$ & $2.42$ & $0.77$  & $0.91$  & $0.56$\\
		$\sigma^2 = 4$ & $0.98$ & $1.54$ & $0.89$ & $0.75$ & $0.56$\\
		\hline
	\end{tabular}
	\caption{Example 2 for AMLMC: The expected ratio of the computed variance $\textrm{Var}(Q_{\ell} - Q_{\ell - 1})$ to the corresponding estimated variance derived from~\eqref{eq:optimal_variance_mesh_size_fam}. The ratios are of $\mathcal{O}(1)$, indicating a sharp enough error estimate. }
	\label{tab:ratio_variance_emperical}
\end{table}

In Table~\ref{tab:ex2_sig1_standard_adaptive_MLMC_work} and \ref{tab:ex2_sig2_standard_adaptive_MLMC_work}, the values $\sqrt{V_\ell W_\ell}$ are the largest in level $\ell = 0$. The estimated MLMC work $W = \tol^{-2} (\sum_{\ell=0}^{L} \sqrt{V_{\ell} W_\ell})^2$. In both cases, the quantity $\sqrt{V_\ell C_\ell}$ in AMLMC decays like a geometric sequence, while such quantity is asymptotically equally distributed in SMLMC, indicating the logarithmic dependence of $\tol$ in the SMLMC complexity. The logarithmic term in the complexity is revealed in Figures~\ref{fig:matern_field_sig_10e-1_cost_log} and \ref{fig:matern_field_sig_20e-1_cost_log}. This {confirms} the theoretical complexity prediction of $\tol^{-2}$ for  AMLMC and $\tol^{-2}\log^2(\tol^{-1})$ for SMLMC. 

{To confirm} that our asymptotic complexity predictions are accurate, cf. Corollary \ref{cor:AMLMCComplexity},  the computational 
 variances across levels {must} be close to their corresponding ones based on \eqref{eq:optimal_variance_mesh_size_fam}. Table~\ref{tab:ratio_variance_emperical} shows that for Example 2 in Section \ref{sec:example_2} the estimated variances are {sufficiently} sharp, even for coarse levels.
\section{Conclusions}
\label{sec:concl}
We proposed an AMLMC for goal-oriented approximation in a linear elliptic PDE with geometric singularities and a lognormal random diffusivity coefficient. 
We {built} our AMLMC on~\cite{adFEM_our}, which developed weak convergence rates for an adaptive algorithm {using} isoparametric d-linear quadrilateral finite element approximations and the dual weighted residual error representation in {a} deterministic setting. This theory provided us with sharp error estimates and indicators that guide the creation of a locally refined mesh sequences, tuned to the singularity at hand.
First, as a preparation phase, a sequence of deterministic $h$-adaptive auxiliary meshes {was} precomputed to satisfy a sequence of tolerances, {using} a deterministic approximation of the coefficient ${a}$, for instance, $\exp(\E{\log(a)})$. This step has a negligible cost and creates deterministic refinements that capture the geometry-driven singularities in the problem.

Then, each sample on a level in the MLMC hierarchy {selects} the coarsest auxiliary mesh such that a sample-dependent, \textit{a posteriori} error estimate-based stopping criterion is satisfied. The error contributions are {in general} not equally distributed among samples to achieve optimal approximation. Rather, they are subject to a scaling factor that may be unbounded in the lognormal case. Compared {with a} batch-adaptive MLMC, for instance in \cite{MLMC_AdFEM_Eigel_2016, oden_2018}, our method circumvents the cost of generating meshes {using} a batch of realizations, based on the assumption that the error density blows up in a fixed location, {as} in our case of geometrically driven singularities. Furthermore, using stochastic meshes, {our approach} is more efficient than a batch adaptive algorithm when there is a lack of uniform coercivity, {such as} in the lognormal case.
{Moreover}, our approach specifically avoids the costly overhead of building adaptive meshes from scratch for each realization. Instead, we only need to check a stopping criterion within the sequence of auxiliary meshes, which are already available.

Furthermore, {when} the solving cost is higher than the assembly cost, we discuss the use of iterative solvers and compare their efficiency with direct ones. To save computational work, we proposed a goal-oriented stopping criterion for the iterative solver, the realization of the diffusivity coefficient, and the desired level of AMLMC approximation.

Theoretically, the error density characterizes the cases where adaptivity provides a noticeable advantage, namely, those where the error density blows up in $L^1_P(\mathcal{D}\times \Omega)$, as 
we refine the mesh around the singularities, cf Theorem \ref{thm:MLMC_compl} and Corollary \ref{cor:AMLMCComplexity}.

We worked with a finite Fourier expansion for the log diffusivity coefficient in our {2D} numerical experiments.
{Note that} this expansion is not the only way to represent a field. Our adaptive method can also work with other representations, see, for instance, the work \cite{Croci2018}. 
In our numerical experiments, we observed efficiency improvements compared {with} both
 standard MC and MLMC for a singularity similar to that at the tip of a slit modeling a crack. 
We compared AMLMC with SMLMC for log fields with variances $\sigma^2 = 1$ and $\sigma^2 = 4,$ respectively. 
In these examples, the error density-based error estimates are sharp. They accurately predict both the bias error and the variances across levels, indicating that the asymptotic theory is relevant to analyzing the observed behavior of our AMLMC.
{Although the observed complexity of the SMLMC deteriorates with the increase in variance, our AMLMC complexity remains stable between these two variances.} {Consistent} with the theoretical results, the observed complexity for those 2D cases is $\tol^{-2}$ for AMLMC, $\tol^{-2} \log^2(\tol)$ for SMLMC, $\tol^{-3}$ for AMC, and $\tol^{-4}$ for MC. We predict further gains favoring AMLMC for larger values of dimension $d$. We note in passing that the uniform refinements fail to produce a sharp error estimate, even when coupled with the same error estimate that is sharp when using adaptive refinements. {This makes} the adaptive approach more attractive {because} it produces more robust and accurate results.
Furthermore, when achieving the $\tol^{-2}$ complexity, our AMLMC can be directly coupled with the popular unbiased MLMC approach introduced {in} \cite{McLeish2011,Rhee2015}, simply by randomizing over the sequence of tolerances used to create the adaptive levels.

The pointwise convergence of the error density proved in \cite{adFEM_our} is the crucial theoretical component of our approach.
The theoretical analysis of the error density convergence relies on proper local averaging and the multilinear structure of the isoparametric d-linear quadrilateral finite elements with hanging nodes. Hanging nodes are {important} to provide {sufficient} flexibility in the mesh refinements. The efficient computation of local averages in the error estimate is a contribution to this work, improving on \cite{adFEM_our}.
The pointwise convergence of the error density is {fundamental to demonstrating the} theoretical results on stopping, asymptotic accuracy, and efficiency, as proved in \cite{adFEM_our}, and {is} inherited by our AMLMC. 
{However}, the theoretical convergence of the error density also requires the solution and its dual to have $C^3(\mathcal{D})$ pathwise regularity.
As {reported} in \cite{adFEM_our}, this is not a practical limitation. A sequence of regularized geometries (for example by rounding a reentrant corner with a small radius) {can be produced} that satisfies the $C^3(\mathcal{D})$ assumptions and converges to our problem, showing up the $L^1(\mathcal{D})$ explosion of the error density, which governs the complexity of SMLMC, and the uniform boundedness of the quasi-norm, corresponding to the AMLMC. In practice, regularization {is not required} to produce accurate error estimates and a stable estimate of error density.
\section*{Acknowledgements}

This publication is based on work supported by the Alexander von Humboldt Foundation and the King Abdullah
University of Science and Technology (KAUST) office of sponsored research (OSR) under Award
No. OSR-2019-CRG8-4033. 
The authors thank Daniele Boffi and Alexander Litvinenko for fruitful discussions.
We also acknowledge the use of the following open-source software
packages: \texttt{deal.II}~\cite{dealII92}, \texttt{PAPI}~\cite{PAPI_paper}.

\appendix
\section{Algorithm Listings}
\label{sec:algos}

We list the algorithms and corresponding notations below. 

\begin{algorithm}
	\SetAlgoLined
	\KwIn{Iteration solver tolerance $\tol_{\mathrm{iter}}$, primal iteration operator $\mathscr{A}_p$,  dual iteration operator $\mathscr{A}_d$,  initial vector $u^{(0)}$ and $\varphi^{(0)}$, primal linear system $A_p$ and right-hand side $b_p$, dual linear system $A_d$ and right-hand side $b_d$. }
	\KwOut{Primal solution $u_{\textrm{iter}}$ and dual solution $\varphi_{\textrm{iter}}$}
	
	Set $\textnormal{PrimalSolver} = true$\\
	Set $\textnormal{DualSolver} = true$\\
	Set iteration step $k = 0$\\
	\While{$\textnormal{PrimalSolver}~||~\textnormal{DualSolver}$}
	{
		Proceed the primal iteration $u^{(k+1)} = \mathscr{A}_p u^{(k)}$\\
		Proceed the dual iteration $\varphi^{(k+1)} = \mathscr{A}_d \varphi^{(k)}$\\
		
		Compute primal residual $r_p^{(k+1)} = b_p - A_p u^{(k+1)}$\\
		Compute dual residual $r_d^{(k+1)} = b_d - A_d \varphi^{(k+1)}$\\
	
		\If{$\lvert (r_p^{(k+1)}, \varphi^{(k+1)} ) \rvert < \textnormal{TOL}_{\textrm{iter}} $}
		{
			Set $\textnormal{PrimalSolver} = false$
		}
		\If{$\lvert (r_d^{(k+1)}, u^{(k+1)} ) \rvert < \textnormal{TOL}_{\textrm{iter}}$}
		{
			Set $\textnormal{DualSolver} = false$
		}
		Set $k = k+1$\\
	}
	
	\caption{Primal and dual solvers iteration }
	\label{alg:primal_dual_solver}
\end{algorithm}

\begin{algorithm}
	\SetAlgoLined
	\KwIn{$\{\textnormal{TOL}_\ell\}_{\ell = 0}^{L}$, $\mathcal{H}_{-1}$, $L$, $\overline{\mathcal{N}}_{-1}$, $C_R$, $C_S$, $c$, $a$ }
	\KwOut{Mesh hierarchy $\{\mathcal{H}_\ell\}_{\ell=0}^L$}
	Initialize the mesh $\mathcal{H}_{-1}$ and the parameter $\mathcal{N}_{-1} = N_{-1}$, with $N_{-1}$ the number of cells on mesh $\mathcal{H}_{-1}$\\
	\For{$\ell = 0, 1,\dotsc,L$}
	{	
		$\text{keepRefining} = \text{TRUE}$\\
		$\mathcal{H}_\ell = \mathcal{H}_{\ell-1}$\\
		$\overline{\mathcal{N}}_\ell = c \overline{\mathcal{N}}_{\ell-1} $\\
		\While{\textnormal{keepRefining}}{
			Assemble the primal and dual linear system using the coefficient $a$ on mesh $\mathcal{H}_\ell$\\
			Solve the primal and dual problem~\eqref{eq:bvp_general} using Algorithm~\ref{alg:primal_dual_solver}, with $\tol_{\textrm{iter}} = \frac{1}{10} \tol_{\ell}$\\
			Compute approximate error density $\bar{\rho}$~\eqref{eq:approximate_error_density}, error indicator $r$ \eqref{eq:error_indicator} and the error estimate $e_{\textnormal{est}}$~\eqref{eq:error_estimate}\\
			\eIf{Stopping condition~\eqref{eq:stap_cond} is not satisfied}{
				\ForEach{$\textnormal{cell } k \in \mathcal{T}$}{
					\If{$r(k)$ meets the refinement criterion~\eqref{eq:ref_crit}}{
						Mark the cell $k$ for refinement\\
					}
				}
				Mark the cells for refinement to satisfy the mesh constraints\\
				Refine the marked cells on mesh $\mathcal{H}_\ell$\\
				Update the number of cells $N_\ell$ on mesh $\mathcal{H}_\ell$ and set $\overline{\mathcal{N}}_\ell = \max \{\overline{\mathcal{N}}_\ell, N_\ell\}$\\
				
			}
			{
				Set $\text{keepRefining} = \text{FALSE}$\\
			}
		}
	}
	
	\caption{Adaptive generation of an auxilliary mesh hierarchy}
	\label{alg:mesh_generation}
\end{algorithm}

The error indicator for cell $K$, using the approximated error density $\bar{\rho}$ \eqref{eq:approximate_error_density}, is given by,
\begin{align}
	\label{eq:error_indicator}
	r(K) &= \bar{\rho}_{K} h_K^{2+d}.
\end{align}

The error estimate $e_{\textnormal{est}}$ is the summation of all local error indicators,
\begin{align}
	\label{eq:error_estimate}
	e_{\textnormal{est}} &= \sum_K r(K).
\end{align}

The absolute error estimate $e_{\textnormal{est, abs}}$ sums up all unsigned local error indicators,
\begin{align}
	\label{eq:abs_error_estimate}
	e_{\textnormal{est, abs}} &= \sum_K \lvert r (K) \rvert.
\end{align}

The stopping condition in Algorithm~\ref{alg:mesh_generation} is to control the largest error indicator,
\begin{align}
	\max_{1\leq K \leq N_\ell}  r(K) &< C_S \frac{\textnormal{TOL}_{\ell}}{\mathcal{N}_\ell}.
	\label{eq:stap_cond}
\end{align}
The stopping condition controls the maximal error indicator for a fast decay~\cite{adFEM_our, adOverview_our}). {Note} that $N_\ell$ is the number of cells on mesh level $\ell$, whereas the variable $\mathcal{N}_\ell$ is introduced to control the refinement and stopping. 

The refinement criterion in Algorithm~\ref{alg:mesh_generation} for cell $K$, $K = 1,2,\dotsc,N_\ell$, is to refine the cells with error indicators greater or equal to a threshold,
\begin{align}
	r(K)  &\geq C_R \frac{\textnormal{TOL}_{\ell}}{\mathcal{N}_\ell}.
	\label{eq:ref_crit}
\end{align}

In Algorithm~\ref{alg:sample_DQ}, the computed stochastic scaling factor $K_k(\omega)$ on mesh $k$ is given by,
\begin{align}
	\label{eq:computed_scaling_factor}
	K_k(\omega) = \frac{\int_\mathcal{D} \bar{\rho}_k(x;\omega)^{\frac{d}{p+d}}\,dx}{\int_\mathcal{D} \E{\bar{\rho}_0^{\frac{d}{p+d}}}},
\end{align}
where in Section~\ref{sec:numex} the denominator was estimated on the coarsest mesh 0. {Generally,} we propose to use an MLMC estimate of $\int_\mathcal{D} \E{\bar{\rho}_0^{\frac{d}{p+d}}}$ using an iterative update procedure. We listed Algorithm~\ref{alg:sample_DQ} {as selecting} specific solvers but in general, the user can choose the solvers in another way to optimize efficiency.

\begin{algorithm}
	\SetAlgoLined
	\KwIn{$\ell$, $\tol_{\mathrm{bias},\ell}$, $\{(\mathcal{H}_k,a_k)\}_{k=0}^\infty$}
	\KwOut{One sample of $(Q_\ell,Q_{\ell-1})$}
	Set $k=0$\\
	Set error estimate $e_{\textnormal{est}}=+\infty$\\
	Set $\textnormal{CoarseFlag} = true$\\
	\While{$true$}
	{
		Generate the random field ${a}_k$ on mesh $\mathcal{H}_k$\\
		
		\eIf{$k == 0$}{
			Compute $u_{k}$ and $\varphi_{k}$ from~\eqref{eq:bvp_general} and a direct solver UMFPACK~\cite{umfpack_davis2004algorithm}\\
		}{
			Set the iterative solver tolerance $\tol_{\textrm{iter}} = \frac{1}{10} K_{k-1}(\omega) \cdot \tol_{\textrm{bias}, \ell} $\\
			Compute $u_{k}$ and $\varphi_{k}$ from~\eqref{eq:bvp_general} and Algorithm~\ref{alg:primal_dual_solver}, with the iterative solver tolerance $\textrm{TOL}_{\mathrm{iter}}$\\
		}
		
		Compute $e_{\textnormal{est}}$ using~\eqref{eq:approximate_error_density}, \eqref{eq:error_indicator} and \eqref{eq:error_estimate}\\
		Compute $K_k(\omega)$ using~\eqref{eq:computed_scaling_factor}\\
		
		\If{$\ell=0$}
		{ Set $Q_{\ell-1}=0$ \\
			Set $CoarseFlag = false$\\
		}
		\If{$e_{\textnormal{est}} < K_k(\omega) \cdot \tol_{\mathrm{bias},\ell - 1} \wedge CoarseFlag$}
		{
			
			Set $Q_{\ell - 1} = Q(u_k)$\\
			Set $CoarseFlag = false$\\
		}
		\If{$e_{\textnormal{est}} < K_k(\omega) \cdot \tol_{\mathrm{bias},\ell}$}
		{
			Set $Q_{\ell} = Q(u_k)$\\
			break;
		}
		Set $k = k + 1$\\
	}
	\caption{One sample of $\Delta Q_\ell = Q_\ell-Q_{\ell-1}$ in AMLMC}
	\label{alg:sample_DQ}
\end{algorithm}

\begin{algorithm}
	\SetAlgoLined
	\KwIn{Quadtree of cells $\mathscr{T}$}
	\KwOut{Dictionary $\mathscr{D}$ of vertices lines $\mathscr{L}$, where the vertices of each line have the same $y$ coordinate}
	
	\SetKwFunction{assembleline}{\textsc{Assemble}-\textsc{line}}
	\SetKwProg{Fn}{Function}{:}{end}
	\Fn{\assembleline{$\mathscr{V}$}}{
	\SetAlgoNoLine
	Create empty vector of cells $\mathscr{V}_{\mathrm{Left}}$ and $\mathscr{V}_{\mathrm{Right}}$\\
	\ForEach{$ \text{Cell} \in \mathscr{V} $}{
		\eIf{Cell $\textrm{has children}$}{
		$\mathscr{V}_{\mathrm{Left}}$ push back $Cell\rightarrow$ child(0)\\
		$\mathscr{V}_{\mathrm{Left}}$ push back $Cell\rightarrow$ child(1)\\
		$\mathscr{V}_{\mathrm{Right}}$ push back $Cell\rightarrow$ child(2)\\
		$\mathscr{V}_{\mathrm{Right}}$ push back $Cell\rightarrow$ child(3)\\
}{$\mathscr{V}_{\mathrm{Left}}$ push back $Cell$}
	}
	Create empty lines of vertices $\mathscr{L}_{\textrm{Bottom}}$ and $\mathscr{L}_{\textrm{Top}}$\\
	\eIf{All cells $\in \mathscr{V}$ have no children}{
		\eIf{All cells $\in \mathscr{V}$ have the same size}{\ForEach{$ \text{Cell} \in \mathscr{V} $}{$\mathscr{L}_{\textrm{Bottom}}$ push back $Cell\rightarrow$ vertex(0)\\
		$\mathscr{L}_{\textrm{Top}}$ push back $Cell\rightarrow$ vertex(2)}
		}{\ForEach{$ \text{Cell} \in \mathscr{V} $}{$\mathscr{L}_{\textrm{Bottom}}$ push back $Cell\rightarrow$ vertex(0)}}
	\eIf{$\mathscr{L}_{\textrm{Bottom}}$ or $\mathscr{L}_{\textrm{Top}}$ has $y$ coordinate already in dict. $\mathscr{D}$}{Merge with the existing line in $\mathscr{D}$}{Save $\mathscr{L}_{\textrm{Bottom}}$ and $\mathscr{L}_{\textrm{Top}}$ in the dictionary $\mathscr{D}$}
\KwRet{}}{
	\assembleline{$\mathscr{V_{\textrm{Left}}}$}\\
	\If{$\mathscr{V_{\textrm{Right}}}$ is non-empty}
	{\assembleline{$\mathscr{V_{\textrm{Right}}}$}}
	}
}
	Create vector $\mathscr{V}$, containing the tree root cell $\mathscr{T} \rightarrow root$\\
	\assembleline{$\mathscr{V}$}
	\caption{Assemble $y$-line structures from a quadtree}
	\label{alg:tree_lines}
\end{algorithm}
\section{Mesh function optimizations}
\label{app:opt}

This appendix contains the calculations resulting in 
the optimal stochastic mesh size functions in~\eqref{eq:opt_h_stoch}

\paragraph{Optimal stochastic, adaptive, mesh functions}

We {examine} the ideal stochastic mesh size function, $h(x;\omega)$, 
which minimizes the work model~\eqref{eq:nelem_model} subject to 
the constraint that the bias model~\eqref{eq:bias_model} equals 
$\tol_\mathrm{bias}$. {For this purpose,} we introduce the Lagrangian,
\begin{align}
	\label{eq:Lagrangian1}
	\mathcal{L}(h;\lambda) & = \int_\mathcal{D} \E{h^{-d}} 
	+\lambda\left(\int_\mathcal{D} \E{\rho h^p}-\tol_\mathrm{bias}\right).
\end{align}
By {considering} the variation of the Lagrangian with respect to $h$, we obtain the optimality
condition
\begin{align*}
	-\frac{d}{h^{d+1}} + \lambda\rho p h^{p-1}  \equiv 0 && \Leftrightarrow &&
	h(x;\omega) = \left(\frac{1}{\lambda}\frac{d}{p}\rho(x;\omega)^{-1}\right)^{\frac{1}{p+d}}.
\end{align*}
Then, using the bias constraint, {we have}
\begin{align*}
	\tol_\mathrm{bias} & = 
	\int_\mathcal{D} \E{\rho\left(\frac{1}{\lambda}\frac{d}{p}\rho^{-1}\right)^{\frac{p}{p+d}}},
\end{align*}
we can determine the Lagrange multiplier,
\begin{align*}
	\lambda^{-\frac{1}{p+d}} & =
	\left(\frac{p}{d}\right)^{\frac{1}{p+d}}
	\frac{\tol_\mathrm{bias}^{1/p}}{\left(\int_\mathcal{D} \E{\rho^{\frac{d}{p+d}}}\right)^{1/p}},
\end{align*}
which after substituting back into the expression for $h(x;\omega)$
gives 
\begin{align*}
	h(x;\omega) & = 
	\frac{\tol_\mathrm{bias}^{1/p}}{\left(\int_\mathcal{D} \E{\rho^{\frac{d}{p+d}}}\right)^{1/p}}
	\rho(x;\omega)^{-\frac{1}{p+d}},
\end{align*}
as stated in~\eqref{eq:opt_h_stoch}.

\paragraph{Optimal stochastic, uniform, mesh functions}

Restricting ourselves to constant mesh size functions, 
$h(x;\omega)\equiv h(\omega)$, corresponding to idealized
uniform meshes, the Lagrangian in~\eqref{eq:Lagrangian1}
becomes
\begin{align*}
	\mathcal{L}(h;\lambda) & = \E{h^{-d}\int_\mathcal{D} 1\,dx} 
	+\lambda\left(\E{ h^p\int_\mathcal{D} \rho(x;)\,dx}-\tol_\mathrm{bias}\right).
\end{align*}
Introducing the constant $A=\int_\mathcal{D} 1\,dx$, the area of $D$, and the 
random variable $R=\int_\mathcal{D} \rho(x;\cdot)\,dx$, the $L^1(D)$-norm of 
$\rho(\cdot;\omega)$, and again taking variations with respect to $h$, the
optimality condition is {shown} to be
\begin{align*}
	-\frac{d}{h^{d+1}}A + \lambda R p h^{p-1}  \equiv 0 && \Leftrightarrow &&
	h(\omega) = \left(\frac{1}{\lambda}\frac{d}{p}\frac{A}{R(\omega)}\right)^{\frac{1}{p+d}}.
\end{align*}
The bias constraint becomes
\begin{align*}
	\tol_\mathrm{bias} & = 
	\E{\left(\frac{1}{\lambda}\frac{d}{p}\frac{A}{R}\right)^{\frac{p}{p+d}}R},
\end{align*}
from which we conclude
\begin{align*}
	\lambda^{-\frac{1}{p+d}} & =
	\left(\frac{p}{d}\frac{1}{A}\right)^{\frac{1}{p+d}}
	\frac{\tol_\mathrm{bias}^{1/p}}{\E{R^{\frac{d}{p+d}}}^{1/p}}.
\end{align*}
Substitution back {in} $h(\omega)$ yields,
\begin{align*}
	h(\omega) & = 
	\tol_\mathrm{bias}^{1/p}\frac{R(\omega)^{-\frac{1}{p+d}}}{\E{R^{\frac{d}{p+d}}}^{1/p}},
\end{align*}
which using the definitions of $A$ and $R$ becomes~\eqref{eq:opt_h_uniform}.
\section{Proof of Corollary~\ref{cor:AMLMCComplexity}}
\label{app:proof_corollary}

{

Associated with the work model~\eqref{eq:work_sum} is the upper 
bound~\eqref{eq:total_work_MLMC} where, since Assumption~\ref{assu:MLMC} is
satisfied with $\kappa=\tolbsep^{-1}$, $\gamma=d/p$, $q_s=2$, and $q_w=1$, 
the second term is asymptotically negligible as $\tol\to 0$.
For the first term we obtain, 
assuming that $V_0$ is constant and $W_0=\E{K_2}\tol_0^{-d/p}$
while using~\eqref{eq:optimal_variance_mesh_size_fam} 
and~\eqref{eq:optimal_work_mesh_size_fam} for $\ell\geq 1$,
in $\sum_{\ell=0}^L \sqrt{V_\ell W_\ell}$, that
\begin{multline}
	\label{eq:term_1}
	\left( \frac{\confpar}{ \splitting\tol} \right)^2
  \left( \sum_{\ell=0}^L \sqrt{\work_\ell V_\ell} \right)^2
  \\ = 
  \frac{K_3 K_5}{\tol^2}
  \underbrace{
  \tol_0^{-d/p}
  \left(
  	\sqrt{\frac{V_0}{\var{K_1}}}
		\frac{1}{\left(\tolbsep^{-1}-1\right) \sqrt{1+\tolbsep^{d/p}}}
		+
		\tol_0 S
  \right)^2}_{T}, 
\end{multline}
where
\begin{align*}
	K_3 & = \E{K_2}\var{K_1}, \\ 
	K_5 & = 
		\left( \frac{\confpar}{ \splitting} \right)^2
		\left(\tolbsep^{-1}-1\right)^2 \left(1+\tolbsep^{d/p}\right),\\
	S & = \sum_{\ell=1}^L \left(\tolbsep^{1-\frac{d}{2p}}\right)^\ell
		= 
		\begin{cases}
			L & \text{, if $d=2p$},\\
			\tolbsep^{1-\frac{d}{2p}}
			\frac{1-\left(\tolbsep^{1-\frac{d}{2p}}\right)^L}
			{1-\tolbsep^{1-\frac{d}{2p}}} 
			&  \text{, if $d\neq 2p$.}
		\end{cases}
\end{align*}
To satisfy the bias constraint, 
$\left(1-\splitting\right)\tol\geq\tol_L=\tolbsep^L\tol_0$, 
choose
\begin{align*}
	L & = \left\lceil L^\ast \right\rceil, &
	L^\ast & = \log_C{\left(\frac{1-\splitting}{\tol_0}\tol\right)}, 
\end{align*}
provided $\tol_0>\left(1-\splitting\right)\tol$. 
Then $L^\ast\in\rset^+$, since $\tolbsep\in(0,1)$.

Keeping $\tol_0$ fixed while letting $\tol\to 0$, one has three different
situations depending on the asymptotic behavior of $S$: 

If $d<2p$, then $\left(\tolbsep^{1-\frac{d}{2p}}\right)^L\to 0$
as $\tol\to 0$, so that $T$ in~\eqref{eq:term_1} remains bounded and
\begin{align*}
	\frac{\work_\mathrm{MLMC}}{\tol^{-2}} & \to
	\frac{K_3 K_5 }
	{\tol_0^{d/p}}
  \left(
  	\sqrt{\frac{V_0}{\var{K_1}}}
  	\frac{1}{\left(\tolbsep^{-1}-1\right) \sqrt{1+\tolbsep^{d/p}}}
  	+
  	\tol_0\frac{C^{1-\frac{d}{2p}}}{1-C^{1-\frac{d}{2p}}}
  \right)^2,
\end{align*}
as $\tol\to 0$.

If $d=2p$, then
\begin{align*}
	\frac{L}{\log{\left(\tol^{-1}\right)}}
	& = 
	\frac{1}{\log{\left(\tol^{-1}\right)}}
	\left\lceil
		\frac{\log{\left(\tol^{-1}\right)}+\log{\left(\frac{\tol_0}{1-\splitting}\right)}}
		{\log{\left(\tolbsep^{-1}\right)}}
	\right\rceil
	\to \frac{1}{\log{\left(\tolbsep^{-1}\right)}},
\end{align*}
and
\begin{align*}
	\frac{\work_\mathrm{MLMC}}{\tol^{-2}\log{\left(\tol^{-1}\right)}^2} & \to
	K_3 K_5 \tol_0^{-d/p}
  \log{\left(\tolbsep\right)}^{-2},
\end{align*} 
as $\tol\to 0$.

If $d>2p$, then 
\begin{align*}
	\tolbsep^{\left(1-\frac{d}{2p}\right)L^\ast} & = 
	\left(\frac{1-\splitting}{\tol_0}\tol\right)^{\left(1-\frac{d}{2p}\right)}.
\end{align*}
Thus, allowing for the non-integer approximation $L^\ast$, 
\begin{align*}
	\frac{S}{\tol^{1-\frac{d}{2p}}} & \to 
	\left(\frac{1-\splitting}{\tol_0}\right)^{\left(1-\frac{d}{2p}\right)}
	\frac{1}{1-\tolbsep^{\frac{d}{2p}-1}},\\
	\frac{T}{\tol^{2-d/p}} & \to 
	\left(1-\splitting\right)^{2-d/p}
	\frac{1}{\left(1-\tolbsep^{\frac{d}{2p}-1}\right)^2},
\end{align*}
and
\begin{align*}
	\frac{\work_\mathrm{MLMC}}{\tol^{-d/p}} & \to
	K_3 K_5 \frac{\left(1-\splitting\right)^{2-d/p}}{\left(1-\tolbsep^{\frac{d}{2p}-1}\right)^2},
\end{align*}
as $\tol\to 0$.

Together, the previous three cases show the limits of Corollary~\ref{cor:AMLMCComplexity}, 
allowing for non-integer $L$ when $d>2p$. Taking the integer constraint into 
account, it follows from $L^\ast\leq L<L^\ast+1$ that
\begin{align*}
	\liminf_{\tol\to 0}\frac{\work_\mathrm{MLMC}}{\tol^{-d/p}} & \geq
	K_3 K_5 \frac{\left(1-\splitting\right)^{2-d/p}}{\left(1-\tolbsep^{\frac{d}{2p}-1}\right)^2},\\
	\limsup_{\tol\to 0}\frac{\work_\mathrm{MLMC}}{\tol^{-d/p}} & \leq
	K_3 K_5\frac{\left(1-\splitting\right)^{2-d/p}}{\left(1-\tolbsep^{\frac{d}{2p}-1}\right)^2} \tolbsep^{1-\frac{d}{2p}}.
\end{align*}
}

\bibliographystyle{siam}
\bibliography{mlmc_afem_bibliography}

\end{document}